\documentclass[12pt]{amsart}
\usepackage[top=72pt,bottom=96pt,left=72pt,right=72pt]{geometry}
\usepackage{amsmath}
\usepackage{mathtools}
\usepackage{amssymb}
\usepackage{bbm}
\usepackage{extarrows}
\def\qA{\mathrm{qAss}}
\def\G{\mathcal{G}}
\def\C{\mathbb{C}}

\def\N{\mathbb{N}}

\def\F{\mathfrak{F}}
\def\f{\mathfrak{f}}
\def\qGG{\mathfrak{qGG}}
\def\Hor{\mathrm{Hor}}
\def\Ad{\mathrm{Ad}}
\def\ad{\mathrm{ad}}
\def\id{\mathrm{id}}
\def\Id{\mathrm{Id}}
\def\Ker{\mathrm{Ker}}
\def\Im{\mathrm{Im}}
\def\inv{\mathrm{inv}}
\def\alt{\mathrm{alt}}
\def\c{\mathrm{c}}

\def\Mor{\textsc{Mor}}
\def\N{\mathbb{N}}

\def\FD{\textsc{FD}}
\def\Z{\mathbb{Z}}
\def\Rep{\mathbf{Rep}}
\def\triv{\mathrm{triv}}
\def\qtrs{\mathrm{qtrs}}
\def\l{\mathrm{L}}
\def\c{\mathrm{c}}
\def\End{\mathrm{End}}
\def\r{\mathrm{R}}

\def\G{\mathcal{G}}

\def\T{\mathcal{T}}

\newtheorem{Lemma}{Lemma}[section]
\newtheorem{Remark}[Lemma]{Remark}
\newtheorem{Example}[Lemma]{Example}
\newtheorem{Theorem}[Lemma]{Theorem}
\newtheorem{Corollary}[Lemma]{Corollary}
\newtheorem{Definition}[Lemma]{Definition}
\newtheorem{Proposition}[Lemma]{Proposition}

\begin{document}
\date{\today}
\title{Geometry of Associated Quantum Vector Bundles and the Quantum Gauge Group}
\author{Gustavo Amilcar Salda\~na Moncada}
\address{Gustavo Amilcar Salda\~na Moncada\\
Mathematics Research Center, CIMAT}
\email{gamilcar@ciencias.unam.mx}
\begin{abstract}
It is well--known that, given a principal $G$--bundle equipped with a principal connection, one can associate to every unitary finite--dimensional representation of $G$ a linear connection and a compatible Hermitian structure on the corresponding associated vector bundle. Moreover, the gauge group acts both on the space of principal connections and on the space of induced linear connections defined on the associated vector bundles. The aim of this paper is to present the \emph{non--commutative geometrical} counterpart of these \emph{classical} facts within the framework of quantum principal bundles and quantum principal connections.

 \begin{center}
  \parbox{300pt}{\textit{MSC 2010:}\ 46L87, 58B99.}
  \\[5pt]
  \parbox{300pt}{\textit{Keywords:}\ Quantum Connections, Hermitian Structures, Quantum Gauge Group.}
 \end{center}
\end{abstract}
\maketitle
\section{Introduction}

Non--commutative geometry, also known as {\it quantum geometry}, arises as a kind of algebraic generalization of geometrical concepts \cite{con, prug, woro0}.  There are a variety of reasons to believe that this branch of mathematics may be able to solve some of the Standard Model's fundamental problems \cite{con}.

Pursuing this philosophy, in \cite{micho1, micho2, micho3} M. Durdevich developed a formulation of the theory of principal bundles and principal connections in the non--commutative geometry framework\footnote{Called quantum principal bundles (qpb's) and quantum principal connections (qpc's).}. This theory uses the concept of a quantum group, as presented by S. L. Woronowicz in \cite{woro1, woro2}, which plays the role of the structure group on the bundle. However, it uses a more general differential calculus on the quantum group that allows one to extend the complete $\ast$--Hopf algebra structure, reflecting the {\it classical} fact that the tangent bundle of every Lie group is a Lie group as well \cite{nodg}. Furthermore, Durdevich's formulation embraces other {\it classical} concepts, such as characteristic classes and classifying spaces \cite{micho6, micho7}.

The paper \cite{sald} develops a categorical equivalence between principal bundles with principal connections over a fixed base space $M$ and the category of associated functors called {\it gauge theory sectors}; Durdevich's theory allows one to recreate this result for quantum principal bundles and real, regular quantum principal connections \cite{sald1}. In concrete, there is a categorical equivalence between quantum principal bundles with real, regular quantum principal connections over a fixed quantum base space $B$ and the category of contravariant functors between the category of finite--dimensional corepresentations and the category of quantum vector bundles with quantum linear connections.  This provides a clear motivation to further develop the theory. The question, what other interesting results from differential geometry hold in the non--commutative setting? was the starting point of this paper.

In this way, the purpose of this work is to extend the theory of associated quantum vector bundle and induced quantum linear connections\footnote{Associated qvb's and induced qlc's} presented in \cite{sald1}, following the line of research of M. Durdevich \cite{micho7} and works by other authors \cite{qvbH, lrz, z}, to add canonical Hermitian structures as well as to study the relationship between this new structure and the induced qlc's for any real qpc. Moreover, we will introduce and analyze an {\it ad hoc} definition of the quantum gauge group for a given qpb with a differential calculus and we will study its natural action on the space of qpc's and the space of qlc's.

We believe that the approach presented is important not only because of the results that we will prove, which reflect the analogy with the {\it classical} case and extend the theory (the reader should pay particular attention to equations (\ref{3.f8}), (\ref{3.f8.1}), (\ref{3.f10.3}), (\ref{3.f10.6}); Definitions \ref{lhs}, \ref{rhs}, \ref{qgg}; Proposition \ref{4.4} and Theorems \ref{hil}, \ref{hil1}, \ref{fgs}, \ref{4.3}), but also because our approach opens the door to many other research lines, such as the moduli space of quantum connections and the Yang--Mills models and field theory, in accordance with the work presented in \cite{lrz, z}. Applying our methods to the study of the Yang--Mills models and field theory (see \cite{sald2, sald3, sald4}) is the ultimate goal of our research.

The paper is organized into six sections. Following this introduction, in the second section, we present preparatory material, broken down into two subsections. The first subsection is about the compact matrix quantum group, its corepresentations, and the universal differential envelope $\ast$--calculus of a first--order differential $\ast$--calculus. In the second subsection, we will present all the basic notions of qpb's and qpc's; however, we will change the standard definition of qpc's in order to embrace a more general theory. The third section also consists of two subsections. In the first subsection, we will develop the general theory of associated qvb's and induced qlc's. In the second one, we will introduce a canonical Hermitian structure which is compatible with induced qlc's. The fourth section is about the quantum gauge group and its action on quantum connections. In the fifth section we present two different classes of examples of our theory: trivial qpb's (in the sense of \cite{micho2}) and homogeneous qpb's. The final section contains some concluding comments.

We shall follow the notation introduced in \cite{sald1}. For instance, every (compact) quantum space will be identified with its $\ast$--algebra of smooth $\C$--valued functions. In other words, 
all  quantum spaces will be formally given by associative 
unital $\ast$--algebras over $\C$ $$(X,\cdot,\mathbbm{1},\ast).$$ In general, we are going to omit the words \emph{associative} and \emph{unital}. 
Furthermore, all of our $\ast$--algebra morphisms will be unital, and 
throughout the paper we use Sweedler notation: for the coproduct 
$\Delta$ of a quantum group, we write $$\Delta(g) = g^{(1)} \otimes 
g^{(2)},$$  and for a (right) corepresentation, such as $\delta^V$, $\Ad$, $\ad$, 
$\Delta_P$, $\Delta_{\Omega^\bullet(P)}$ or $\Delta_{\Hor}$, we denote the image of $a$ by 
$$a^{(0)} \otimes a^{(1)}.$$

In the literature, there are other viewpoints on qpb’s (see 
\cite{bm, bu, pl}), all of which are intrinsically related by the theory 
of Hopf--Galois extensions \cite{kt}. We have chosen to employ 
Durdevich’s formulation of qpb’s because of its purely 
geometric--algebraic framework, where differential calculus, 
connections, their curvature, and their covariant derivatives are the 
most relevant objects.

It is worth noting that, throughout this paper, we will use the word \emph{classical} or the expression \emph{classical case} to refer to differential geometry. In this way, we will explain how our definitions and constructions reflect the \emph{classical case}, allowing the reader to appreciate the \emph{naturalness} of our results by comparing them with their counterparts in differential geometry. Throughout the text, we use the terms \emph{non--commutative} and \emph{quantum} interchangeably. Although this work is based on Durdevich's theory, our definition of the quantum gauge group follows that introduced in \cite{br}, but formulated at the level of differential calculus. Unfortunately, not all results valid at degree~$0$ can be extended to the level of differential calculus. A concrete example illustrating this limitation will be presented in the following section.

\section{Preparatory Material}

As in the \emph{classical} case, the notion of a group plays a fundamental 
role in the theory of principal bundles. Therefore, we review the key aspects 
of the framework developed by S. L. Woronowicz \cite{woro1, woro2}, as well 
as the universal differential envelope $\ast$--calculus presented in 
\cite{micho1, stheve}.

\subsection{Quantum Groups}

A compact matrix quantum group (for this paper, we shall refer to it simply as a quantum group) will be denoted by $\G$, and its dense $\ast$--Hopf (sub)algebra will be denoted by
\begin{equation}
    \label{2.f0}
    H^\infty:=(H,\cdot,\mathbbm{1},\Delta,\epsilon,S,\ast),
\end{equation}
where  $\Delta$ is the coproduct, $\epsilon$ is the counity and $S$ is the coinverse. The space $H^\infty$ shall be treated as the algebra of all {\it polynomial functions} defined on $\G$.  In the same way, a (smooth right) $\G$--corepresentation on a $\C$--vector space $V$ is a linear map $$\delta^V: V \longrightarrow V \otimes H$$ such that
\begin{equation} 
\label{2.f1}
(\id_V\otimes \epsilon )\circ \delta^V= \id_V \qquad \mbox{ and } \qquad (\id_V\otimes \Delta )\circ \delta^V=(\delta^V\otimes \id_H ) \circ \delta^V.
\end{equation}

\noindent We say that the corepresentation is {\it finite--dimensional} if $\dim_{\C}(V)< |\N|$. The map $\delta^V$ is often referred to as \emph{(right) coaction} of $\G$ on $V$. It is worth mentioning that in the general theory (\cite{woro1}), the first part of equation (\ref{2.f1}) is not necessary. 

Given two $\G$--corepresentations $\delta^V$, $\delta^W$,  a {\it corepresentation morphism} is a linear map 
\begin{equation}
\label{2.f3}
T:V\longrightarrow W \qquad \mbox{ such that }\qquad (T\otimes \id_H)\circ \delta^V=\delta^W \circ T.
\end{equation}
The notions of monomorphism, epimorphism and isomorphism of corepresentations should be clear. The set of all corepresentation morphisms between two corepresentations $\delta^V$, $\delta^W$ will be denoted as 
\begin{equation}
\label{2.f4}
\Mor(\delta^V,\delta^W),
\end{equation}
and the set of all finite--dimensional $\G$--corepresentations will be denoted by
\begin{equation}
\label{2.f6}
\FD(\Rep_{\G}).
\end{equation}

A $\G$--corepresentation $\delta^V$ is {\it reducible} if there exists a non--trivial subspace $L$ $(L\not=\{0\}, V)$ such that $\delta^V(L)\subseteq L\otimes H$ and $\delta^V$ is {\it unitary} if viewed as an element of $B(V)\otimes H$ (with $B(V):=\{f:V\longrightarrow V\mid f \mbox{ is linear} \}$) is unitary. Of course, for the last definition, it is necessary an inner product $\langle-| -\rangle$ on $V$ and a corepresentation is said to be {\it irreducible} if it is not reducible. In \cite{woro1}, Woronowicz proved that every finite--dimensional $\G$--corepresentation on $V$ admits an inner product $\langle -|-\rangle$ (not necessarily unique) that makes the corepresentation unitary. Henceforth, we will assume that every finite--dimensional $\G$--corepresentation is unitary.

A proof of the following theorem can be found in \cite{woro1}.

\begin{Theorem}
\label{rep}
Let $\T$ be a complete set of mutually non--equivalent, irreducible, finite--dimensional $\G$--corepresentations with $\delta^C_\triv$ $\in$ $\T$, where $\delta^C_\triv$ denotes the trivial corepresentation of $H$ on $\C$ (\cite{woro1}). For any $\delta^V$ $\in$ $\T$ that coacts on $V$,
\begin{equation}
\label{2.f8}
\delta^V(e_j)=\sum^{n_{V}}_{i=1} e_i\otimes g^{V}_{ij},
\end{equation}
where $\{ e_i\}^{n_{V}}_{i=1}$ is an orthonormal basis of $V$ (with respect to the inner product that makes $\delta^V$ unitary) and $\{g^{V}_{ij}\}^{n_{V}}_{i,j=1}$ $\subseteq$ $H$. Then $\{g^{V}_{ij}\}_{\delta^V,i,j}$ is a linear basis of $H$, where the index $\delta^V$ runs on $\T$ and $i$, $j$ run from $1$ to $n_{V}=\mathrm{dim}_\C(V)$.
\end{Theorem}

For every $\delta^V$ $\in$ $\T$, the set $\{g^{V}_{ij}\}^{n_{V}}_{ij=1}$ satisfies
\begin{equation}
\label{2.f9}
\begin{aligned}
    \Delta(g^{V}_{ij})= \sum^{n_{V}}_{k=1} &g^{V}_{ik}\otimes g^{V}_{kj},\qquad S(g^{V}_{ij})=g^{V\,\ast}_{ji},\qquad \epsilon(g^{V}_{ij})=\delta_{ij}
    \\ &\sum^{n_{V}}_{k=1}S(g^{V}_{ik})\,g^{V}_{kj}=\sum^{n_{V}}_{k=1}g^{V}_{ik}\,S(g^{V}_{kj})=\delta_{ij}\mathbbm{1}
\end{aligned}
\end{equation}
with $\delta_{ij}$ being the Kronecker delta, among other properties \cite{woro1}.

We now present a brief overview of first--order differential $\ast$--calculus 
($\ast$--FODC). For additional details, see \cite{woro2, stheve}. 
A FODC over $\G$ or over $H$, is a pair $(\Gamma,d)$, where $\Gamma$ 
is an $H$--bimodule and $d: H \longrightarrow \Gamma$ is a linear map, referred to as the \emph{differential}, satisfying the following conditions:
\begin{enumerate}
    \item The Leibniz rule.
    \item For every $\vartheta \in \Gamma$, there exist (not necessarily unique) 
          elements $g_k, h_k \in H$ such that $\vartheta=\displaystyle \sum_{k}g_k(dh_k)$.      
\end{enumerate}   
If there exists an antilinear involution $$\ast : \Gamma \longrightarrow \Gamma$$ such that $(\vartheta g)^\ast=g^\ast \vartheta^\ast$, $(g\vartheta)^\ast=\vartheta^\ast g^\ast$ and $(dg)^\ast=d(g^\ast)$ for all $\vartheta$ $\in$ $\Gamma$, $g$ $\in$ $H$, we say that the FODC is actually a $\ast$--FODC \cite{stheve}. It can be proven that if such an involution exists, then it is unique, and in general, one does not need a quantum group or a $\ast$--Hopf algebra to define a $\ast$--FODC; a $\ast$--algebra is sufficient \cite{stheve}.

A $\ast$--FODC $(\Gamma,d)$ over $\G$ is left covariant if there exists a linear map 
\begin{equation}
    \label{2.f9.2}
    \Phi_\Gamma:\Gamma \longrightarrow H\otimes \Gamma
\end{equation}
such that
\begin{enumerate}
    \item  $\Phi_\Gamma$ preserves the $\ast$--structure and $\Phi_\Gamma(g\, \vartheta)=\Delta(g)\Phi_\Gamma(\vartheta)$ for all $\vartheta$ $\in$ $\Gamma$ and all $g$ $\in$ $H$. Here, the $\ast$--structure of $H\otimes \Gamma$ is given by $(g\otimes \vartheta)^\ast:=g^\ast\otimes\vartheta^\ast $.
    \item $\Phi_\Gamma$ satisfies $(\epsilon\otimes \id_{\Gamma})\circ \Phi_\Gamma = \id_{\Gamma}$ and $(\Delta\otimes \id_{\Gamma})\circ \Phi_\Gamma=(\id_H\otimes \Phi_\Gamma)\circ \Phi_\Gamma$.
    \item $\Phi_\Gamma\circ d=(\id_H\otimes d)\circ \Delta$.
\end{enumerate}
Similarly, a $\ast$--FODC $(\Gamma,d)$ over $\mathcal{G}$ is said to be right covariant if there exists a linear map
\begin{equation}
    \label{2.f9.1}
    {}_{\Gamma}\Phi:\Gamma \longrightarrow \Gamma\otimes H
\end{equation}
satisfying properties analogous to those of $\Phi_{\Gamma}$ (see Section 3 of \cite{stheve}). Finally, we say that a $\ast$--FODC $(\Gamma, d)$ is \emph{bicovariant} if it is both left covariant and right covariant.

The reader can find the exposition of the following example in \cite{libro}. It consists solely of straightforward calculations, relying only on basic properties of $\ast$–Hopf algebras such as coassociativity and the identities $m \circ (S \otimes \id_H) \circ \Delta = \mathbbm{1}\, \epsilon$, $m \circ (\id_H \otimes S) \circ \Delta = \mathbbm{1} \, \epsilon,$ where $m: H\otimes H\longrightarrow H$ is the product map.
\begin{Example}
    \label{innese1}
    Let $\G$ be a quantum group and consider the $\C$--vector space $$\Gamma_U:=H\otimes \Ker(\epsilon)$$ with the $\ast$--$H$--bimodule structure given by $$g\,(a\otimes b):=ga\otimes b,\qquad (a\otimes b)\,g:= ag^{(1)}\otimes bg^{(2)},\qquad (a\otimes b)^\ast=-a^{(1)\ast}\otimes S(b)^\ast a^{(2)\ast},$$ for $g$ $\in$ $H$ and $a\otimes b$ $\in$ $\Gamma_U$,  where $\Delta(g)=g^{(1)}\otimes g^{(2)}$. With this structure, the pair
    \begin{equation}
    \label{2.f9.3}
    (\Gamma_U, D),
\end{equation}
where  
\begin{equation}
    \label{2.f9.3.1}
   D: H \longrightarrow \Gamma_U,\qquad g\longmapsto \Delta(g)-g\otimes \mathbbm{1},
\end{equation}
is a $\ast$--FODC. Furthermore, by considering the linear maps
\begin{equation}
    \label{2.f9.3.2}
{_{\Gamma_U}}\Phi:\Gamma_U\longrightarrow \Gamma_U\otimes H,\qquad \Phi_{_{\Gamma_U}}:\Gamma_U\longrightarrow  H \otimes \Gamma_U
\end{equation}
given by $${_{\Gamma_U}}\Phi(a\otimes b)=a^{(1)}\otimes b^{(2)}\otimes a^{(2)}S(b^{(1)})b^{(3)} \quad \mbox{ and }\quad \Phi_{\Gamma_U}(a\otimes b)=a^{(1)} \otimes a^{(2)}\otimes b,$$  we obtain that $(\Gamma_U,D)$ is bicovariant.
\end{Example}

The bicovariant $\ast$--FODC $(\Gamma_U, D)$ is called {\it the universal  $\ast$--FODC} because of the following proposition. The reader can check a proof of it in \cite{woro2,libro,stheve}.

\begin{Proposition}
    \label{0U}
    Let  $\mathcal{R}$ $\subseteq$ $\Ker(\epsilon)$ be a right $H$--ideal such that $S(\mathcal{R})^\ast\subseteq \mathcal{R}$ and $\Ad(\mathcal{R})\subseteq \mathcal{R}\otimes H$, where
\begin{equation}
    \label{2.f9.5}
    \Ad:H\longrightarrow H\otimes H, \qquad g\longmapsto g^{(2)}\otimes S(g^{(1)})\,g^{(3)}
\end{equation}
is the right adjoint coaction. Then 
    \begin{equation}
    \label{2.f9.4}
    (\Gamma:=H\otimes {\Ker(\epsilon)\over \mathcal{R}},d_{\mathcal{R}})
\end{equation}
defines a $\ast$–FODC over $\G$, where the $H$–bimodule structure on $\Gamma$ is the quotient bimodule structure inherited from $\Gamma_U$, and $d_\mathcal{R}$ is the map defined by $D$ and the canonical projection from $\Gamma_U$ to $\Gamma$. Moreover, this $\ast$--FODC is bicovariant, where the maps ${_{\Gamma}}\Phi$, $\Phi_{{\Gamma}}$ of the bicovariant structure are the ones induced by ${_{\Gamma_U}}\Phi$ and $\Phi_{_{\Gamma_U}}$ on the quotient space $\Gamma$, respectively.

Reciprocally, every bicovariant $\ast$--FODC $(\Gamma,d)$ over $\G$ is isomorphic to the one in equation (\ref{2.f9.4}) for some right $H$--ideal $\mathcal{R}$ $\subseteq$ $\Ker(\epsilon)$ such that $S(\mathcal{R})^\ast\subseteq \mathcal{R}$ and $\Ad(\mathcal{R})=\mathcal{R}\otimes H$.
\end{Proposition}

It is worth mentioning that there are similar propositions for left covariant $\ast$--FODC's and right covariant $\ast$--FODC's \cite{stheve}. 

Let $(\Gamma,d)$ be a bicovariant $\ast$--FODC over $\G$ and consider the $\C$--vector space given by
\begin{equation}
\label{2.f13}
\mathfrak{qg}^\#:=\{\theta \in \Gamma \mid \Phi_\Gamma(\theta)=\mathbbm{1}\otimes \theta  \}   =\mathbbm{1}\otimes {\Ker(\epsilon)\over \mathcal{R}}\cong {\Ker(\epsilon)\over \mathcal{R}}.
\end{equation}
This space allows to consider the quantum germs map (see Section 6.4 of \cite{stheve})
\begin{equation}
\label{2.f14}
\pi:H \longrightarrow \mathfrak{qg}^\#,\qquad 
g \longmapsto S(g^{1})dg^{2}.
\end{equation}
 The map $\pi$ has several useful properties, for example, the restriction map $\pi|_{\Ker(\epsilon)}$ is surjective and 
\begin{equation}
\label{properties}
    \begin{aligned}
\ker(\pi)=\mathcal{R}\oplus \C\mathbbm{1}, \quad & \qquad dg=g^{(1)}\pi(g^{(2)}),\quad  \qquad   \pi(g)^\ast=-\pi(S(g)^\ast)\\
\pi(g)=-(d&S(g^{(1)}))g^{(2)},\qquad dS(g)=-\pi(g^{(1)})S(g^{(2)})&
    \end{aligned}
\end{equation}
for all $g$ $\in$ $H$ \cite{stheve}. It is worth mentioning that $${_\Gamma}\Phi(\mathfrak{qg}^\#)\subseteq \mathfrak{qg}^\#\otimes H;$$ so 
\begin{equation}
\label{2.f15}
\ad:={_\Gamma}\Phi|_{\mathfrak{qg}^\#}: \mathfrak{qg}^\#\longrightarrow \mathfrak{qg}^\#\otimes H
\end{equation}
is a $\G$--corepresentation and it fulfills (\cite{stheve})
\begin{equation}
\label{2.f15.1}
\ad \circ \pi= (\pi\otimes \id_H) \circ \Ad.
\end{equation}

On the other hand, there is a right $H$--module structure on $\mathfrak{qg}^\#$ given by 
\begin{equation}
\label{2.f16}
\theta \diamondsuit g:=\pi(hg-\epsilon(h)g)
\end{equation}
for every $\theta=\pi(h)$ $\in$ $\mathfrak{qg}^\#$. This right $H$--action satisfies $$\theta\diamondsuit g=S(g^{(1)})\theta g^{(2)},\qquad (\theta \diamondsuit g)^\ast=\theta^\ast\diamondsuit S(g)^\ast$$ as the reader can verify in Section 6.4 of \cite{stheve}.

Let $(\Gamma,d)$ be a $\ast$--FODC over $\G$ and consider the graded vector space
$$\otimes^\bullet_H\Gamma:=\bigoplus_k (\otimes^k_H\Gamma)\quad \mbox{ with } \quad \otimes^0_H\Gamma=H \quad \mbox{ and } \quad \otimes^k_H\Gamma:=\underbrace{\Gamma\otimes_H\cdots\otimes_H \Gamma}_{k\; times}$$ endowed with its canonical graded  $\ast$--algebra structure (\cite{stheve}). Now, let us consider the quotient graded space
  \begin{equation}
      \label{udtensor}     \Gamma^\wedge:=\otimes^\bullet_H\Gamma/\mathcal{Q},
  \end{equation}
  where $\mathcal{Q}$ is the two--side ideal of $\otimes^\bullet_H\Gamma$ generated by elements
  \begin{equation}
      \label{udtensor1}
     \sum_i dg_i\otimes_H dh_i \quad \mbox{ such that } \quad \sum_i g_i\,dh_i=0,
  \end{equation}
  for all $g_i$, $h_i$ $\in$ $H$. According to \cite{micho1,stheve}, the graded $\ast$--algebra structure of $\otimes^\bullet_H \Gamma$ endows $\Gamma^\wedge$ with structure of graded $\ast$--algebra. The product in $\Gamma^\wedge$ is  denoted simply by juxtaposition of elements (not by $\wedge$), as the reader can verify in \cite{micho1,micho2,micho3,stheve}; so we will adopt this notation.  

  The linear map 
  \begin{equation}
      \label{diff12}
     d:\Gamma^\wedge\longrightarrow \Gamma^\wedge 
  \end{equation}
  given by the differential of the $\ast$--FODC on $\Gamma^{\wedge\,0}= H$ (the degree--zero component), and for $t=\vartheta_1\cdots \vartheta_n$ $\in$ $\Gamma^{\wedge\,n }$ with $\vartheta_1$,..., $\vartheta_n$ $\in$ $\Gamma$, given by $$d(t)=d(\vartheta_1\cdots \vartheta_n)=\displaystyle \sum^n_{j=1}(-1)^{j-1}\vartheta_1\cdots \vartheta_{j-1}\cdot d\vartheta_j\cdot \vartheta_{j+1}\cdots \vartheta_n \; \in \; \Gamma^{\wedge\,n+1 },$$ where $$d\vartheta_j=\displaystyle \sum_l dg_l\,dh_l \qquad \mbox{ if }\qquad \vartheta_j=\displaystyle \sum_l g_l\,(dh_l),$$  is well--defined, satisfies the graded Leibniz rule, $d^2=0$ and $d(t^\ast)=(dt)^\ast$ \cite{micho1,stheve}. In other words, $d$ is the {\it natural extension} in $\Gamma^\wedge$ of the differential of $\Gamma$. In this way, 
  \begin{equation}
    \label{2.f10.1}
     (\Gamma^\wedge,d,\ast)
\end{equation}
  is a graded differential $\ast$--algebra generated by $\Gamma^{\wedge\,0}= H$  and it is called {\it the universal differential envelope $\ast$--calculus} \cite{micho1,stheve}. In general, one does not need a quantum group to define $(\Gamma^\wedge,d,\ast)$, a $\ast$--algebra is sufficient. For more details about the universal differential envelope $\ast$--calculus, see references \cite{micho1,stheve}.  
  
  In the rest of this section, we will state some results concerning the universal differential envelope $\ast$--calculus. Full details can be found in \cite{micho1,stheve}.

\begin{Proposition}
\label{1U}
Suppose $(\Omega,d_{\Omega})$ is a graded differential $\ast$--algebra generated by $\Omega^0=H$ and $(\Gamma,d)$ is a $\ast$--FODC over $\G$. Let $\varphi^{0}:\Gamma^{\wedge\,0}=H\longrightarrow  \Omega^0=H$ be a $\ast$--algebra morphism and $\varphi^{1}:\Gamma\longrightarrow \Omega $ be a linear map such that $\varphi^{1}(g\,dh)=\varphi^{0}(g)\,d_{\Omega}(\varphi^{0}(h))$ for all $g$, $h$ $\in$ $H$. Then there exists a unique family of linear maps $\varphi^{k}:\Gamma^{\wedge k}\longrightarrow \Omega$ such that $$\varphi:=\displaystyle\bigoplus_k \varphi^{k}: \Gamma^{\wedge}\longrightarrow \Omega $$ is a graded differential $\ast$--algebra morphism.
\end{Proposition}

\begin{Proposition}
\label{2U}
Suppose $(\Omega,d_{\Omega})$ is a graded differential $\ast$--algebra generated by $\Omega^0=H$ and $(\Gamma,d)$ is a $\ast$--FODC over $\G$. Let $\hat{\varphi}^{0}:\Gamma^{\wedge\,0}=H\longrightarrow \Omega^{0}=H$ be a $\ast$--antimultiplicative linear morphism and $\hat{\varphi}^{1}:\Gamma\longrightarrow \Omega^{1} $ be a linear map such that $\hat{\varphi}^{1}(g\,dh)=d_{\Omega}(\hat{\varphi}^{0}(h))\,\hat{\varphi}^{0}(g)$ for all $g$, $h$ $\in$ $H$. Then there exists a unique family of linear maps $\hat{\varphi}^{k}:\Gamma^{\wedge k}\longrightarrow \Omega$ such that $$\hat{\varphi}:=\bigoplus_k \hat{\varphi}^{k}: \Gamma^{\wedge}\longrightarrow \Omega $$ is a graded $\ast$--antimultiplicative morphism and $\hat{\varphi}\circ d=d_{\Omega}\circ \hat{\varphi}$.
\end{Proposition}

In accordance with \cite{micho2,stheve}, it can be proven that for a given $\ast$--FODC $(\Gamma, d)$ over a $\ast$--Hopf algebra, its maximal prolongation, i.e., the biggest graded differential $\ast$–algebra generated by its degree--zero elements (elements of $H$) and whose degree--one component is $\Gamma$, is $(\Gamma^\wedge, d, \ast)$.

Let $(\Gamma,d)$ be a bicovariant $\ast$--FODC over a quantum group $\G$. Then the $\ast$--Hopf algebra $H$ can be endowed with a graded differential $\ast$--algebra structure by defining $d=0$. In this way, we can consider the following tensor product of graded differential $\ast$--algebras (\cite{stheve}) $$(\Gamma^\wedge\otimes H,d_\otimes,\ast),\qquad (H\otimes \Gamma^\wedge,d_\otimes,\ast).$$  Define ${_{\Gamma^\wedge}}\Phi^{0}=\Phi^{0}_{\Gamma^\wedge}=\Delta$, ${_{\Gamma^\wedge}}\Phi^{1}={_\Gamma}\Phi$  and $\Phi^{1}_{\Gamma^\wedge}=\Phi_\Gamma$. Thus, by Proposition \ref{1U} we obtain graded differential $\ast$--algebra morphisms 
\begin{equation}
   \label{2.f10.2}
   {_{\Gamma^\wedge}}\Phi: \Gamma^\wedge \longrightarrow \Gamma^\wedge \otimes H,\qquad \Phi{_{\Gamma^\wedge}}: \Gamma^\wedge \longrightarrow H \otimes \Gamma^\wedge. 
\end{equation}
 Similarly, consider now the tensor product of $(\Gamma^\wedge,d,\ast)$ with itself (\cite{stheve}) $$ (\Gamma^\wedge\otimes \Gamma^\wedge,d_\otimes,\ast)$$ and by setting $\Delta^{0}=\Delta$ and $\Delta^{1}={_\Gamma}\Phi+ \Phi_{\Gamma}$, we can use Proposition \ref{1U} to extend the coproduct  $\Delta:H\longrightarrow H\otimes H$  to a graded differential $\ast$--algebra morphism
 \begin{equation}
\label{2.f3.5}
\Delta: \Gamma^\wedge \longrightarrow \Gamma^\wedge \otimes  \Gamma^\wedge.
\end{equation}
In particular, in accordance with \cite{micho1,stheve}, we have 
\begin{equation}
    \label{coproduc.1}
    \Delta(\theta)=\mathbbm{1}\otimes \theta+\ad(\theta).
\end{equation}
for all $\theta$ $\in$ $\mathfrak{qg}^\#$.

The counit and the coinverse can also be extended. In fact,  consider the linear map 
\begin{equation}
\label{2.f3.6}
\epsilon:\Gamma^\wedge\longrightarrow \C 
\end{equation}
defined by $\epsilon|_{H}:=\epsilon$ and $\epsilon|_{\Gamma^{\wedge k}}:=0$ for $k\geq 1$. 

On the other hand, for any $g$ $\in$ $\Ker(\epsilon)$, define $$S^1(g):=-\pi(g^{(2)})S(g^{(3)})S(S(g^{(1)})).$$ Since $\Ad(\mathcal{R})\subseteq \mathcal{R}\otimes H$, we obtain $S^1(g)=0$ for all $g$ $\in$ $\mathcal{R}$. Hence, there exists a well--defined linear map  $S^1:\mathfrak{qg}^\# \longrightarrow \Gamma.$  In accordance with \cite{micho1}, we can extend $S^1$ to the whole $\Gamma$ in such a way that the following equalities hold
 $$S^1(h\,\pi(g))=S^1(\pi(g))\,S(h)\quad \mbox{ and }\quad S^{1}(h\,dg)=d(S(g))\,S(h) $$ for all $g$, $h$ $\in$ $H$. Since Proposition \ref{2U} also holds without considering the $\ast$ operation (\cite{micho1,stheve}), we can use this statement to extend $S^0:=S$ and $S^1$ to a graded antimultiplicative linear map 
\begin{equation}
\label{2.f3.7}
S:\Gamma^\wedge\longrightarrow \Gamma^\wedge
\end{equation}
which commutes with the differential (\cite{micho1}).  These maps define a graded differential $\ast$--Hopf algebra structure 
\begin{equation}
\label{2.f3.8}
\Gamma^{\wedge\,\infty}:=(\Gamma^\wedge,\cdot,\mathbbm{1},\Delta,\epsilon, S,d,\ast)
\end{equation}
on $\Gamma^\wedge$ which extends $H^\infty=(H,\cdot,\mathbbm{1},\Delta,\epsilon, S,d,\ast)$ \cite{micho1}.  

Now it is possible to consider the right adjoint coaction of $\Gamma^\wedge$ by taking
\begin{equation}
\label{2.f11}
\Ad:\Gamma^\wedge \longrightarrow \Gamma^\wedge \otimes \Gamma^\wedge
\end{equation}
such that $$\Ad(t)=(-1)^{\partial t^{(1)}\partial t^{(2)}} t^{(2)}\otimes S(t^{(1)})t^{(3)},$$ where $\partial x$ denotes the grade of $x$ and   $(\id_{\Gamma^\wedge}\otimes \Delta)\Delta(t)=(\Delta\otimes \id_{\Gamma^\wedge})\Delta(t)=t^{(1)}\otimes t^{(2)}\otimes t^{(3)}.$  Clearly, $\Ad$ extends the right adjoint coaction $\Ad$ of $H$.

Let us define
\begin{equation}
    \label{2.f12.1}
    \begin{aligned}
    \mathfrak{qg}^{\#\wedge}&=\otimes^\bullet \mathfrak{qg}^\#/A^\wedge, \qquad \otimes^\bullet\mathfrak{qg}^\#:=\bigoplus_k (\otimes^k\mathfrak{qg}^\#)\qquad \mbox{ with }\\
    &\otimes^0\mathfrak{qg}^\#=\C\mathbbm{1},\qquad \otimes^k\mathfrak{qg}^\#:=\underbrace{\mathfrak{qg}^\#\otimes\cdots\otimes \mathfrak{qg}^\#}_{k\; times},
     \end{aligned}
\end{equation}
where $A^\wedge$ is the two--side ideal of $\otimes^\bullet\mathfrak{qg}^\#$ generated by elements $$\pi(g^{(1)})\otimes \pi(g^{(2)})\qquad \mbox{ for all }\qquad g \,\in\, \mathcal{R}.$$  In accordance with \cite{micho1}, we have that $\mathfrak{qg}^{\#\wedge}=\{ t \in \Gamma^{\wedge}\mid \Phi_{\Gamma^{\wedge}}(t)=\mathbbm{1}\otimes t\}$ and this space  is a graded differential $\ast$--subalgebra of $(\Gamma^\wedge,d,\ast)$. In particular, the product on $\mathfrak{qg}^{\#\wedge}$ will also be denoted simply by juxtaposition of elements. Furthermore, it is possible to extend the right $H$--module structure of $\mathfrak{qg}^\#$  to $\mathfrak{qg}^{\#^\wedge}$ by means of (see equation (\ref{2.f16}))
\begin{equation}
    \label{2.f12.4}
    \mathbbm{1}\diamondsuit g:=\epsilon(g),\quad (\theta_1\theta_2)\diamondsuit g:=(\theta_1\diamondsuit g^{(1)})(\theta_2\diamondsuit g^{(2)}).
\end{equation}
It is worth mentioning that (\cite{micho1,stheve})
\begin{equation}
    \label{2.f12.2}
    d\pi(g)=-\pi(g^{(1)})\pi(g^{(2)})
\end{equation}
for all $g$ $\in$ $H$. According to \cite{micho1}, the following identification holds:
\begin{equation}
    \label{2.f12.3}
    \Gamma^{\wedge}= H\otimes \mathfrak{qg}^{\#\wedge}.
\end{equation}

Let $G$ $\subset$ $M_n(\C)$  be a compact matrix Lie group  and $\G$ its associated quantum group. If the bicovariant $\ast$--FODC of $\G$ is defined by $\mathcal{R}=\Ker^2(\epsilon)=\displaystyle\{\sum^n_{i=1} a_i\,b_i\mid a_i,\,b_i\,\in\,\Ker(\epsilon),\;  n\, \in\, \N \}$, then, according to \cite{woro2,appendix} we have 
\begin{equation}
    \label{classical}
    \Gamma=H\otimes \mathfrak{g}^\#_\C,\qquad \mbox{ where }\qquad {\Ker(\epsilon)\over \mathcal{R}}={\Ker(\epsilon)\over \Ker^2(\epsilon)}=\mathfrak{g}^\#_\C
\end{equation}
is the complexification of the dual space of the Lie algebra $\mathfrak{g}$ of $G$. Moreover, $\mathfrak{qg}^{\#\wedge}=\bigwedge \mathfrak{g}^\#_\C$ is the exterior algebra of $\mathfrak{g}^\#_\C$ and by equation (\ref{2.f12.3}) we obtain $\Gamma^\wedge=H \otimes \bigwedge \mathfrak{g}^\#_\C$ \cite{appendix}. In other words, $(\Gamma^\wedge,d,\ast)$ is a subalgebra of the algebra of $\C$--valued differential forms of $G$, and by considering convergent sequences, one can recover the full algebra, even though this falls outside our purely geometric–algebraic approach. Hence, we can conclude that the universal differential envelope $\ast$--calculus is a proper generalization of the algebra of $\C$--valued differential forms of $G$ in non--commutative geometry.  In this way, for a given quantum group $\G$ and a bicovariant $\ast$--FODC $(\Gamma,d)$ over $\G$, the triplet $$(\Gamma^\wedge,d,\ast)$$ will be interpreted as the $\ast$--algebra of {\it quantum differential forms} of $\G$. In this sense,  the space 
$$\mathfrak{qg}^\#={\Ker(\epsilon)\over \mathcal{R}}$$ plays the role of the {\it quantum dual Lie algebra} and the $\ad$ corepresentation plays the role of the {\it dualization} of the right adjoint action of $G$ on $\mathfrak{g}$.

Unfortunately, not all the conditions on $H^\infty$ can be extended to $\Gamma^{\wedge\,\infty}$. In fact, we have the following example.

\begin{Example}
    \label{e.1}
    Let $\Z_2=\{0,1\}$ be the group of integers modulo $2$ and let $\mathcal{G}$ be its associated quantum group \cite{woro1}. Then $$H=\{\phi:\Z_2\longrightarrow \C \mid \phi \mbox{ is a function} \}$$ and its $\ast$--algebra structure is given by point--wise operations of function on $\mathbb{C}$ which of course, is commutative. Furthermore, the Hopf algebra structure of $H$ is the pull--back of the group structure of $\Z_2$.  The space $H$ has a canonical linear basis given by \begin{equation*}
\phi_0(x)=\left\{ 0 \;\;\mbox{ if }\;\; x=1
\atop
1 \;\; \mbox{ if }\;\; x=0,
\right. \qquad \phi_1(x)=\left\{ 0 \;\;\mbox{ if }\;\; x=0
\atop
1 \;\; \mbox{ if }\;\; x=1.
\right.
\end{equation*}
In terms of this basis, we have
\begin{equation}
    \label{coprodz_2}
    \Delta(\phi_0)=\phi_0\otimes \phi_0+\phi_1\otimes \phi_1,\,\qquad \Delta(\phi_1)=\phi_0\otimes \phi_1+ \phi_1\otimes \phi_0,
\end{equation}
\begin{equation}
    \label{coinitz_2}
    \epsilon(\phi_0)=1,\qquad \epsilon(\phi_1)=0,\qquad S(\phi_0)=\phi_0,\qquad S(\phi_1)=\phi_1.
\end{equation}

Since $\Ker(\epsilon)=\mathrm{span}_\C\{\phi_1 \}$, by Proposition \ref{0U},  the quantum group $\G$ has only two bicovariant $\ast$--FODC's: the trivial one defined by $\mathcal{R}=\Ker(\epsilon)$ and the universal one defined by $\mathcal{R}=\{0 \}$.

Let us take the universal $\ast$--FODC $(\Gamma_U,D)$ of $\G$ (see Example \ref{innese1}) and consider $(\Gamma^\wedge,d,\ast)$ its  universal differential envelope $\ast$--calculus. It is well--known that the maximal prolongation of the universal $\ast$--FODC  is the universal graded differential calculus (\cite{stheve,dv}) and therefore, the space $$
\Gamma^\wedge = \otimes_H^\bullet \Gamma_U / \mathcal{Q}$$ 
is obtained with $\mathcal{Q} = \{0\}$.

Let $R : H \otimes H \longrightarrow \mathbb{C}$ be a linear map such that $H$ is coquasitriangular with respect to $R$ (see \cite{landi}), and assume that $R$ can be extended to a map
$R : \Gamma^{\wedge} \otimes \Gamma^{\wedge} \longrightarrow \mathbb{C}$ in such a way that $\Gamma^{\wedge}$ is also coquasitriangular. However, the only possible extension of the map $R : H \otimes H \longrightarrow \mathbb{C}$ to differential forms necessarily vanishes on elements of positive degree because of the grade of the elements. Hence, identities as the following
\begin{equation}
\label{coqua}
   \phi_0\, \pi(\phi_1)=R(\pi(\phi_1)^{(1)}\otimes \phi^{(1)}_0)\,\pi(\phi_1)^{(2)}\,\phi^{(2)}_0\,\overline{R}(\pi(\phi_1)^{(3)}\otimes \phi^{(3)}_0), 
\end{equation}
can never be satisfied. For example, a straightforward calculation using the fact that
\begin{equation}
    \label{adjoinphi}
    \ad(\pi(\phi_1))=(\pi\otimes \id_H)\Ad(\phi_1)=\pi(\phi_1)\otimes \mathbbm{1},
\end{equation}
\begin{equation}
    \label{deltapi}
    \Delta(\pi(\phi_1))=\mathbbm{1}\otimes \pi(\phi_1)+\ad(\pi(\phi_1))=\mathbbm{1}\otimes \pi(\phi_1)+\pi(\phi_1)\otimes \mathbbm{1},
\end{equation}
and
$$\pi(\phi_1)^{(1)}\otimes \pi(\phi_1)^{(2)}\otimes \pi(\phi_1)^{(3)}=\pi(\phi_1)\otimes \mathbbm{1}\otimes\mathbbm{1}+\mathbbm{1}\otimes\pi(\phi_1)\otimes\mathbbm{1}+\mathbbm{1}\otimes\mathbbm{1}\otimes\pi(\phi_1)$$ proves that equation (\ref{coqua}) is incorrect. We conclude that the coquasitriangular structure cannot be extended to $\Gamma^\wedge$. 
\end{Example}

 The previous example explicitly illustrates that even for the simplest non--trivial $\ast$–Hopf algebra equipped with its universal graded differential calculus, it is impossible to extend the coquasitriangular structure to the level of differential forms.

\subsection{Basic Theory of Quantum Principal Bundles}
At this point, we introduce the notions of quantum principal bundles and 
quantum principal connections. As we have mentioned earlier, our work is based on the theory developed by M. Durdevich, although we modify the definition of quantum principal connections and adopt the notation used in \cite{sald1}. For the purposes of this paper, we will present all the necessary basics of Durdevich's formulation while explaining the motivation behind some of the definitions. Of course, there are similarities between Durdevich's formulation and the one presented in \cite{libro}, especially in degrees $0$ and $1$; however, there are also differences, particularly for degrees $2$ and higher. The reader is encouraged to consult the original works \cite{micho1, micho2, micho3} for further details.

Let $(B,\cdot,\mathbbm{1},\ast)$ be a quantum space and let $\G$ be a quantum group. A {\it quantum principal $\G$--bundle} over $B$ (abbreviated ``qpb") is a quantum structure formally represented by the triplet 
\begin{equation}
\label{2.f17}
\zeta=(P,B,\Delta_P),
\end{equation}
where $(P,\cdot,\mathbbm{1},\ast)$ is called the {\it quantum total space}, and $(B,\cdot,\mathbbm{1},\ast)$ is a $\ast$--subalgebra, which receives the name {\it quantum base space}. Furthermore, $$\Delta_P:P \longrightarrow P\otimes H$$ is a $\ast$--algebra morphism that satisfies
\begin{enumerate}
\item $\Delta_P$ is a $\G$--corepresentation.
\item $\Delta_P(x)=x\otimes \mathbbm{1}$ if and only if $x$ $\in$ $B$.
\item The linear map 
\begin{equation}
    \label{porqueyo.3}
    \beta:P\otimes P\longrightarrow P\otimes H
\end{equation}
given by $$\beta(x\otimes y):=x\cdot \Delta_P(y):=(x\otimes \mathbbm{1})\cdot \Delta_P(y) $$ is surjective. 
\end{enumerate}

A motivation for this definition comes from the notion of principal bundles in differential geometry. A  detailed explanation of this can be found in \cite{stheve,libro}. In general, there is no need to work with a quantum group; a $\ast$--Hopf algebra is sufficient (\cite{micho2,stheve}) and it is worth mentioning that the previous definition is intrinsically related with the notion of Hopf--Galois extension, as the reader can verify in \cite{libro}.

Given $\zeta$ a qpb over $B$, a {\it differential calculus} on it is:
 \begin{enumerate}
 \item A graded differential $\ast$--algebra $(\Omega^\bullet(P),d,\ast)$ generated by its degree--zero elements $\Omega^0(P)=P$ ({\it quantum differential forms of $P$}).
 \item  A bicovariant $\ast$--FODC $(\Gamma,d)$ over $\G$ and its universal differential envelope $\ast$--calculus $(\Gamma^\wedge,d,\ast)$.
 \item The map $\Delta_P$ is extendible to a graded differential $\ast$--algebra morphism $$\Delta_{\Omega^\bullet(P)}:\Omega^\bullet(P)\longrightarrow \Omega^\bullet(P)\otimes \Gamma^{\wedge}.$$ Here we have considered that $\otimes$ is the tensor product of graded differential $\ast$--algebras.
 \end{enumerate}

Notice that if $\Delta_{\Omega^\bullet(P)}$ exists, then it is unique because all our graded differential $\ast$--algebras are generated by their degree--zero elements. Furthermore,  $\Delta_{\Omega^\bullet(P)}$ is a graded differential $\Gamma^\wedge$--corepresentation on $\Omega^\bullet(P)$ \cite{micho2}. In this way, the space of horizontal forms is defined as 
\begin{equation}
\label{2.f18}
\Hor^\bullet P\,:=\{\varphi \in \Omega^\bullet(P)\mid \Delta_{\Omega^\bullet(P)}(\varphi)\, \in \, \Omega^\bullet(P)\otimes H \},
\end{equation}
and it is a graded  $\ast$--subalgebra of $\Omega^\bullet(P)$ \cite{stheve}. Since $$\Delta_{\Omega^\bullet(P)}(\Hor^\bullet P)\subseteq \Hor^\bullet P\otimes H,$$ the map 
\begin{equation}
\label{2.f19}
\Delta_\Hor:=\Delta_{\Omega^\bullet(P)}|_{\Hor^\bullet P}: \Hor^\bullet P \longrightarrow \Hor^\bullet P\otimes H
\end{equation}
is a $\G$--corepresentation on $\Hor^\bullet P$. Also, one can define the space of {\it base} forms ({\it quantum differential forms of $B$}) as 
\begin{equation}
\label{2.f20}
\Omega^\bullet(B):=\{\mu \in \Omega^\bullet(P)\mid \Delta_{\Omega^\bullet(P)}(\mu)=\mu\otimes \mathbbm{1}\}.
\end{equation}
The space of base forms is a graded differential $\ast$--subalgebra of $(\Omega^\bullet(P),d,\ast)$. In general, it is not generated by $\Omega^0(B)=B$ and an explicit example of this fact can be found in \cite{appendix}. 

It is worth mentioning that we can define the graded differential 
$\ast$--algebra of vertical forms using the graded vector space 
$P\otimes \mathfrak{qg}^{\#\,\wedge}$. However, for the purposes of this paper, 
it is unnecessary to cover this aspect of the theory. Interested readers are 
encouraged to consult \cite{micho2, stheve}.

By {\it dualizing} the notion of principal connections in differential geometry (\cite{nodg}), we introduce the concept of quantum principal connection (abbreviated ``qpc") on a quantum principal $\G$--bundle $\zeta$ with a differential calculus. It is defined as a linear map
\begin{equation}
\label{2.f24}
\omega:\mathfrak{qg}^\#\longrightarrow \Omega^{1}(P) 
\end{equation}
satisfying  
\begin{equation}
    \label{2.f24.qpc}
    \Delta_{\Omega^\bullet(P)}(\omega(\theta))=(\omega\otimes \id_H)\ad(\theta)+\mathbbm{1}\otimes\theta
\end{equation}
for all $\theta$ $\in$ $\mathfrak{qg}^\#$, where $\ad$ is the $\G$--corepresentation given in equation (\ref{2.f15}).

In analogy with the {\it classical} case, it can be proved that the set 
\begin{equation}
\label{2.f24.1}
\mathfrak{qpc}(\zeta):=\{\omega:\mathfrak{qg}^\#\longrightarrow \Omega^1(P)\mid \omega \mbox{ is a qpc on }\zeta \}
\end{equation}
is not empty, and it is an affine space modeled by the vector space of {\it connection displacements} (\cite{micho2,stheve}) 
\begin{equation}
    \label{2.f24.2}
    \overrightarrow{\mathfrak{qpc}(\zeta)}:=\Mor^1(\ad,\Delta_\Hor),
\end{equation}
where
\begin{equation*}
\Mor^1(\ad,\Delta_\Hor)=\{\lambda:\mathfrak{qg}^\# \longrightarrow \Hor^1 P\mid \lambda \mbox{ is linear and } (\lambda\otimes \id_H)\circ \ad=\Delta_\Hor \circ \lambda \}.
\end{equation*}

There is a canonical involution on $\mathfrak{qpc}(\zeta)$ given by
\begin{equation}
\label{2.f24.3}
\widehat{\omega}:=\ast\circ \omega\circ \ast,
\end{equation}
and we define the dual qpc of $\omega$ as $\widehat{\omega}$. A qpc $\omega$ is real if 
\begin{equation}
\label{2.f24.6}
\widehat{\omega}=\omega
\end{equation}
and we say that it is imaginary if  $\widehat{\omega}=-\omega.$  It is worth mentioning that in order to embrace a more general theory, our definition of qpc's is different that the standard one in the literature for example  \cite{micho1, micho2, micho3, stheve}, because in these references, qpc's are always real, i.e., in these references, equation (\ref{2.f24.6}) is part of the definition of qpc's. 

A qpc is called {\it regular} if for all $\varphi$ $\in$ $\Hor^{k}P$ and $\theta$ $\in$ $\mathfrak{qg}^\#$, we have 
\begin{equation}
\label{2.f25}
\omega(\theta)\,\varphi=(-1)^{k}\varphi^{(0)}\omega(\theta\diamondsuit\varphi^{(1)}), 
\end{equation}
where $\Delta_\Hor(\varphi)=\varphi^{(0)}\otimes\varphi^{(1)}.$ A qpc $\omega$ is called {\it multiplicative} if for all $g$ $\in$ $\mathcal{R}$ we obtain
\begin{equation}
\label{2.f26}
\omega(\pi(g^{(1)}))\omega(\pi(g^{(2)}))=0.
\end{equation} 

It is worth noting that every qpc that comes from the {\it dualization} of a {\it classical} principal connection is  regular and multiplicative.

By {\it dualizing} the notion of the covariant derivative of a principal connection in differential geometry (\cite{nodg}), the covariant derivative of a qpc is defined as the projection of $d$ onto the space of horizontal 
forms. In accordance with \cite{micho2,stheve},  for all $\varphi$ $\in$ $\Hor^k P$, the difference $d\varphi-(-1)^{k}\varphi^{(0)}\omega(\pi(\varphi^{(1)}))$ lies in $\Hor^{k+1} P$. In this way, for a given qpc $\omega$, we define its {\it covariant derivative} as the first--order linear map 
\begin{equation}
\label{2.f30}
D^{\omega}: \Hor^\bullet P \longrightarrow \Hor^\bullet P
\end{equation}
such that for every $\varphi$ $\in$ $\Hor^k P$ $$
D^{\omega}(\varphi)=  d\varphi-(-1)^{k}\varphi^{(0)}\omega(\pi(\varphi^{(1)})).$$ 
On the other hand, the first--order linear map 
\begin{equation}
\label{2.f30.1}
\widehat{D}^{\omega}:=\ast\circ D^{\omega}\circ\ast: \Hor^\bullet P\longrightarrow \Hor^\bullet P
\end{equation}
is called the {\it dual covariant derivative} of $\omega$. Explicitly, we have $$\widehat{D}^{\omega}(\varphi)=d\varphi +\widehat{\omega} (\pi(S^{-1}(\varphi^{(1)})))\varphi^{(0)}$$ for every $\varphi$ $\in$ $\Hor^k P$. It should be clear that, in general, $\widehat{D}^\omega\not= D^{\widehat{\omega}}$. In addition, we have
\begin{equation}
    \label{covd}
    D^\omega\not=\widehat{D}^\omega.
\end{equation}
In fact, 
\begin{equation}
    \label{covd.1}
    \widehat{D}^{\omega}(\varphi)=D^{\omega}(\varphi)+\ell^{\widehat{\omega}}(\pi(S^{-1}(\varphi^{(1)})),\varphi^{(0)})+(-1)^k\varphi^{(0)}(\omega-\widehat{\omega})(\pi(\varphi^{(1)})),
\end{equation}
where 
\begin{equation*}
\begin{aligned}
\ell^{\omega}:\mathfrak{qg}^\#\times \Hor^\bullet P &\longrightarrow \Hor^\bullet P \\
(\theta\;,\;\varphi)\qquad&\longmapsto \omega(\theta)\varphi-(-1)^k \varphi^{(0)}\omega(\theta\diamondsuit \varphi^{(1)}).
\end{aligned}
\end{equation*}
The map $\ell^{\omega}$ {\it measures the degree of non--regularity} of $\omega$, in the sense of $\ell^{\omega}=0$ if and only if $\omega$ is regular. In this way, for real and regular qpc's we obtain $D^\omega=\widehat{D}^{\omega}=D^{\widehat{\omega}}$, which is the situation for qpc's arising from {\it classical} principal connections. In other words, $D^\omega$ and $\widehat{D}^\omega$ are two different horizontal operators that generalize the covariant derivative of a principal connection in differential geometry. In the next section, we will work with both operators. 

Direct calculations prove that (\cite{micho2,micho3})
\begin{equation}
\label{2.f31}
D^{\omega}, \, \widehat{D}^{\omega}\, \in\,\Mor(\Delta_\Hor,\Delta_\Hor)\,,\quad D^{\omega}|_{\Omega^\bullet(B)}=\widehat{D}^{\omega}|_{\Omega^\bullet(B)} =d|_{\Omega^\bullet(B)}
\end{equation}
and 
\begin{equation}
\label{2.f32}
D^{\omega}(\varphi\psi)=D^{\omega}(\varphi)\psi+(-1)^k\varphi D^{\omega}(\psi)+(-1)^k \varphi^{(0)}\ell^{\omega}(\pi(\varphi^{(1)}),\psi),
\end{equation}
\begin{equation}
\label{2.f32.1}
\widehat{D}^{\omega}(\varphi\psi)=\widehat{D}^{\omega}(\varphi)\psi+(-1)^k\varphi \widehat{D}^{\omega}(\psi)+ \ell^{\widehat{\omega}}(\pi(S^{-1}(\psi^{(1)}))\diamondsuit S^{-1}(\varphi^{(1)}),\varphi^{(0)})\psi^{(0)},
\end{equation}
for all $\varphi$ $\in$ $\Hor^k P$, $\psi$ $\in$ $\Hor^\bullet P$. Moreover, for real qpc's
\begin{equation}
\label{2.f32.5}
D^{\omega}(\psi)^\ast=D^{\omega}(\psi^\ast)+\ell^{\omega}(\pi(S(\psi^{(1)})^\ast),\psi^{(0)\ast})=\widehat{D}^{\omega}(\psi^\ast).
\end{equation}

Notice that $\omega$ is real and regular if and only if $D^{\omega}$ and $\widehat{D}^{\omega}$ satisfy the graded Leibniz rule; this is the main reason to study real and regular qpc's \cite{micho2}. 

Let us define the map $r^\omega$ given by
\begin{equation}
\label{falsecur}
\begin{aligned}
r^\omega:H &\longrightarrow \Hor^2 P \\
g&\longmapsto d(\omega(\pi(g)))+\omega(\pi(g^{(1)}))\omega(\pi(g^{(2)})).
\end{aligned}
\end{equation}
Then we have (\cite{micho3})
$$D^{\omega\,2}(\varphi)=-\varphi^{(0)}r^{\omega}(\varphi^{(1)})\quad \mbox{ and }\quad  \widehat{D}^{\omega\,2}(\varphi)=-r^{\omega}(S^{-1}(\varphi^{(1)}))\varphi^{(0)}.$$

Let $\pi:P\longrightarrow B$ be a {\it classical} principal $G$--bundle ($P$ is the total space, $B$ is the base space and $\pi$ is the bundle projection) with a principal connection $\omega$. Then, the curvature $R^\omega$ of $\omega$ is defined as the $\mathfrak{g}$--valued differential $2$--form of $P$ given by \cite{nodg}
\begin{equation}
    \label{classcur}
    R^\omega=d\omega + {1\over 2}[\omega\wedge \omega], \quad R^\omega_x(X_x,Y_x)=d\omega_x(X_x,Y_x)+[\omega(X_x),\omega(Y_x)]_{\mathfrak{g}},
\end{equation}
where  $[-,-]_{\mathfrak{g}}$ denotes the Lie bracket of $\mathfrak{g}$, the Lie algebra of $G$, and $X_x$, $Y_x$ $\in$ $T_xP$, $x$ $\in$ $P$. In addition,  the square of the covariant derivative is related to $R^\omega$ \cite{nodg}.  In this way, one could take equation (\ref{falsecur}) as the definition of the curvature for a qpc.

However, in differential geometry, the curvature is $\mathfrak{g}$--valued, so it is natural to expect that in the {\it non--commutative geometrical} setting the curvature be defined on $\mathfrak{qg}^\#$, as the {\it dualization} of the {\it classical} case indicates. Moreover, in  differential geometry, the curvature is a basic form of type $\ad$ \cite{nodg}; therefore, for qpc's, the curvature must be 
an element of $$\Mor^2(\ad,\Delta_\Hor)=\{\lambda:\mathfrak{qg}^\# \longrightarrow \Hor^2 P\mid \lambda \mbox{ is linear and } (\lambda\otimes \id_H)\circ \ad=\Delta_{\Hor} \circ \lambda \}.$$ 

Of course, $r^\omega$ induces a well--defined map on $\mathfrak{qg}^\#=\Ker(\epsilon)/ \mathcal{R}$ if and only if $r^\omega(\mathcal{R})=0$, and this relation holds if and only if $\omega$ is multiplicative. This is problematic since, in general, there exist qpc's that are not multiplicative \cite{micho2}.

To define the curvature as a map from $\mathfrak{qg}^\#$ to $\Hor^2(P)$ for every qpc, we begin by defining the following auxiliary map

\begin{Definition}[Embedded differential]
    \label{embeddeddifferential}
    An embedded differential is a linear map $$\Theta:\mathfrak{qg}^\#\longrightarrow \mathfrak{qg}^\#\otimes \mathfrak{qg}^\# $$ such that
\begin{enumerate}
    \item $\ad^{\otimes 2}\circ \Theta=(\Theta\otimes \id_H) \circ \ad,$ where $\ad^{\otimes 2}:=M\circ (\ad\otimes \ad)$, with $$M: \mathfrak{qg}^\#\otimes H \otimes \mathfrak{qg}^\#\otimes H\longrightarrow \mathfrak{qg}^\#\otimes \mathfrak{qg}^\#\otimes H$$ given by $M(\theta_1\otimes g_1,\theta_2 \otimes g_2)=\theta_1\otimes \theta_2 \otimes  g_1g_2$.
    \item If $\Theta(\theta)=\displaystyle\sum^n_{i,j=1}\theta_i\otimes \theta'_j$, then $d\theta=\displaystyle\sum^n_{i,j=1}\theta_i \theta'_j$ and $\Theta(\theta^\ast)=-\displaystyle\sum^n_{i,j=1}\theta'^\ast_j \otimes \theta^\ast_i$ (see equation (\ref{2.f12.2})).
\end{enumerate}
\end{Definition}

In general, an embedded differential can be constructed by choosing a $\ast$--$S$--invariant $\Ad$--invariant complement $L \subset \Ker(\epsilon)$ of $\mathcal{R}$ and taking  $\Theta=-(\pi\otimes \pi)\circ \Delta\circ \pi^{-1}|_{L}.$  Notice that choosing $\Theta$ is choosing a compatible way (with respect to the differential structure) to embed $\mathfrak{qg}^\#$ into $\mathfrak{qg}^\#\otimes \mathfrak{qg}^\#$. 

Fix any such embedded differential $\Theta$. We define
the {\it curvature} of a qpc $\omega$ as the linear map
\begin{equation}
\label{2.f28}
R^{\omega}:=d\omega-\langle \omega,\omega\rangle : \mathfrak{qg}^\#\longrightarrow  \Omega^2(P)
\end{equation}
with $$\langle \omega,\omega\rangle:=m_\Omega\circ(\omega\otimes \omega)\circ\Theta: \mathfrak{qg}^\#\longrightarrow \Omega^2(P),$$
where $m_\Omega:\Omega^\bullet(P)\otimes \Omega^\bullet(P)\longrightarrow \Omega^\bullet(P)$ is the product map. By the properties of $\Theta$, it can be proven that $\Im(R^{\omega})\subseteq \Hor^2 P$ and
\begin{equation}
\label{2.f29}
R^{\omega}\,\in\, \Mor^2(\ad,\Delta_\Hor)
\end{equation}
for every qpc $\omega$ \cite{micho2,stheve}. In other words, $R^\omega$  is a {\it quantum} basic form of type $\ad$.

 If $\omega$ is multiplicative, then $R^\omega$ does not depend on the choice of $\Theta$ and it agrees with the induced map of $r^\omega$ on $\mathfrak{qg}^\#$ (\cite{micho2, stheve}). This is the main reason to study multiplicative qpc's. The presence of $\Theta$ in the definition of $R^\omega$ in the {\it non--commutative geometrical} case can be interpreted as a {\it quantum phenomenon} in which there can be several non--equivalent ways to embed $\mathfrak{qg}^\#$ into $\mathfrak{qg}^\#\otimes \mathfrak{qg}^\#$ and produce horizontal quadratic expressions with $\omega$. Therefore, the definition of $R^\omega$ must involve $\Theta$ if we want $R^\omega$ to be defined on $\mathfrak{qg}^\#$ and to lie in $\Mor^2(\ad,\Delta_\Hor)$.

\begin{Remark}
\label{rema}
From this point onward until the end of the paper, we shall restrict our attention exclusively to qpb's for which the quantum base space $(B,\cdot,\mathbbm{1},\ast)$ can be completed to a $C^\ast$--algebra. 

According to \cite{micho3}, in this case, for every $\delta^V$ $\in$ $\T$ there exists a set  $$\{T^\l_k \}^{d_{V}}_{k=1} \subseteq \Mor(\delta^V,\Delta_P)$$ for some $d_{V}$ $\in$ $\N$ such that
\begin{equation}
    \label{generators}
\sum^{d_{V}}_{k=1}x^{V\,\ast}_{ki}x^{V}_{kj}=\delta_{ij}\mathbbm{1},
\end{equation}
where $x^{V}_{ki}:=T^\l_k(e_i)$ and $\delta_{ij}$ is the Kronecker delta. Here, $\T$  is a complete set of mutually non--equivalent irreducible finite--dimensional $\G$--corepresentations with $\delta^\C_\triv$ $\in$ $\T$, and $\{e_i\}^{n_{V}}_{i=1}$ is the orthonormal basis of $V$ given in Theorem \ref{rep}.
\end{Remark}

\noindent The {\it dualization} of the following proposition motivates Remark \ref{rema}. As we have just mentioned, a proof of 
the existence of the maps $\{T^\l_k\}_{k=1}^{d_V}$ in the {\it non--commutative geometrical} setting can be found in \cite{micho3}.

\begin{Proposition}
\label{difgeo}
    Let $G$ $\subset$ $M_n(\C)$ be a compact matrix Lie group and let $\G$ be its associated quantum group \cite{woro1}. Take a classical principal $G$--bundle $\pi:P\longrightarrow B$, where $P$ is the total space, $B$ is the base space and $\pi$ is the bundle projection. Assume $P$, $B$ are compact. If $\T$ is a complete set of mutually non--equivalent irreducible $\G$--corepresentations, then for every  $\delta^V$ $\in$ $\T$, there exists $$\{T^\l_k\}^{d_V}_{k=1}\;\subseteq\; \Mor(\delta^V, \Delta_P)$$ for some $d_V$ $\in$ $\N$ such that equation (\ref{generators}) holds. Here $$\Delta_P:C^\infty_\C(P)\longrightarrow C^\infty_\C(P\times G)\supset C^\infty_\C(P)\otimes C^\infty_\C(G)$$ is the pull--back of the right $G$--action on $P$; and $C^\infty_\C(P)$,  $C^\infty_\C(G)$, $C^\infty_\C(P\times G)$ denote the spaces of $\C$--valued smooth functions on $P$, $G$ and $P\times G$, respectively.
\end{Proposition}

\begin{proof}
     Let $\delta^V$ $\in$ $\T$ and assume $n=\mathrm{dim}_\C(V)$. Then, the corepresentation  $\delta^V$ induces a $G$--action  $\alpha: G\times V\longrightarrow V $ given by 
     \begin{equation}
         \label{actionclass}
         \alpha(A,e_j)=\displaystyle\sum^n_{i=1} g^V_{ij}(A)\,e_i,
     \end{equation}
      where the elements $\{ g^V_{ij}: G \longrightarrow \C \}$ are the ones given in Theorem \ref{rep} and $\{e_i \}^n_{i=1}$ is the corresponding orthonormal basis of $V$.  
     
     For each $b$ $\in$ $B$, let $(U_b, \Psi_b)$ be a principal $G$--bundle local trivialization of $\pi:P\longrightarrow B$ (\cite{nodg}). Since $\{U_b \}_{b\in B}$ is an open cover, by compactness, there exist points $b_1, \ldots, b_r \in B$ such that $\{\,U_{b_i}\}_{i=1}^r $ remains an open cover of $B$. Let $\{\rho_{b_i} \}_{i=1}^r$ be a partition of unity subordinate to the open cover $\{\,U_{b_i}\}_{i=1}^r$, where each $\rho_{b_i}$ has compact support and admits a smooth square root. Then, the maps
     \begin{equation*}
f^i_j :P \longrightarrow V, \qquad x \longmapsto v=\sqrt{\rho_{b_i}(\pi(x))}\;\alpha(A(x)^{-1},e_j),
\end{equation*}
where $A(x)$ is the unique element of $G$ such that $\Psi_{b_i}(\pi(x),A(x))=x$, are smooth because they are defined by compositions and products of smooth functions (\cite{nodg}).  Define the smooth functions
\begin{equation*}
f^i_{jk} :P \longrightarrow \C, \qquad x \longmapsto \langle e_k \mid f^i_j(x)\rangle,
\end{equation*}
where $\langle - | -\rangle$ is the inner product (antilinear in the second coordinate) that makes $\delta^V$ unitary. Now, let us consider the linear maps $$T^i_j: V \longrightarrow C^\infty_\C(P)$$ given by $T^i_{j}(e_k)=f^i_{jk}$. A direct calculation shows that for all $x$ $\in$ $P$ and for all $A$ $\in$ $G$, $$\left(\sum^{r,n}_{l,k=1}f^{l\ast}_{ki}\,f^{l}_{kj} \right)(x)=\delta_{ij},\qquad (\Delta_P\circ T^i_j)(x,A)=((T^i_j\otimes \id_H)\circ \delta^V)(x,A)$$ and hence, the proposition follows by taking $T^\l_1=T^1_1$, $T^\l_2=T^1_2$,..., $T^\l_n=T^1_n$, $T^\l_{n+1}=T^2_1$,..., $T^\l_{d_V}=T^r_n$ with $d_V=rn$.  
\end{proof}

 It is worth mentioning that the maps $T^i_j$ agree with the pull--back of $f^i_j$, once the dual space $V^\#$ of $V$ is identified with $V$. In addition, it can be proven that $\{f^i_j\}$ is a set of generators of the $C^\infty_\C(B)$--bimodule of $G$--equivariant maps $$C^\infty_{\C}(P,V)^G= \{f:P\longrightarrow V \mid f \mbox{ is smooth and } f(xA)=\alpha(A^{-1})f(x) \mbox{ for all } x \,\in \, P,\, A \,\in\, G \},$$ as the reader can check in \cite{appendix}. Finally, it is worth mentioning that the {\it non--commutative geometrical} counterpart of  $C^\infty_{\C}(P,V)^G$  is the space  $\Mor(\delta^V, \Delta_P).$ 
 
 In the next section, we will prove that, for every qpb (which of course, includes the dualization of a {\it classical} principal bundle), the maps $\{T^\l_i\}_{i=1}^{d_V}$ of Remark \ref{rema} always form a set of left $B$--generators of $\Mor(\delta^V, \Delta_P)$. The superscript $\l$ in $T^\l_k$ is purely symbolic, used to indicate that these maps are left $B$--generators.

 In the context of the theory of Hopf--Galois extensions, equation (\ref{generators}) implies that $P$ is principal \cite{kt}. Moreover, under equation (\ref{generators}), real qpc's always exist for every qpb with a differential calculus \cite{micho2}.  In the final section, we present the explicit form of the maps $\{T^\l_k \}^{d_{V}}_{k=1}$ in some concrete examples.

\section{Associated Quantum Vector Bundles, Induced Quantum Linear Connections and Hermitian Structures}

The primary purpose of this paper is to present some of the essential aspects of associated quantum vector bundles, induced quantum linear connections and the definition of the canonical Hermitian structure, as well as the relationships among these structures, illustrating their analogy with the {\it classical} case. In this section, we shall deal with all of these topics. Since we are not interested in the categorical point of view, we can weaken some conditions imposed in \cite{sald1}.

\subsection{Associated Quantum Vector Bundles and Induced Quantum Linear Connections}

Let us begin by considering a quantum $\G$--bundle $\zeta=(P,B,\Delta_P)$ and a $\G$--corepresentation $\delta^V$ $\in$ $\T$. Notice that the $\C$--vector space $\Mor(\delta^V,\Delta_P)$ has a natural $B$--bimodule structure given by multiplication with elements of $B$, i.e., for $b$ $\in$ $B$ and $T$ $\in$ $\Mor(\delta^V,\Delta_P)$, the operations
$$(b,T)\longmapsto b\,T,\quad \mbox{ where }\quad b\,T: V\longrightarrow P $$ is given by $(b\,T)(v)=b\,T(v)$ for all $v$ $\in$ $V$,  and 
$$(T,b)\longmapsto T\,b,\quad \mbox{ where }\quad T\,b: V\longrightarrow P $$ is given by $(T\,b)(v)=T(v)\,b$ for all $v$ $\in$ $V$, induces a  $B$--bimodule structure on $\Mor(\delta^V,\Delta_P)$.

Our first objective is to show that $\Mor(\delta^V,\Delta_P)$ is finitely generated and projective, as both a left $B$--module and a right $B$--module \cite{micho2}. 

Let $b$ $\in$ $B$ and consider the element $\displaystyle\sum^{n_{V}}_{i=1}x^{V}_{ki}\,b\,x^{V\,\ast}_{li}$, with $T^\l_k(e_i)=x^V_{ij}$ (see Remark \ref{rema}). By equation (\ref{2.f9}) it is easy to check that $$\Delta_P\left(\sum^{n_{V}}_{i=1}x^{V}_{ki}\,b\,x^{V\,\ast}_{li}\right)=\sum^{n_{V}}_{i=1}x^{V}_{ki}\,b\,x^{V\,\ast}_{li}\otimes \mathbbm{1};$$ so  $\displaystyle\sum^{n_{V}}_{i=1}x^{V}_{ki}\,b\,x^{V\,\ast}_{li}$ $\in$ $B$ and we can define
\begin{equation}
\label{3.f1}
\varrho^{V}_{kl}:B \longrightarrow B,\qquad b \longmapsto \varrho^{V}_{kl}(b)=\sum^{n_{V}}_{i=1}x^{V}_{ki}\,b\,x^{V\,\ast}_{li},
\end{equation}
where $k$, $l$ $\in$ $\{1,...,d_{V} \}$ (recall that $\delta^V$ is irreducible, see Remark \ref{rema}). 

\begin{Proposition}
    \label{innessari2}
    For all $b$, $a$ $\in$ $B$ we have $$\varrho^{V}_{kl}(b)^\ast=\varrho^{V}_{lk}(b^\ast)\quad \mbox{ and }\quad  \displaystyle \sum^{n_{V}}_{i=1}\varrho^{V}_{ki}(b)\varrho^{V}_{il}(a)=\varrho^{V}_{kl}(ba).$$
\end{Proposition}
\begin{proof}
    Let $b$ $\in$ $B$. Then $$\varrho^{V}_{lk}(b^\ast)=\sum^{n_{V}}_{i=1}x^{V}_{li}\,b^\ast\,x^{V\,\ast}_{ki}=\left( \sum^{n_{V}}_{i=1}x^{V}_{ki}\,b\,x^{V\,\ast}_{li} \right)^\ast=\varrho^{V}_{kl}(b)^\ast.$$

  On the other hand, by equation (\ref{generators}), we get
  \begin{eqnarray*}
        \sum^{n_{V}}_{i=1}\varrho^{V}_{ki}(b)\varrho^{V}_{il}(a)=\sum^{n_{V}}_{i,j,s=1}x^V_{kj}\,b\,x^{V\,\ast}_{ij}\,x^V_{is}\,a\,x^{V\,\ast}_{ls}= \sum^{n_{V}}_{j,s=1}x^V_{kj}\,b\,\delta_{js}\,a\,x^{V\,\ast}_{ls}&=& \sum^{n_{V}}_{j=1}x^V_{kj}\,b\,\,a\,x^{V\,\ast}_{lj}
        \\
        &=&
        \varrho^{V}_{kl}(ba)
\end{eqnarray*}
for all $b$, $a$ $\in$ $B$.
\end{proof}

In light of the last proposition, there exists a linear, multiplicative, $\ast$--preserving (in general not--unital) map 
\begin{equation}
\label{3.f2}
\varrho^{V}: B  \longrightarrow M_{d_{V}}(B),\qquad 
b  \longmapsto \varrho^{V}(b)=(\varrho^{V}_{kl}(b)),
\end{equation}
    \noindent where $M_{d_{V}}(B)$ denotes the space of $d_{V}\times d_{V}$ matrices with entries in $B$ \cite{micho2}. The $\ast$ operation of $M_{d_{V}}(B)$ will be denoted by $\dagger$ and it is defined as the composition of the $\ast$ operation on $B$ (applied entry--wise on $M_{d_{V}}(B)$) with the usual matrix transposition. Next, consider the free left $B$--module $B^{d_{V}}$ with its canonical basis $\{\overline{e}_1,...,\overline{e}_{d_{V}} \}$ and the left $B$--submodule $$B^{d_{V}}\cdot \varrho^{V}(\mathbbm{1})\;\subseteq\;B^{d_{V}}.$$

Define the left $B$--module morphism
\begin{equation}
\label{3.f4}
H':B^{d_{V}}  \longrightarrow \Mor(\delta^V,\Delta_P)
\end{equation}
such that $H'(\overline{e}_k)=T^\l_k.$ We claim that for all $b$ $\in$ $B$ and every $\overline{b}=(b_1,...,b_{d_V}) \in B^{d_V}$, the following identity holds:
 $$H'(\overline{b}\cdot\varrho^{V}(b))=H'(\overline{b})\,b. $$ In fact, for each  $\overline{e}_i$ $\in$ $B^{d_V}$, we have $$H'(\overline{e}_i\cdot \varrho^{V}(b))=\displaystyle \sum^{d_{V}}_{j=1}H'(\varrho^{V}_{ij}(b)\,\overline{e}_j)=\sum^{V}_{j=1}\varrho^{V}_{ij}(b)H'(\overline{e}_j)=\sum^{d_{V}}_{j=1}\varrho^{V}_{ij}(b)T^\l_j.$$  Evaluating this expression on the orthonormal basis $\{e_i \}^{n_V}_{i=1}$ of $V$ (see Theorem \ref{rep} and   Remark \ref{rema})  and using equation (\ref{generators}) we obtain
\begin{eqnarray*}
\left(\sum^{d_{V}}_{j=1}\varrho^{V}_{ij}(b)T^\l_j\right)(e_s) \;=\; \sum^{{d_{V}},n_{V}}_{j,k=1}x^{V}_{ik}\,b\, x^{V\, \ast}_{jk}T^\l_j(e_s) \;=\; \sum^{{d_{V}},n_{V}}_{j,k=1}x^{V}_{ik}\,b\, x^{V\, \ast}_{jk}x^{V}_{js}  &= & \sum^{n_{V}}_{k=1}x^{V}_{ik}\,b\, \delta_{ks}
   \\
  &= & 
x^{V}_{is}\,b=T^\l_i(e_s)\,b,  
\end{eqnarray*}
so $H'(\overline{e}_i\cdot \varrho^{V}(b))=T^\l_i\,b =H'(\overline{e}_i)\,b$ and by $B$--linearity we conclude  $H'(\overline{b}\cdot\varrho^{V}(b))=H'(\overline{b})\,b$. In particular  
\begin{equation}
\label{6.f1.11}
H'(\overline{b}\cdot \varrho^{V}(\mathbbm{1}))=H'(\overline{b})
\end{equation}
for all $\overline{b}$ $\in$ $B^{d_V}$.

\begin{Proposition}
\label{6.1.6}
The map $\widetilde{H}=H'|_{B^{d_{V}}\cdot \varrho^{V}(\mathbbm{1})}:B^{d_{V}}\cdot \varrho^{V}(\mathbbm{1})\longrightarrow \Mor(\delta^V,\Delta_P) $ is a bijection.
\end{Proposition}

\begin{proof}
Let $\overline{b}\cdot \varrho^{V}(\mathbbm{1})$ $\in$ $\Ker(\widetilde{H})$ with  $\overline{b}=\displaystyle \sum^{d_{V}}_{i=1}b_i\overline{e}_i $. Then,  by equation (\ref{6.f1.11}) we have that $$\widetilde{H}(\overline{b}\cdot \varrho^{V}(\mathbbm{1}))=H'(\overline{b})=\displaystyle \sum^{d_{V}}_{k=1}b_kT^\l_k=0$$ and evaluating in the basis $\{e_i\}^{n_V}_{i=1}$ we obtain $\displaystyle \sum^{d_{V}}_{k=1}b_k\,x^{V}_{ki}=0$ for every $i=1,...,n_V$; so   $\displaystyle \sum^{d_{V}}_{k,i=1}b_k\,x^{V}_{ki}\,x^{V\,\ast}_{li}=0$  for every $l$ $\in$ $\{1,...,{d_{V}} \}$. In other words  $$0=\displaystyle \sum^{d_{V}}_{k,i=1}b_k\,x^{V}_{ki}\,x^{V\,\ast}_{li}=\overline{b}\cdot \varrho^{V}(\mathbbm{1})$$  and thus, $\Ker(\widetilde{H})=0$. Let $T$ $\in$ $\Mor(\delta^V,\Delta_P)$ and $b^{_T}_k=\displaystyle \sum^{n_{V}}_{i=1}T(e_i)x^{V\,\ast}_{ki}$. By equation (\ref{2.f9}) we have 
\begin{eqnarray}
\label{ecbt}
\Delta_P(b^{_T}_k)=\sum^{n_{V}}_{i=1}\Delta_P(T(e_i))\,\Delta(x^{V\,\ast}_{ki})&=& \sum^{n_{V}}_{i=1} [(T\otimes \id_H)\delta^V(e_i)]\,[(T^\l_k\otimes \id_H)\delta^V(e_i)]^\ast \nonumber
        \\
        &=&
        \sum^{n_{V}}_{i,j,l=1}T(e_j)\,x^{V\,\ast}_{kl}\otimes g^V_{ji}\,g^{V\,\ast}_{li} 
        \\
        &=&
        \sum^{n_{V}}_{j,l=1}T(e_j)\,x^{V\,\ast}_{kl}\otimes \delta_{jl}\mathbbm{1}=\sum^{n_{V}}_{j=1}T(e_j)\,x^{V\,\ast}_{kj}\otimes \mathbbm{1}=b^{_T}_k\otimes \mathbbm{1}.\nonumber
\end{eqnarray}
So, $b^{_T}_k$ $\in$ $B$ for every $k$ $\in$ $\{1,...,d_{V}\}$. Similarly, according to equation (\ref{generators}) we obtain 
\begin{equation}
\label{inverdiff}
    \left(\sum^{d_{V}}_{k=1}b^{_T}_k\,T^\l_k\right)(e_j)=\sum^{d_{V}}_{k=1}b^{_T}_k x^{V}_{kj}=\sum^{{d_{V}},n_{V}}_{k,i=1}T(e_i)\,x^{V \,\ast}_{ki}\,x^{V}_{kj}=\sum^{n_{V}}_{i=1}T(e_i)\delta_{ij}=T(e_j).
\end{equation}
We conclude that $$T=\displaystyle \sum^{d_{v}}_{k=1}b^{_T}_k\,T^\l_k.$$ 
Finally, if $\overline{b}=\displaystyle \sum^{d_{V}}_{k=1}b^{_T}_k\,\overline{e}_k$ $\in$ $B^{d_{V}}$, by equation (\ref{6.f1.11}) we get $$\widetilde{H}(\overline{b}\cdot \varrho^{V}(\mathbbm{1}))=H'(\overline{b})=\sum^{d_{v}}_{k=1}b^{_T}_k\,T^\l_k=T.$$  Therefore, $\widetilde{H}$ is a bijection. 
\end{proof}

It follows from Proposition \ref{6.1.6} that 
\begin{equation}
\label{3.f4.1}
E^V_\l:=\Mor(\delta^{V},\Delta_P)
\end{equation}
is a finitely generated projective left $B$--module for $\delta$ $\in$ $\T$ \cite{micho2}. In particular, for every $T$ $\in$ $E^V_\l$ we have 
\begin{equation}
\label{3.f5}
T=\sum^{d_{V}}_{k=1}b^{_T}_k\,T^\l_k \qquad \mbox{ with } \qquad b^{_T}_k=\sum^{n_{V}}_{i=1}T(e_i)\,x^{V\,\ast}_{ki}\;\in\;B.
\end{equation}

 Let $\delta^V$ $\in$ $\FD(\Rep_{\G})$. In accordance with \cite{woro1}, there exists $\delta^{V_1}$,..., $\delta^{V_m}$ $\in$ $\T$ such that $\delta^V \cong \bigoplus^m_{j=1} \delta^{V_j}$. Since $$\Mor(\delta^U\bigoplus \delta^W, \Delta_P)=\Mor(\delta^U,\Delta_P)\bigoplus \Mor(\delta^W,\Delta_P)$$ for any $\delta^U$, $\delta^W$ $\in$ $\T$,  it follows that 
\begin{equation}
    \label{dirsumqvb}
    E^V_\l:=\Mor(\delta^V,\Delta_P)\cong \Mor(\bigoplus^m_{j=1} \delta^{V_j},\Delta_P)=\bigoplus^m_{j=1}\Mor( \delta^{V_j},\Delta_P)=\bigoplus^m_{j=1} E^{V_j}_\l. 
\end{equation}
This shows that $E^{V}_\l$ is a finitely generated projective left $B$--module. Moreover, 
every element of $E^V_\l$ is a sum of elements of $E^{V_j}_\l$, so we can apply equation (\ref{3.f5}) to each summand of this sum and hence,  equation (\ref{3.f5}) naturally extends to every element of $E^V_\l$. In particular, the union of the left $B$--generators $\{ T^\l_k\}$ of each $E^{V_j}_\l$ forms a set of left $B$--generators of $E^V_\l$. As the reader should have already noticed, the subscript $\l$ in the notation $E^V_\l$ is to indicate that we are considering the space $\Mor(\delta^V,\Delta_P)$ as left $B$--module.

Let $\delta^V$ $\in$ $\FD(\Rep_{\G})$. Then $\delta^{\overline{V}}$ $\in$ $\FD(\Rep_{\G})$ and hence $\Mor(\delta^{\overline{V}},\Delta_P)$ is a finitely generated projective left $B$--module, where $\delta^{\overline{V}}$ denotes the conjugate corepresentation of $\delta^V$ \cite{woro1}. Endowing $\Mor(\delta^{\overline{V}},\Delta_P)$  with the right 
$B$--module structure given by $$T\cdot b:=b^\ast\, T\qquad \mbox{ where } \qquad b^\ast\, T:\overline{V}\longrightarrow P $$ is given by $(b^\ast\, T)(\overline{v})=b^\ast\, T(\overline{v})$ for all $\overline{v}$ $\in$ $\overline{V}$, it becomes a finitely generated projective right $B$--module. Furthermore, the map 
\begin{equation*}
\ast:\Mor(\delta^{V},\Delta_P) \longrightarrow \Mor(\delta^{\overline{V}},\Delta_P),\qquad
T  \longmapsto T^\ast
\end{equation*}
is a right $B$--module isomorphism, where $T^\ast$ is defined as $(T^\ast)(\overline{v}):=T(v)^\ast$ for all $\overline{v}$ $\in$ $\overline{V}$. Thus
\begin{equation}
\label{3.f5.1}
E^V_\r:=\Mor(\delta^{V},\Delta_P)
\end{equation} 
is a finitely projective right $B$--module. The subscript $\r$ in the notation $E^V_\r$ is to indicate that we are considering the space $\Mor(\delta^V,\Delta_P)$ as right $B$--module.

To obtain an equation for the right structure analogous to equation (\ref{3.f5}), note that if  $T$ $\in$ $\Mor(\delta^V,\Delta_P)$, then $T^\ast$ $\in$ $\Mor(\delta^{\overline{V}},\Delta_P)$ and therefore $$T^\ast=\displaystyle \sum b^{_{T^\ast}}_kT^\l_k.$$ Here, the maps $\{T^\l_k\}$ are the corresponding left $B$--generators of $E^{\overline{V}}_\l=\Mor(\delta^{\overline{V}},\Delta_P)$. Hence 
\begin{equation}
\label{3.f5.2}
T=\displaystyle \sum_k T^\r_k\, (b^{_{T^\ast}}_k)^\ast
\end{equation} 
with $$T^\r_k:=T^{\l\,\ast}_k \;\in\; \Mor(\delta^V,\Delta_P).$$ The superscript $\r$ in $T^\r_k$ is purely symbolic, used to indicate that these maps are right $B$--generators of $\Mor(\delta^V,\Delta_P)$.

In differential geometry, it is well--known that given a principal $G$--bundle  $\pi:P\longrightarrow B$
and a linear representation $\alpha: G\longrightarrow GL(V)$, the space of global smooth sections $\Gamma(E^V)$ of the associated vector bundle $E^V:=P\times_\alpha V$ is isomorphic to the space of $G$--equivariant maps $C^\infty_\C(P,V)^G$ as $C^\infty_\C(B)$--bimodules \cite{nodg}. Notice that the Serre--Swan theorem (\cite{serre}) allows us to identify $E^V$ with $\Gamma(E^V)$ and hence, with $C^\infty_\C(P,V)^G$. Furthermore, by the Serre--Swan theorem and in light of \cite{con, dv}, in non--commutative geometry, we define a quantum vector bundle as a finitely generated projective left or right module. 

In this way, since the {\it non--commutative geometrical} counterpart of $C^\infty_\C(P,V)^G$ is the space $\Mor(\delta^V,\Delta_P)$, then  $E^V_\l$ can be interpreted as the {\it associated left quantum vector bundle} (abbreviated ``associated left qvb") and  $E^V_\r$ can be interpreted as the {\it associated right quantum vector bundle} (abbreviated ``associated right qvb"). Our notation is analogous to {\it classical} case, but takes into account the left/right structures. 

Let $\delta^V \in \FD(\Rep_{\mathcal{G}})$. Fix a differential calculus on $\zeta=(P,B,\Delta_P)$ and fix any qpc $\omega$. Notice that $\Mor(\delta^V,\Delta_{\Hor})$ can be endowed with a $\Omega^\bullet(B)$--bimodule structure analogous to the $B$--bimodule structure of $\Mor(\delta^V,\Delta_P)$.

In this way, the map
$$\Upsilon^{-1}_{V}:\Omega^\bullet(B)\otimes_{B}E^V_\l\longrightarrow \Mor(\delta^V,\Delta_\Hor)$$ given by $\Upsilon^{-1}_{V}(\mu\otimes_{B} T)=\mu\, T$ is a left $\Omega^\bullet(B)$--module morphism.
\begin{Proposition}
    \label{innecesari3}
   The map $\Upsilon^{-1}_{V}$ is an isomorphism.
\end{Proposition}
\begin{proof}
    Let $\tau$ $\in$ $\Mor(\delta^V,\Delta_\Hor)$ and consider
    \begin{equation}
        \label{munecesito}
        \mu^\tau_k=\displaystyle \sum^{n_{V}}_{i=1} \tau(e_i)\,x^{V\,\ast}_{ki},
    \end{equation}
    where $T^\l_k(e_i)=x^V_{ij}$. If we substitute $b^{_T}_k$ with $\mu^\tau_k$ and $T$ with $\tau$ in all the calculations presented in equation (\ref{ecbt}), we now obtain $\mu^\tau_k \in \Omega^\bullet(B)$. Additionally, by
    substituting $b^{_T}_k$ with $\mu^\tau_k$ and $T$ with $\tau$ in all the calculations presented in equation (\ref{inverdiff}), we now obtain that $$\displaystyle \sum^{d_{V}}_{k=1}\mu^\tau_k\,T^\l_k=\tau.$$ In this way, we define
    \begin{equation}
\label{3.f7}
\Upsilon_{V}: \Mor(\delta^V,\Delta_\Hor)\longrightarrow \Omega^\bullet(B)\otimes_B E^V_\l,\qquad  \tau\longmapsto \sum^{d_{V}}_{k=1}\mu^\tau_k\otimes_B T^\l_k
\end{equation}
and it follows that   $\Upsilon^{-1}_{V}\circ\Upsilon_{V}=\id_{\Mor(\delta^V,\Delta_\Hor)}$. 

Let $\eta\otimes_B T$ $\in$ $\Omega^\bullet(B)\otimes_B E^V_\l$. By equation (\ref{3.f5}) we have
\begin{eqnarray*}
\Upsilon_{V}(\Upsilon^{-1}_{V}(\eta\otimes_B T))=\sum^{d_{V}}_{k=1} \mu^{\eta\,T}_k\otimes_B T^\l_k= \sum^{d_{V},n_V}_{k,i=1} \eta\,T(e_i)\,x^{V\,\ast}_{ki}\otimes_B T^\l_k  &=&
       \eta \,\sum^{d_{V}}_{k=1} b^{_T}_k\otimes_B T^\l_k
       \\
        &=&
        \eta \left( \mathbbm{1}\otimes_B \sum^{d_{V}}_{k=1} b^{_T}_k T^\l_k\right)
        \\
        &=&
        \eta  \otimes_B T
\end{eqnarray*}
and it follows that $\Upsilon_{V}\circ\Upsilon^{-1}_{V}=\id_{\Omega^\bullet(B)\otimes_B E^V_\l}$.
\end{proof}

Elements of $\Omega^\bullet(B)\otimes_{B}E^V_\l$ can be interpreted as {\it left qvb--valued differential forms} of $B$. Thus, by the first part of equation (\ref{2.f31}), $D^\omega\circ T$ $\in$ $\Mor(\delta^V,\Delta_\Hor)$ for every $T$ $\in$ $E^V_\l$ and we obtain that the linear map
\begin{equation}
\label{3.f8}
\nabla^{\omega}_{V}:E^V_\l \longrightarrow \Omega^{1}(B)\otimes_B E^V_\l,\qquad
T  \longmapsto \Upsilon_{V}\circ D^{\omega}\circ T,
\end{equation}
is a {\it quantum linear connection} on $E^V_\l$, in the sense of \cite{dv}, i.e., $\nabla^{\omega}_{V}$ satisfies the left Leibniz rule: for every $b$ $\in$ $B$ and every $T$ $\in$ $E^V_\l$ we have $$\nabla^{\omega}_V(b\,T)=db\otimes_B T+b\,\nabla^\omega_V(T).$$

Now, consider map $$\widehat{\Upsilon}^{-1}_{V}:E^V_\r\otimes_{B}\Omega^\bullet(B)\longrightarrow\Mor(\delta^V,\Delta_\Hor) $$ given by $ \widehat{\Upsilon}^{-1}_{V}(T\otimes_{B}\mu)=T \mu$. Clearly, $\widehat{\Upsilon}^{-1}_{V}$ is a right $\Omega^\bullet(B)$--module morphism,
\begin{Proposition}
    \label{innecesari4}
   The map $\widehat{\Upsilon}^{-1}_{V}$ is an isomorphism.
\end{Proposition}
\begin{proof}
     Let $\tau$ $\in$ $\Mor(\delta^V,\Delta_\Hor)$. Then $\tau^\ast$ $\in$ $\Mor(\delta^{\overline{V}},\Delta_\Hor)$, where $\delta^{\overline{V}}$ is the conjugate corepresentation of $\delta^V$ and $\tau^\ast: \overline{V}\longrightarrow \Hor^\bullet P$ is given by $\tau^\ast(\overline{v})=\tau(v)^\ast$. According to the proof of Proposition \ref{innecesari3} we obtain $$\tau^\ast=\sum_k \mu^{\tau\ast}_k\,T^\l_k.$$ Here, the maps $\{T^\l_k\}$ are the left $B$--generators of $\Mor(\delta^{\overline{V}},\Delta_P)$. Hence 
     \begin{equation}
         \label{loquesa5}
         \tau=\sum_k T^\r_k\,(\mu^{\tau\ast}_k)^\ast,
     \end{equation}
     remembering that $T^\r_k=T^{\l\,\ast}_k$.  Thus, by defining
\begin{equation}
\label{3.f7.5}
  \widehat{\Upsilon}_{V}: \Mor(\delta^V,\Delta_\Hor)\longrightarrow E^V_\r\otimes_B \Omega^\bullet(B),\qquad \tau\longmapsto \sum_k T^\r_k\otimes_B (\mu^{\tau\ast}_k)^\ast  
\end{equation}
it follows that $\widehat{\Upsilon}^{-1}_{V}\circ\widehat{\Upsilon}_{V}=\id_{\Mor(\delta^V,\Delta_\Hor)}$.  Additionally, a direct calculation using equation (\ref{3.f5.2}) as in the proof of Proposition \ref{innecesari3} shows that  $\widehat{\Upsilon}_{V}\circ\widehat{\Upsilon}^{-1}_{V}=\id_{ E^V_\r\otimes_B\Omega^\bullet(B)}$.
\end{proof}

 Elements of $E^V_\r\otimes_{B}\Omega^\bullet(B)$ can be interpreted as {\it right qvb--valued differential forms} of $B$. Hence, by the first part of equation (\ref{2.f31}), $\widehat{D}^\omega\circ T$ $\in$ $\Mor(\delta^V,\Delta_\Hor)$ for every $T$ $\in$ $E^V_\r$ and we obtain that the linear map
\begin{equation}
\label{3.f8.1}
\widehat{\nabla}^{\omega}_{V}:E^V_\r \longrightarrow E^V_\r\otimes_{B}\Omega^{1}(B),\qquad
T  \longmapsto \widehat{\Upsilon}_{V}\circ \widehat{D}^{\omega}\circ T,
\end{equation}
is a {\it quantum linear connection} on $E^V_\r$, i.e., $\widehat{\nabla}^{\omega}_{V}$ satisfies the right Leibniz rule: for every $b$ $\in$ $B$ and every $T$ $\in$ $E^V_\r$ we have $$\widehat{\nabla}^{\omega}_{V}(T\,b)=\widehat{\nabla}^{\omega}_{V}(T)\,b+T\otimes_B db .$$ The maps  $\nabla^{\omega}_{V}$ and $\widehat{\nabla}^{\omega}_{V}$ receive the name of {\it induced quantum linear connections} of $\omega$ (abbreviated ``induced qlc's"). 

It is worth mentioning that by defining 
\begin{equation}
\label{3.f10}
\sigma_{V}:=\widehat{\Upsilon}_{V}\circ \Upsilon^{-1}_{V},
\end{equation} 
we obtain $\sigma_{V}\circ \nabla^{\omega}_{V}=\widehat{\nabla}^{\omega}_{V}$, when $\omega$ is real and regular (recall that in this situation, $D^\omega=\widehat{D}^\omega$). This is the case of the theory presented in \cite{sald1}.

Extending $\nabla^{\omega}_{V}$ to the exterior covariant derivative  $$d^{\nabla^{\omega}_{V}}: \Omega^\bullet(B)\otimes_B E^V_\l\longrightarrow \Omega^\bullet(B)\otimes_B E^V_\l$$ by means of the graded Leibniz rule  
\begin{equation}
\label{3.f10.1}
d^{\nabla^{\omega}_{V}}(\mu\otimes_B T)=d\mu\otimes_B T +(-1)^k\mu \,\nabla^{\omega}_{V}(T) 
\end{equation}
for every $\mu$ $\in$ $\Omega^k(B)$ and every $T$ $\in$ $E^V_\l$, the curvature of $\nabla^{\omega}_{V}$ is defined as
\begin{equation}
\label{3.f10.2}
R^{\nabla^{\omega}_{V}}:=d^{\nabla^{\omega}_{V}}\circ \nabla^{\omega}_{V}:E^V_\l\longrightarrow \Omega^2(B)\otimes_B E^V_\l.
\end{equation}

\begin{Proposition}
    \label{innecesari6}
The following relation holds
\begin{equation}
\label{3.f10.3}
d^{\nabla^{\omega}_{V}}= \Upsilon_{V}\circ D^{\omega}\circ \Upsilon^{-1}_{V}.
\end{equation}
\end{Proposition}

\begin{proof}
    Let $\mu\otimes_B T$ $\in$ $\Omega^\bullet(B)\otimes_B E^V_\l$. Then $$(D^{\omega}\circ \Upsilon^{-1}_{V})(\mu\otimes_B T)=D^{\omega}(\mu\,T).$$ So, for all $v$ $\in$ $V$, by equations (\ref{2.f31}), (\ref{2.f32}) and the facts that $\Delta_\Hor(\Omega^\bullet(B))=\Omega^\bullet(B)\otimes \mathbbm{1}$ and $\pi(\mathbbm{1})=0$, we have 
    \begin{eqnarray*}
      (D^{\omega}(\mu\,T))(v)=D^{\omega}(\mu\,T(v))&=&D^{\omega}(\mu)\,T(v)+(-1)^k\mu\, D^{\omega}(T(v))+(-1)^k \mu\,\ell^{\omega}(\pi(\mathbbm{1}),T(v))  
        \\
        &=&
        d\mu\,T(v)+(-1)^k\mu\, D^{\omega}(T(v)).
     \end{eqnarray*}
Thus $$D^{\omega}(\mu\,T)=d\mu\,T+(-1)^k\mu\, \left(D^{\omega}\circ T\right)$$ and therefore 
\begin{eqnarray*}
      (\Upsilon_{V}\circ D^{\omega}\circ \Upsilon^{-1}_{V})(\mu\otimes_B T)= \Upsilon_{V}((D^{\omega}(\mu\,T)))&=&\Upsilon_{V}(d\mu\,T+(-1)^k\mu\, \left(D^{\omega}\circ T\right)) 
         \\
        &=&
        d\mu\otimes_B T +(-1)^k\mu\, \nabla^{\omega}_{V}(T)
        \\
        &=&
       d^{\nabla^{\omega}_{V}}(\mu\otimes_B T). 
     \end{eqnarray*}   
\end{proof}

Similarly, extending $\widehat{\nabla}^{\omega}_{V}$ to the exterior covariant derivative  $$d^{\widehat{\nabla}^{\omega}_{V}}: E^V_\r \otimes_B \Omega^\bullet(B)\longrightarrow E^V_\r \otimes_B \Omega^\bullet(B)$$ by means of the graded Leibniz rule 
\begin{equation}
\label{3.f10.4}
d^{\widehat{\nabla}^{\omega}_{V}}(T\otimes_B \mu)=\widehat{\nabla}^{\omega}_{V}(T)\, \mu +T \otimes_B d\mu 
\end{equation}
for every $\mu$ $\in$ $\Omega^\bullet(B)$ and every $T$ $\in$ $E^V_\r$, the curvature of $\widehat{\nabla}^{\omega}_{V}$ is defined as
\begin{equation}
\label{3.f10.5}
R^{\widehat{\nabla}^{\omega}_{V}}:=d^{\widehat{\nabla}^{\omega}_{V}}\circ \widehat{\nabla}^{\omega}_{V}:E^V_\r\longrightarrow E^V_\r\otimes_B \Omega^2(B).
\end{equation}

The proof of the next proposition is analogous to that of Proposition \ref{innecesari6} and will therefore be omitted. 
\begin{Proposition}
    \label{innecesari7}
The following relation holds
\begin{equation}
\label{3.f10.6}
d^{\widehat{\nabla}^{\omega}_{V}}= \widehat{\Upsilon}_{V}\circ \widehat{D}^{\omega}\circ \widehat{\Upsilon}^{-1}_{V}.
\end{equation}
\end{Proposition}

By equation (\ref{dirsumqvb}), all these constructions extend naturally to any $\delta^V$ $\in$ $\FD(\Rep_{\G})$. It is worth remarking that our formulation holds for every qpc $\omega$: it is not necessary to impose any condition on $\omega$ (as {\it reality} or {\it regularity}) to define $\nabla^\omega_V$ and $\widehat{\nabla}^\omega_V$.

The theory of connections on left/right quantum vector bundles (finitely generated projective left/right modules according to the Serre--Swan theorem) has been studied over many years, for example in \cite{lan}, and we will follow this line of research. In particular, in \cite{lan} a type of Bianchi identity is proved that all connections satisfy, but only when $B$ is commutative.

In the {\it classical} case, given a principal $G$--bundle $\pi:P\longrightarrow B$ and  a linear representation $\alpha:G\longrightarrow GL(V)$, there is a canonical isomorphism between associated vector bundle--valued differential forms of $B$ and basic differential forms of $P$ of type $\alpha$ \cite{nodg}. Moreover, this isomorphism allows to define the exterior derivative of the induced linear connection in terms of the covariant derivative of a principal connection \cite{nodg}. The definitions of $E^V_\l$, $E^V_\r$, the fact that $\Upsilon_{V}$ and $\widehat{\Upsilon}_{V}$ are isomorphisms, and equations (\ref{3.f10.3}) and (\ref{3.f10.6}) are all {\it non--commutative geometrical} counterparts of these results in differential geometry.

In accordance with \cite{br2}, $P\,\square^{H}\, V^{\#} \cong E^V_\l$ (for the natural left coaction on $V^{\#}$, the dual space of $V$). This construction is the commonly accepted one for associated qvb's. Nevertheless, we have chosen to work with $E^V_\l$ and $E^V_\r$ because, in this way, the definitions of $\nabla^{\omega}_V$ and $\widehat{\nabla}^{\omega}_V$ (and their exterior covariant derivatives) become completely analogous to their {\it classical} counterparts (in differential geometry, both connections coincide). Furthermore, they are easier to work with, as they permit explicit calculations, as the reader will verify in the remainder of this paper. In addition, by using intertwining maps, the definition of the canonical Hermitian structure appears more natural.

\subsection{The Canonical Hermitian Structure}

One of the purposes of this paper is to introduce a Hermitian structure on associated qvb's compatible with induced qlc's, and showing some of the properties of this structure. The main result will be Theorem \ref{fgs}. 

\begin{Definition}[Hermitian structures]
    \label{hs0}
    Let $(B,\cdot,\mathbbm{1},\ast)$ be a quantum space and let $M$ be a left quantum vector bundle on $B$ (a finitely generated projective left $B$--module). A Hermitian structure on $M$ is defined as a $B$--valued sesquilinear map (antilinear in the second coordinate)
    $$\langle-,-\rangle: M\times M \longrightarrow B $$ such that for all $x_1$, $x_2$ $\in$ $M$
    \begin{enumerate}
\item $\langle x_1,b\,x_2\rangle=\langle x_1,x_2\rangle \,b^\ast$;
\item $\langle x_1,x_2\rangle^\ast=\langle x_2,x_1\rangle$;
\item $\langle x_1,x_1\rangle$ $\in$ $B^+$, where $B^+$ is the pointed convex cone generated by elements of the form $\{b\,b^\ast\}$.
\end{enumerate}

In the same way, if $M$ is a right quantum vector bundle on $B$ (a finitely generated projective right $B$--module), a Hermitian structure on $M$ is a $B$--valued sesquilinear map (antilinear in the first coordinate) $$\langle-,-\rangle: M\times M \longrightarrow B $$ such that for all $x_1$, $x_2$ $\in$ $M$
    \begin{enumerate}
\item $\langle  x_1\,b,x_2\rangle=b^\ast\,\langle x_1,x_2\rangle $;
\item $\langle x_1,x_2\rangle^\ast=\langle x_2,x_1\rangle$;
\item $\langle x_1,x_1\rangle$ $\in$ $B_+$, where $B_+$ is the pointed convex cone generated by elements of the form $\{b^\ast\, b\}$.
   \end{enumerate}
\end{Definition}

For example, let us take the free left $B$--module $B^{d_V}$  (by the Serre--Swan theorem, it can be considered as a trivial left qvb). There is a canonical  Hermitian structure on $B^{d_V}$ defined by
\begin{equation}
\label{canb}
\langle-,-\rangle^{d_V}_\l:B^{d_V}\times B^{d_V} \longrightarrow B,\qquad
(\overline{b}\,,\, \overline{a})\longmapsto  \sum^{d_V}_{i=1} b_i\,a^\ast_i
\end{equation}
where $\overline{b}=(b_1,...,b_{d_V}) $ and $\overline{a}=(a_1,...,a_{d_V})$. It is worth mentioning that $\langle-,-\rangle^{d_V}_\l$ is non--degenerate, i.e., there is a Riesz representation theorem in terms of left $B$--modules, as the reader can verify in \cite{lan}. Similarly, there is a canonical  Hermitian structure on $B^{d_V}$ as a free right $B$--module (by the Serre--Swan theorem, it can be considered as a trivial right qvb) given by
\begin{equation}
\label{canb2}
\langle-,-\rangle^{d_V}_\r:B^{d_V}\times B^{d_V} \longrightarrow B,\qquad
(\overline{b}\,,\, \overline{a})\longmapsto  \sum^{d_V}_{i=1} b^\ast_i\,a_i
\end{equation} 
and it is non--degenerate, i.e., there is a Riesz representation theorem in terms of right $B$--modules. For more details, see \cite{lan}.

Let $\delta^V$ $\in$ $\T$ and let $\zeta=(P,B,\Delta_P)$ be a quantum principal $\G$--bundle. According to Proposition \ref{6.1.6}, $$E^V_\l=\Mor(\delta^V,\Delta_P)\,\cong\, B^{d_{V}}\cdot \varrho^{V}(\mathbbm{1})$$ as left $B$--modules by the map $\widetilde{H}$. 

\begin{Proposition}
    \label{caher}
    $(B^{d_V},\langle-,-\rangle^{d_V}_\l)$ induces a non--degenerate Hermitian structure on $E^V_\l$.
\end{Proposition}
\begin{proof}
    In light of \cite{lan}, the Hermitian structure $\langle-,-\rangle^{d_V}_\l$ on $B^{d_V}$ induces a non--degenerate Hermitian structure on the left $B$--submodule $B^{d_V}\cdot \varrho^V(\mathbbm{1})$  provided that the matrix 
$\varrho^V(\mathbbm{1})$ is idempotent and self--adjoint. This follows immediately from Proposition \ref{innessari2}. Thus, the induced Hermitian structure $$\langle-,-\rangle^{d_V}_\l\Big|_{B^{d_V}\cdot \varrho^V(\mathbbm{1})\times B^{d_V}\cdot \varrho^V(\mathbbm{1})}: B^{d_V}\cdot \varrho^V(\mathbbm{1}) \times B^{d_V}\cdot \varrho^V(\mathbbm{1}) \longrightarrow B$$ is non--degenerate. Finally, considering the isomorphism $\widetilde{H}$ we obtain a non--degenerate Hermitian 
structure
$$ \langle-,-\rangle_\l: E^V_\l\times E^V_\l\longrightarrow B
$$ on $E^V_\l$.
\end{proof}

 Explicitly, by equation (\ref{3.f5}) we have $\widetilde{H}^{-1}(T)=\overline{b}^{_T}\cdot \varrho^{V}(\mathbbm{1}) $ with $\overline{b}^{_T}=(b^{_T}_1,...,b^{_T}_{d_V})$ $\in$ $B^{d_V}$; and by equation (\ref{generators}) one can get $\overline{b}^{_T}\cdot \varrho^{V}(\mathbbm{1})=\overline{b}^{_T}$. Thus $\widetilde{H}^{-1}(T)=\overline{b}^{_T}$ and hence 
\begin{equation}
\label{3.f10.7}
\langle T_1,T_2\rangle_\l=\langle \overline{b}^{_{T_1}},\overline{b}^{_{T_2}}\rangle^{d_V}_\l=\sum^{d_V}_{k=1} b^{_{T_1}}_k\,(b^{_{T_2}}_k)^\ast=\sum^{n_V}_{i=1}T_1(e_i)T_2(e_i)^\ast.
\end{equation}
Furthermore, for all $b$ $\in$ $B$ we have
\begin{equation}
\label{3.f11}
\langle T_1\,b,T_2\rangle_\l=\langle T_1,T_2\,b^\ast\rangle_\l.
\end{equation}
It is worth mentioning that $\langle-,-\rangle_\l$ does not depend on the orthonormal basis $\{e_k\}^{n_{V}}_{k=1}$.

Let $\delta^{V}$ $\in$ $\FD(\Rep_{\G})$. Then there exists $\delta^{V_i}$ $\in$ $\T$ such that $\delta^{V} \, \cong \, \bigoplus^m_{i=1}\delta^{V_i}$ \cite{woro1}. Assume that $f$ is a corepresentation isomorphism between $\delta^{V}$ and $\bigoplus^m_{i=1}\delta^{V_i}$. Thus
\begin{equation}
\label{3.f12}
A_f :\bigoplus^m_{i=1}E^{V_i}_\l  \longrightarrow E^{V}_\l,\qquad
T  \longmapsto T\circ f
\end{equation}
is a left $B$--module isomorphism and its inverse is $A_{f^{-1}}$ \cite{sald1}. We can define a Hermitian structure on $E^{V}_\l$ given by

\begin{equation}
\label{3.f13}
\langle-,-\rangle_\l:E^{V}_\l\times E^{V}_\l \longrightarrow B,\qquad
(T_1\,,\, T_2)\longmapsto  \sum(T_1\circ f^{-1})(v_k)\,(T_2\circ f^{-1})(v_k)^\ast,
\end{equation}
with $\{v_k\}$ an orthonormal basis of $\bigoplus^m_{i=1}V_i$. For any unitary corepresentation morphism $f$, the previous equation agrees with the Hermitian structure on $E^V_\l$ induced by the direct sum of $(E^{V_i}_\l,\langle-,-\rangle_\l)$. So we can take equation  (\ref{3.f13}) as our definition for the Hermitian structure on $E^V_\l$ for every $\delta^{V}$ $\in$ $\FD(\Rep_{\G})$; especially since unitary corepresentation morphisms always exist. In fact, according to \cite{woro1}, $V$ decomposes into an orthogonal direct sum of subspaces $W_i$ such that $\delta^V|_{W_i}\cong \delta^{V_i}$, and $\delta^V|_{W_i}$ is unitary and irreducible. Consequently, it is sufficient to find a unitary corepresentation morphism  between  $\delta^V|_{W_i}$ and $\delta^{V_i}$. In accordance with \cite{woro1}, we have
$$\Mor(\delta^V|_{W_i},\delta^{V_i})=\{\hat{z}\, \hat{f}\mid \hat{z} \in \C  \}\qquad \mbox{ and }\qquad  \Mor(\delta^{V_i},\delta^V|_{W_i})=\{\hat{z} \hat{f}^{-1}\mid \hat{z} \in \C  \}$$
where $\hat{f}:W_i\longrightarrow V_i$ is a corepresentation isomorphism. Moreover, if $\hat{f}^\ast$ denotes the adjoint operator of $\hat{f}$, it is well--known that $\hat{f}^\ast$ $\in$ $\Mor(\delta^{V_i},\delta^V|_{W_i})$ \cite{woro1}. So, there exists $z$ $\in$ $\C$ such that $\hat{f}^\ast=z \,\hat{f}^{-1}$. Thus $\hat{f}^\ast \circ \hat{f}=z\, \id_{V_i}$ and due to the fact that $\hat{f}^\ast \circ \hat{f}$ is a positive operator and $z$ is one of its eigenvalues, we obtain that $z$ is a positive real number. Now, it is easy to check that the desired unitary corepresentation morphism is $f=\displaystyle {1\over \sqrt{z} }\, \hat{f}$.

\begin{Definition}[Canonical Hermitian structure] 
\label{lhs}
For every $\delta^V$  $\in$ $\FD(\Rep_{\G})$, we define the canonical Hermitian structure on the associated left qvb $E^V_\l$ as the sesquilinear map given by 
\begin{equation}
\langle-,-\rangle_\l:E^V_\l\times E^V_\l \longrightarrow B,\qquad
(T_1\,,\, T_2) \longmapsto  \sum^{n_{V}}_{k=1}T_1(e_k)T_2(e_k)^\ast,
\end{equation}
where $\{e_i\}^{n_{V}}_{i=1}$ is any orthonormal basis of $V$.
\end{Definition}

It is worth mentioning that, despite of the presence of the word {\it canonical} in its name, $\langle-,-\rangle_\l$ depends on the inner product $\langle-|-\rangle$ of $V$ for which $\delta^V$ is unitary, as in the {\it classical} case.

Let $\delta^V_1$, $\delta^V_2$ $\in$ $\T$. Then, the corepresentation tensor product $\delta^V_1\otimes \delta^V_2$ is finite--dimensional \cite{woro1}. Take the canonical left $B$--module isomorphism $$A_{V_1,V_2}:E^{V_1}_\l\otimes_B E^{V_2}_\l \longrightarrow E^{V_1\otimes V_2}_\l$$ defined as $$T_1\otimes_B T_2 \longmapsto A_{V_1,V_2}(T_1\otimes_B T_2),$$ where $$A_{V_1,V_2}(T_1\otimes_B T_2):V_1\otimes V_2\longrightarrow P$$ is given by $$A_{V_1,V_2}(T_1\otimes_B T_2)(v_1\otimes v_2)=T_1(v_1)T_2(v_2).$$ For more details on the morphism $A_{V_1,V_2}$, see \cite{sald1}. The proof of the following proposition is a straightforward calculation, so it will be omitted.

\begin{Proposition}
\label{tensor}
By considering the tensor product Hermitian structure on $E^{V_1}_\l \otimes_B E^{V_2}_\l$ given by  $$\langle T_1\otimes_B T_2,U_1\otimes_B U_2\rangle^\otimes_\l=\langle T_1\langle T_2,U_2\rangle_\l, U_1\rangle_\l$$ for $T_1\otimes_B T_2,$ $U_1\otimes_B U_2$ $\in$ $E^{V_1}_\l \otimes_B E^{V_2}_\l$,  the map $A_{V_1,V_2}$ is an isometry.
\end{Proposition}

Since $\langle-,-\rangle_\l$ does not depend on the orthonormal basis used to calculate it, the next proposition follows directly by another straightforward calculation; so it will be omitted as well.

\begin{Proposition}
\label{uni}
If $f: V\longrightarrow W$ is a unitary corepresentation morphism between $\delta^V$ and $\delta^W$, then the left $B$--module morphism
\begin{equation*}
A_f: E^W_\l \longrightarrow E^V_\l,\qquad T \longmapsto T\circ f
\end{equation*}
is an isometry.
\end{Proposition}

In the context of \cite{sald1}, Propositions \ref{tensor} and \ref{uni} show that the functor $\qA$ (the contravariant functor that sends $\delta^V$ to $E^V_\l$) can be
defined by endowing $\delta^V$ with an inner product that makes it a unitary corepresentation and by incorporating the canonical Hermitian structure, at least for degree $0$ morphisms. 

The introduction of Hermitian structures on associated left qvb's opens the door to the study of adjointable operators (\cite{lan})
$$\End(E^V_\l)$$ and unitary operators (\cite{lan}) $$U(E^V_\l).$$

Consider a qpb $\zeta=(P,B,\Delta_P)$. Now, let us complete   $B$ into a $C^\ast$--algebra (see Remark \ref{rema}). In light of \cite{lan}, $(B^{d_V},\langle-,-\rangle^{d_V}_\l)$ is a left Hilbert $C^\ast$--module. 

\begin{Theorem}
\label{hil}
In the previous situation, the pair $$(E^V_\l,\langle-,-\rangle_\l)$$ is a left Hilbert $C^\ast$--module  for every $\delta^V$ $\in$ $\FD(\Rep_{\G})$.
\end{Theorem}
\begin{proof}
    Equation (\ref{dirsumqvb}) shows that it is enough to prove the theorem for $\delta^V$ $\in$ $\T$. By construction, the pair $(B^{d_V},\langle-,-\rangle^{d_V}_\l)$ induces a left pre--Hilbert $C^\ast$--module structure on (see Proposition (\ref{caher})) 
    \begin{equation}
        \label{hilbermodule}
        (B^{d_V}\cdot \varrho^V(\mathbbm{1}),\langle-,-\rangle^{d_V}_\l|_{(B^{d_V}\cdot \varrho^V(\mathbbm{1}))\times (B^{d_V}\cdot \varrho^V(\mathbbm{1}))})\cong (E^V_\l,\langle-,-\rangle_\l).
    \end{equation}
    Thus, we only have to prove the completion property.
    
    It is well--known that $\End(B^{d_V})= M_{d_V}(B)$ and that every element of $\End(B^{d_V})$ is continuous \cite{lan}. Notice $\varrho^V(\mathbbm{1})$ $\in$  $M_{d_V}(B)=\End(B^{d_V})$. 
    
    Let $\{\overline{b}^i \cdot \varrho^V(\mathbbm{1})\}^\infty_{i=1}$ be a Cauchy sequence in $B^{d_V}\cdot \varrho^V(\mathbbm{1})$, where $\overline{b}^i=(b^i_1,...,b^i_{d_V})$ $\in$ $B^{d_V}$. Then $\{\overline{b}^i \cdot \varrho^V(\mathbbm{1})\}^\infty_{i=1}$ is a Cauchy sequence in $B^{d_V}$, so there exists $\overline{a}$ $\in$ $B^{d_V}$ such that  $$\overline{a}=\lim_{i\rightarrow \infty}(\overline{b}_i\cdot \varrho^V(\mathbbm{1}))=\lim_{i\rightarrow \infty}(\overline{b}_i\cdot \varrho^V(\mathbbm{1})\cdot \varrho^V(\mathbbm{1}))=\lim_{i\rightarrow \infty}(\overline{b}_i\cdot \varrho^V(\mathbbm{1}))\cdot \varrho^V(\mathbbm{1})=\overline{a}\cdot \varrho^V(\mathbbm{1}),$$ where we have used the fact that $\varrho^V(\mathbbm{1})$ is an idempotent element. The last equation implies that $\overline{a}$ $\in$ $B^{d_V}\cdot \varrho^V(\mathbbm{1})$ and hence, the left--hand side of the equation (\ref{hilbermodule}) has structure of left Hilbert $C^\ast$--module. Finally, by considering the isomorphism $\widetilde{H}$ the theorem follows. 
\end{proof}

 Theorem (\ref{hil}) is particularly important because it allows one to apply the full theory of left Hilbert $C^\ast$--modules to associated left qvb's.

All the theory developed in this subsection holds for associated right qvb's $E^V_\r$. In particular, we have

\begin{Definition}[Canonical Hermitian structure] 
\label{rhs}
For every $\delta^V$  $\in$ $\FD(\Rep_{\G})$, we define the canonical Hermitian structure on the associated right qvb $E^V_\r$ as the sesquilinear map (now antilinear in the first coordinate) given by 
\begin{equation*}
\langle-,-\rangle_\r:E^V_\r\times E^V_\r \longrightarrow B,\qquad
(T_1\,,\, T_2)\; \longmapsto  \sum^{n_{V}}_{k=1}T_1(e_k)^\ast\, T_2(e_k),
\end{equation*}
where $\{e_i\}^{n_{V}}_{i=1}$ is any orthonormal basis of $V$.
\end{Definition}
\noindent and
\begin{Theorem}
\label{hil1}
Completing $B$ into a $C^\ast$--algebra, the pair  $$(E^V_\r,\langle-,-\rangle_\r)$$ is a right Hilbert $C^\ast$--module for every $\delta^V$ $\in$ $\FD(\Rep_{\G})$.
\end{Theorem}

Before continuing, it is important to note that our definition of the canonical 
Hermitian structure on associated left/right qvb's is based on the general theory presented in \cite{lan} and on the dual observation that, in differential geometry, given a principal $G$--bundle $\pi: P\longrightarrow B$, one can define a Hermitian structure on the associated vector bundle $E^V$ (arising from a unitary $G$--representation on $V$) by using the inner product of $V$. In addition, in  the {\it classical} case, the space of smooth  sections of an associated vector bundle can be completed to a Hilbert $C^\ast$--bimodule \cite{lan}.  Theorems (\ref{hil}), (\ref{hil1}) are the {\it non--commutative geometrical} counterparts of this result.

By taking a differential calculus on the qpb $\zeta=(P,B,\Delta_P)$, the canonical Hermitian structure on $E^V_\l$ can be extended to 
\begin{equation}
\label{algo11}
\begin{aligned}
\langle-,-\rangle_\l:(\Omega^\bullet(B)\otimes_B E^V_\l)\times (\Omega^\bullet(B)\otimes_B E^V_\l)\longrightarrow  \Omega^\bullet(B)
\end{aligned}
\end{equation}
by means of \cite{lan} $$\langle \mu_1\otimes_B T_1,\mu_2\otimes_B T_2\rangle_\l=\mu_1 \,\langle T_1,T_2\rangle_\l\,\mu^\ast_2.$$ 
Similarly, the canonical Hermitian structure on $E^V_\r$ can be extended to  
\begin{equation}
\label{algo12}
\begin{aligned}
\langle-,-\rangle_\r:(E^V_\r\otimes_B \Omega^\bullet(B))\times (E^V_\r\otimes_B \Omega^\bullet(B))\longrightarrow  \Omega^\bullet(B)
\end{aligned}
\end{equation}
by means of \cite{lan} $$ \langle  T_1\otimes_B \mu_1, T_2\otimes_B \mu_2 \rangle_\r=\mu^\ast_1 \,\langle T_1,T_2\rangle_\r\,\mu_2.$$ 

\begin{Definition}[Compatible Quantum Linear Connections]
    \label{definne}
    Let $\nabla: E^V_\l\longrightarrow \Omega^1(B)\otimes_B E^V_\l$ be a qlc, i.e., $\nabla$ is a linear map that satisfies the left Leibniz rule. We say that $\nabla$ is {\it compatible} with the Hermitian structure $\langle-,-\rangle_\l$, or simply that it is a Hermitian qlc if
    \begin{equation}
\label{3.f18}
\langle \nabla(T_1),T_2 \rangle_\l+\langle T_1,\nabla(T_2)\rangle_\l=d\langle T_1,T_2\rangle_\l 
\end{equation}
for all $T_1$, $T_2$ $\in$ $E^V_\l$. 

Similarly, let $\nabla: E^V_\r\longrightarrow E^V_\r\otimes_B\Omega^1(B) $ be a qlc, i.e., $\nabla$ is a linear map that satisfies the right Leibniz rule. We say that $\nabla$ is {\it compatible} with the Hermitian structure $\langle-,-\rangle_\r$, or simply that it is a Hermitian qlc if
    \begin{equation}
\label{3.f18.1}
\langle \nabla(T_1),T_2 \rangle_\r+\langle T_1,\nabla(T_2)\rangle_\r=d\langle T_1,T_2\rangle_\r
\end{equation}
for all $T_1$, $T_2$ $\in$ $E^V_\r$.
\end{Definition}

In the {\it classical} case, if we take a principal connection, then the induced linear connection on an associated vector bundle $E^V$ for a unitary representation on $V$ is compatible with the Hermitian structure induced by the inner product on $V$. The following theorem reflects this important result in non--commutative geometry.

\begin{Theorem}
\label{fgs}
Let $\zeta=(P,B,\Delta_P)$ be a qpb and choose a differential calculus on it. Let $\omega$ be a real qpc (which always exists \cite{micho2}) and $\delta^V$ $\in$ $\FD(\Rep_{\G})$. Then, the induced qlc is compatible with the canonical Hermitian structure. 
\end{Theorem}
\begin{proof}
Let $\delta^V$ $\in$ $\T$. According to equation (\ref{3.f5}), every $T$ $\in$ $E^V_\l$ can be expressed as
$$ T=\sum^{d_{V}}_{k=1}b^{_T}_k\,T^\l_k \quad \mbox{ with } \quad b^{_T}_k=\sum^{n_{V}}_{i=1}T(e_i)\,x^{V\,\ast}_{ki}.$$
Moreover,  by equations (\ref{3.f7}), (\ref{3.f8}) we get
$$\nabla^\omega_{V}(T)=\sum^{d_V}_{k=1}\mu^{D^{\omega}\circ T}_k\otimes_B T^\l_k \qquad \mbox{ with }\qquad  \mu^{D^{\omega}\circ T}_k=\sum^{n_V}_{i=1}D^{\omega}(T(e_i))x^{V\ast}_{ki}.$$

\noindent Then, by equation (\ref{generators}), for all $T_1$, $T_2$ $\in$ $E^V_\l$ we have
\begin{eqnarray*}
\langle \nabla^{\omega}_{V} (T_1),T_2 \rangle_\l+\langle T_1, \nabla^{\omega}_{V} (T_2)\rangle_\l  &= &  \sum^{d_{V},n_{V}}_{k,i=1}\mu^{D^{\omega}\circ T_1}_k \,T^\l_k(e_i)\,T_2(e_i)^\ast
  \\
  &+ &
T_1(e_i)\,T^\l_k(e_i)^\ast\,(\mu^{D^{\omega}\circ T_2}_k)^\ast
  \\
  &= &
\sum^{d_{V},n_{V}}_{k,i,j=1}D^{\omega}(T_1(e_j))\,x^{V\,\ast}_{kj}x^{V}_{ki}\,T_2(e_i)^\ast
  \\
  &+ &
T_1(e_i)\,x^{V\,\ast}_{ki}x^{V}_{kj}\,D^{\omega}(T_2(e_j))^\ast
  \\
  &= &
  \sum^{n_{V}}_{i=1}D^{\omega}(T_1(e_i))\,T_2(e_i)^\ast+T_1(e_i)\,D^{\omega}(T_2(e_i))^\ast.
\end{eqnarray*}  
Notice that for every $T$ $\in$ $E^V_\l$ $$\Delta_P(T(e_i))=\displaystyle \sum^{n_V}_{j=1}T(e_j)\otimes g^V_{ji},$$  where the elements $g^V_{ji}$ $\in$ $H$ are given in Theorem \ref{rep}. Since $\omega$ is real, in accordance with equation (\ref{2.f32.5}) and using the fact that $S(g^V_{ji})=g^{V\ast}_{ij}$ (\cite{micho1}), we obtain $$D^{\omega}(T_2(e_i))^\ast= D^{\omega}(T_2(e_i)^\ast)+ \sum^{n_V}_{j=1} \ell^{\omega}(\pi(S(g^V_{ji})^\ast),T_2(e_j)^\ast)=D^{\omega}(T_2(e_i)^\ast)+ \sum^{n_V}_{j=1} \ell^{\omega}(\pi(g^{V}_{ij}),T_2(e_j)^\ast).$$ Thus
\begin{eqnarray*}
 \sum^{n_{V}}_{i=1}D^{\omega}(T_1(e_i))\,T_2(e_i)^\ast+T_1(e_i)\,D^{\omega}(T_2(e_i))^\ast
  &= &
  \sum^{n_{V}}_{i=1}D^{\omega}(T_1(e_i))\,T_2(e_i)^\ast\\
  &+ &
  T_1(e_i) D^{\omega}(T_2(e_i)^\ast)
  \\
  &+ &
  \sum^{n_V}_{i,j=1} T_1(e_i) \ell^{\omega}(\pi(g^{V}_{ij}),T_2(e_j)^\ast).
\end{eqnarray*}
Now, according to equation (\ref{2.f32}) we have
\begin{eqnarray*}
D^\omega(T_1(e_i)\, T_2(e_i)^\ast)&= & D^\omega(T_1(e_i))\,T_2(e_i)^\ast+ T_1(e_i)\,D^\omega(T_2(e_i)^\ast)
\\
&+ &
 \sum^{n_V}_{j=1}T_1(e_j)\,\ell^{\omega}(\pi(g^{V}_{ji}),T_2(e_i)^\ast); 
\end{eqnarray*}
which implies

\begin{eqnarray*}
\sum^{n_{V}}_{i=1}D^{\omega}(T_1(e_i))\,T_2(e_i)^\ast\;+\; T_1(e_i) D^{\omega}(T_2(e_i)^\ast) &+ &\sum^{n_V}_{i,j=1} T_1(e_i) \ell^{\omega}(\pi(g^{V}_{ij}),T_2(e_j)^\ast)
    \\
  &= &
  \sum^{n_{V}}_{i=1} D^{\omega}(T_1(e_i)\,T_2(e_i)^\ast)
    \\
  &- &  
 \sum^{n_{V}}_{i,j=1}T_1(e_j)\,\ell^{\omega}(\pi(g^{V}_{ji}),T_2(e_i)^\ast) 
  \\
  &+ &
 \sum^{n_{V}}_{i,j=1}T_1(e_i)\,\ell^{\omega}(\pi(g^{V}_{ij})),T_2(e_j)^\ast)  
     \\
  &= &
  \sum^{n_{V}}_{i=1} D^{\omega}(T_1(e_i)\,T_2(e_i)^\ast)\; =\; D^{\omega} \langle T_1,T_2\rangle_\l
\end{eqnarray*}
Finally, by the second part of equation (\ref{2.f31}) we get $D^{\omega}\langle T_1,T_2\rangle_\l=d \langle T_1,T_2\rangle_\l$ and therefore $$\langle \nabla^{\omega}_{V} (T_1),T_2 \rangle_\l+\langle T_1, \nabla^{\omega}_{V} (T_2)\rangle_\l=d \langle T_1,T_2\rangle_\l. $$ Let $\delta^V$ $\in$ $\FD(\Rep_{\G})$. Then $\delta^V\cong \bigoplus^m_{i=1}\delta^{V_i}$ for some $m$ $\in$ $\N$ with $\delta^{V_i}$ $\in$ $\T$ (\cite{woro1}) and by equation (\ref{dirsumqvb}) we have  $E^V_\l\cong \bigoplus^m_{i=1}E^{V_i}_\l$. Since the theorem holds for $E^{V_i}_\l$ and the canonical Hermitian structure of $E^V_\l$ agrees with the direct sum of the Hermitian structures of $E^{V_i}_\l$, it follows that the theorem holds for $E^V_\l$.

Let $\delta^V$ $\in$ $\FD(\Rep_{\G})$. For every $T$ $\in$ $E^V_\r$  we have (see equation (\ref{3.f5.2}))
$$ T=\sum_{k} T^\r_k\,(b^{_{T^\ast}}_k)^\ast \quad \mbox{ with } \quad 
T^\r_k=T^{\l\,\ast}_k \quad \mbox{ and } \quad (b^{_{T^\ast}}_k)^\ast= \left( \sum^{n_{V}}_{i=1}T(e_i)^\ast\,x^{V\,\ast}_{ki} \right)^\ast, $$ where the maps $T^{\l}_k$ are the corresponding left $B$--generators of $E^{\overline{V}}_\l=\Mor(\delta^{\overline{V}},\Delta_P)$. In addition, by  equations (\ref{3.f7.5}), (\ref{3.f8.1}) we get
$$\widehat{\nabla}^\omega_{V}(T)=\sum_{k} T^\r_k\otimes_B (\mu^{(\widehat{D}^{\omega}\circ T)^\ast}_k)^\ast \quad \mbox{ with } \quad  (\mu^{(\widehat{D}^{\omega}\circ T)^\ast}_k)^\ast=\left(\sum^{n_V}_{i=1}(\widehat{D}^{\omega}(T(e_i)))^\ast x^{V\ast}_{ki}\right)^\ast.
$$ In this way, for all $T_1$, $T_2$ $\in$ $E^V_\r$ we obtain
\begin{eqnarray*}
\langle \widehat{\nabla}^{\omega}_{V} (T_1),T_2 \rangle_\r+\langle T_1, \widehat{\nabla}^{\omega}_{V} (T_2)\rangle_\r &= & \sum_{k,i}\mu^{(\widehat{D}^{\omega}\circ T_1)^\ast}_k \,T^{\r}_k(e_i)^\ast\,T_2(e_i) + T_1(e_i)^\ast\,T^\r_k(e_i)\,(\mu^{(\widehat{D}^{\omega}\circ T_2)^\ast}_k)^\ast
\\
  &= &
\sum_{k,i}\mu^{(\widehat{D}^{\omega}\circ T_1)^\ast}_k \,T^{\l}_k(e_i)\,T_2(e_i) + T_1(e_i)^\ast\,T^\l_k(e_i)^\ast\,(\mu^{(\widehat{D}^{\omega}\circ T_2)^\ast}_k)^\ast
  \\
  &= &
\sum_{k,i}\mu^{{D}^{\omega}\circ T^\ast_1}_k \,T^{\l}_k(e_i)\,T_2(e_i) + T_1(e_i)^\ast\,T^\l_k(e_i)^\ast\,(\mu^{D^{\omega}\circ T^\ast_2}_k)^\ast,
\end{eqnarray*}  
where we have used that $\widehat{D}^{\omega}=\ast \circ D^\omega\circ \ast$ (see equation (\ref{2.f30.1})). Notice that $T_3:=T^\ast_1$, $T_4:=T^\ast_2$ $\in$ $E^{\overline{V}}_\l$; so according to the first part of this proof
$$\sum_{k,i}\mu^{{D}^{\omega}\circ T^\ast_1}_k \,T^{\l}_k(e_i)\,T_2(e_i) + T_1(e_i)^\ast\,T^\l_k(e_i)^\ast\,(\mu^{D^{\omega}\circ T^\ast_2}_k)^\ast= $$ $$\sum_{k,i}\mu^{{D}^{\omega}\circ T_3}_k \,T^{\l}_k(e_i)\,T_4(e_i)^\ast + T_3(e_i)\,T^\l_k(e_i)^\ast\,(\mu^{D^{\omega}\circ T_4}_k)^\ast =$$ $$\langle \nabla^{\omega}_{\overline{V}} (T_3),T_4 \rangle_\l+\langle T_3, \nabla^{\omega}_{\overline{V}} (T_4)\rangle_\l=d \langle T_3,T_4\rangle_\l=d \langle T_1,T_2\rangle_\r$$ and therefore $$\langle \widehat{\nabla}^{\omega}_{V} (T_1),T_2 \rangle_\r+\langle T_1, \widehat{\nabla}^{\omega}_{V} (T_2)\rangle_\r=d \langle T_1,T_2\rangle_\r. $$
\end{proof}

The last theorem is the core of the paper. First, it establishes the relationship between induced quantum linear connections and canonical Hermitian structures. Second, although Theorem \ref{fgs} holds only for real qpc's, it will enable us to prove the existence of formally adjoint operators for the exterior covariant derivatives $d^{\nabla^\omega_V}$, $d^{\widehat{\nabla}^\omega_V}$ for every qpc $\omega$ (without the {\it reality} condition). With these formally adjoint operators, one can develop a field theory and a Yang--Mills theory starting from a qpb, exactly as in differential geometry \cite{sald2, sald3, sald4}. 

It is worth mentioning that even when all right $B$--structures are induced by the $\ast$ operation, in general, $\nabla^\omega_{\overline{V}}=\Upsilon_{\overline{V}}\circ D^\omega$ is not equivalent to $\widehat{\nabla}^\omega_{V}=\widehat{\Upsilon}_V\circ \widehat{D}^\omega$ due to equation (\ref{covd.1}).

\section{The Quantum Gauge Group and its Action on Quantum Connections}

In addition to the aim mentioned earlier of studying associated qvb's and its induced structures, the other purpose of this paper is to present an {\it ad hoc} definition of the quantum gauge group for a given qpb with a differential calculus and studying its action on the space of qpc's, as well as on the space of induced qlc's. We address all of these topics in this section. We are particularly interested in the action of the quantum gauge group on qpc's because, as was discussed at the introduction, the final goal of our research is to develop a Yang--Mills theory in Durdevich's framework, and for that it is necessary to classify qpc's via gauge transformations. 

Our definition of the quantum gauge group is inspired by the one presented in \cite{br}, and for that we need to discuss first the {\it quantum translation map} \cite{br, micho4, micho5}.

\subsection{The Quantum Translation Map and The Quantum Gauge Group}

Let $\zeta=(P,B,\Delta_P)$ be a qpb. The  map $\beta$ of equation (\ref{porqueyo.3}) can be used to define the linear isomorphism
\begin{equation*}
\widetilde{\beta}:P\otimes_B P\longrightarrow P\otimes H
\end{equation*}
given by (\cite{stheve}) $$\widetilde{\beta}(x\otimes_B y)=\beta(x\otimes y)=(x\otimes \mathbbm{1})\cdot \Delta_P(y).$$ The degree zero quantum translation map is defined as
\begin{equation}
\qtrs: H \longrightarrow P\otimes_B P
\end{equation}
such that  $$\qtrs(g)=\widetilde{\beta}^{-1}(\mathbbm{1}\otimes g).$$  Explicitly, by taking the linear basis $\{g^{V}_{ij}\}_{\delta^V,i,j}$, we have (see Theorem \ref{rep} and Remark \ref{rema})
\begin{equation}
    \label{qtrs0}
    \qtrs(g^{V}_{ij})=\sum^{d_{V}}_{k=1} T^\l_k(e_i)^\ast\otimes_B T^\l_k(e_j)=\sum^{d_{V}}_{k=1} x^{V\,\ast}_{ki}\otimes_B x^{V}_{kj}.
\end{equation}

The map $\qtrs$ can be extended to 
\begin{equation}
\label{5.f1.3}
\widetilde{\qtrs}:P\otimes H\longrightarrow P\otimes_B P
\end{equation}
by means of  $$\widetilde{\qtrs}(x\otimes g^{V}_{ij})=x\,\qtrs(g^{V}_{ij})=\displaystyle\sum^{d_{\alpha}}_{k=1}x\, x^{V\,\ast}_{ki}\otimes_B x^{V}_{kj}$$  and a direct calculation shows that  $\widetilde{\beta}$ and $\widetilde{\qtrs}$ are mutually inverse (\cite{stheve}).

Let $g$ $\in$ $H$. Throughout the various computations of this paper, we shall use the symbolic notation
\begin{equation}
\label{5.f1.4}
\qtrs(g)=\sum_{i,j}p_i\otimes_B p'_j,=:[g]_1\otimes_B [g]_2,
\end{equation}
the sum is understood.

Now, we assume that $\zeta=(P,B,\Delta)$ is endowed with a differential calculus. In this situation, $\widetilde{\beta}$ has a natural extension to (taking the tensor product of graded differential $\ast$--algebras)
\begin{equation}
\label{6.f1.5}
\widetilde{\beta}:\Omega^\bullet(P)\otimes_{\Omega^\bullet(B)}\Omega^\bullet(P)\longrightarrow \Omega^\bullet(P)\otimes \Gamma^\wedge
\end{equation}
given by $$\widetilde{\beta}(w_1\otimes_{\Omega^\bullet(B)}w_2)=(w_1\otimes\mathbbm{1})\cdot \Delta_{\Omega^\bullet(P)}(w_2).$$  According to \cite{micho4}, this map is bijective. 

On the other hand, for a real qpc $\omega$ (which always exists \cite{micho2}) and in accordance with \cite{micho4}, we can extend $\qtrs$ to
\begin{equation}
\label{6.f1.6}
\qtrs:\mathfrak{qg}^\#\longrightarrow \left(\Omega^\bullet(P)\otimes_{\Omega^\bullet(B)}\Omega^\bullet(P)\right)^{1}:=(\Omega^1(P)\otimes_{\Omega^\bullet(B)} P)\oplus (P\otimes_{\Omega^\bullet(B)}\Omega^1(P))
\end{equation}
by means of 
\begin{eqnarray}
\label{6.f1.6.1}
\qtrs(\theta):=\mathbbm{1}\otimes_{\Omega^\bullet(B)}\omega(\theta)-(m_\Omega\otimes_{\Omega^\bullet(B)}\id_{P})(\omega\otimes \qtrs)\ad(\theta), 
\end{eqnarray}
where $m_\Omega:\Omega^\bullet(P)\otimes\Omega^\bullet(P)\longrightarrow \Omega^\bullet(P)$  is the product map and $\ad(\theta)=\theta^{(0)}\otimes \theta^{(1)}$ $\in$ $\mathfrak{qg}^\#\otimes H$. As before, we are going to use the symbolic notation $$\qtrs(\theta)=:[\theta]_1\otimes_{\Omega^\bullet(B)}[\theta]_2,$$ the sum is understood. According to \cite{micho4}, $\qtrs$ can be extended to
\begin{equation*}
    \label{qtrsgamma}
    \qtrs:\Gamma\longrightarrow \left(\Omega^\bullet(P)\otimes_{\Omega^\bullet(B)}\Omega^\bullet(P)\right)^{1} 
\end{equation*}
by means of
\begin{equation}
\label{6.f1.7}
\qtrs(g\,\theta)=[\theta]_1\,\qtrs(g)\, [\theta]_2,\qquad \qtrs(\theta\,g):= [g]_1\,\qtrs(\theta)\,[g]_2,
\end{equation}
    where $g$ $\in$ $H$, $\theta$ $\in$ $\mathfrak{qg}^\#$.  It is worth mentioning that although the definition 
of $\qtrs$ appears to depend on the choice of the real qpc $\omega$, the uniqueness of the inverse ensures that $\qtrs$ is independent of this choice. Moreover,  since $\widetilde{\beta}$ commutes with the corresponding differential maps, it follows that (\cite{micho4})
 
\begin{Proposition}
\label{6.1.1}
The following identity holds  $$\qtrs\circ d=d_{\otimes^\bullet} \circ \qtrs,$$  where $d_{\otimes^\bullet}$ is the canonical differential map of  $\Omega^\bullet(P)\otimes_{\Omega^\bullet(B)}\Omega^\bullet(P)$.
\end{Proposition}

In light of \cite{micho4}, the quantum translation map can be extended to 
\begin{equation}
\label{6.f1.10}
\qtrs:\Gamma^\wedge\longrightarrow \Omega^\bullet(P)\otimes_{\Omega^\bullet(B)}\Omega^\bullet(P)
\end{equation}
by means of $$\qtrs(\vartheta\,\upsilon)= (-1)^{\partial\vartheta\,\partial[\upsilon]_1} \,[\upsilon]_1\,\qtrs(\vartheta)\,[\upsilon]_2.$$
In addition, according to \cite{micho4}, the corresponding  map
\begin{equation}
\label{6.f1.12}
\widetilde{\qtrs}:\Omega^\bullet(P)\otimes \Gamma^\wedge\longrightarrow \Omega^\bullet(P)\otimes_{\Omega^\bullet(B)}\Omega^\bullet(P)
\end{equation}
defined by  $$\widetilde{\qtrs}(w_1\otimes \vartheta)=(w_1\otimes \mathbbm{1})\cdot \qtrs(\vartheta)$$ is the inverse of $\widetilde{\beta}.$  In accordance with \cite{micho5}, the following properties hold:
   \begin{enumerate}  
\item For all $\vartheta$ $\in$ $\Gamma^\wedge$, we have  $$[\vartheta]_1\,[\vartheta]_2=\epsilon(\vartheta)\mathbbm{1}.$$ 
\item $(\id_{\Omega^\bullet(P)}\otimes_{\Omega^\bullet(B)}\Delta_{\Omega^\bullet(P)})\circ\qtrs=(\qtrs\otimes \id_{\Gamma^\wedge})\circ \Delta$.
\item $(\Delta_{\Omega^\bullet(P)}\otimes_{\Omega^\bullet(B)} \id_{\Omega^\bullet(P)})\circ\qtrs=(\sigma\otimes_{\Omega^\bullet(B)} \id_{\Omega^\bullet(P)}) \circ(S\otimes \qtrs)\circ \Delta$, where $$\sigma:\Gamma^\wedge \otimes \Omega^\bullet(P) \longrightarrow \Omega^\bullet(P)\otimes\Gamma^\wedge $$ is the canonical graded twist map, i.e.,  $$\sigma(\vartheta\otimes w)=(-1)^{kl}\,w\otimes \vartheta$$  if $w$ $\in$ $\Omega^k(P)$ and $\vartheta$ $\in$ $\Gamma^{\wedge l}$.
\item For all $\mu$ $\in$ $\Omega^k(B)$, $\vartheta$ $\in$ $\Gamma^{\wedge l}$, we obtain $$\mu\,\qtrs(\vartheta)=(-1)^{lk}\qtrs(\vartheta)\,\mu.$$ 
\end{enumerate}

Let   $$\f_1,\,\f_2:\Gamma^\wedge\longrightarrow \Omega^\bullet(P)$$ be two graded linear maps. The convolution product of $\f_1$ with $\f_2$ is defined by 
\begin{equation}
    \label{convolutionproduct}
    \f_1\widetilde{\ast} \,\f_2=m_\Omega \circ (\f_1\otimes \f_2)\circ \Delta:\Gamma^\wedge \longrightarrow \Omega^\bullet(P).
\end{equation}

Henceforth, we will just consider graded maps $\f$ such that 
\begin{equation}
    \label{4.f3}
    \f(\mathbbm{1})=\mathbbm{1} \qquad \mbox{and} \qquad (\f\otimes \id_{\Gamma^\wedge})\circ \Ad=\Delta_{\Omega^\bullet(P)}\circ \f,
\end{equation}
where $\Ad: \Gamma^\wedge\longrightarrow \Gamma^\wedge \otimes \Gamma^\wedge$ is the extension of the right $\G$--coaction $\Ad: H\longrightarrow H\otimes H$ (see equation (\ref{2.f11})). We say that $\f$ is a {\it convolution invertible map} if there exists a graded linear map   $$\f^{-1}: \Gamma^\wedge\longrightarrow \Omega^\bullet(P)$$
such that 
\begin{equation}
    \label{loquencesito5}
    \f\,\widetilde{\ast}\, \f^{-1}=\f^{-1}\,\widetilde{\ast}\, \f =\mathbbm{1}\epsilon.
\end{equation}
A direct calculation proves that the set of all convolution invertible maps $$ \{ \f: \Gamma^\wedge\longrightarrow \Omega^\bullet(P)\} $$ is a group with respect to the convolution product. 

\begin{Proposition}
    \label{4.1}
    There exists a group isomorphism between the group of all convolution invertible maps $\{ \f: \Gamma^\wedge\longrightarrow \Omega^\bullet(P)\} $ and the group of all graded left $\Omega^\bullet(B)$--module isomorphisms  $$\F:\Omega^\bullet(P)\longrightarrow \Omega^\bullet(P)$$ that satisfy 
    \begin{equation}
        \label{4.f4}
        \F(\mathbbm{1})=\mathbbm{1} \qquad \mbox{ and }\qquad (\F\otimes \id_{\Gamma^\wedge})\circ \Delta_{\Omega^\bullet(P)}=\Delta_{\Omega^\bullet(P)}\circ \F.
    \end{equation}
    Here, we are considering the group product  $(\F_1\cdot \F_2)(w)=\F_2(\F_1(w)).$  
\end{Proposition}

\begin{proof}
For a map $\F$, consider 
\begin{equation}
\label{4.f5}
\f_\F:=m_{\Omega^\bullet}\circ (\id_{\Omega^\bullet(P)}\otimes_{\Omega^\bullet(B)}\F)\circ \qtrs:\Gamma^\wedge\longrightarrow \Omega^\bullet(P),
\end{equation}
where $m_{\Omega^\bullet}:\Omega^\bullet(P)\otimes_{\Omega^\bullet(B)}\Omega^\bullet(P)\longrightarrow \Omega^\bullet(P)$ is the product map; and for $\f$ define
\begin{equation}
\label{4.f5.56}
\F_\f:=m_{\Omega}\circ (\id_{\Omega^\bullet(P)}\otimes \f)\circ \Delta_{\Omega^\bullet (P)}:\Omega^\bullet(P)\longrightarrow \Omega^\bullet(P).
\end{equation}
Now, the proof is completely analogous to the one presented in \cite{br}, but considering the graded--differential $\ast$--Hopf algebra $\Gamma^{\wedge\,\infty}$ instead of $H^\infty$, as the reader can verify in \cite{appendix}. 
\end{proof}

In this way, we define 

\begin{Definition}[The quantum gauge group]
    \label{qgg}
     Let $\zeta=(P,B,\Delta_P)$ be a quantum principal $\G$--bundle over $B$ with a differential calculus. We define the quantum gauge group $$\qGG$$ as the group of all graded left $\Omega^\bullet(B)$--module isomorphisms $$\F:\Omega^\bullet(P)\longrightarrow \Omega^\bullet(P)$$ that satisfy equation (\ref{4.f4}). Elements of $\qGG$ are referred to as quantum gauge transformations (abbreviated ``qgt's").
\end{Definition}

The {\it a priori} motivation for our definition of $\qGG$ is the fact that
in differential geometry, gauge transformations are vertical principal bundle automorphisms. In this way, Definition \ref{qgg} was derived by {\it dualizing} this {\it classical} fact, while ensuring that $\qGG$ is defined for every degree in the most general manner, without imposing any unnecessary condition.

Moreover, in the \emph{classical} case, gauge transformations can also be described as elements of $C^\infty(P,G)^G $, where the action of $G$ on $G$ is the adjoint one  \cite{nodg}. The {\it dual} result of this {\it classical} fact is consistent with Proposition \ref{4.1}. Additionally, \cite{sald2, sald3, sald4} provides \emph{a posteriori} justification for Definition \ref{qgg} based on the orbits of Yang--Mills qpc's. We will discuss this further in the final section.  

It is worth remarking that $\F$ and $\f$ are only graded linear maps, so, in general, they do not commute with the corresponding differentials. 

\begin{Remark}
   \label{remahor} 
   By our definition of $\Hor^\bullet P$ (see equation (\ref{2.f18})), it follows that  $$\F(\Hor^\bullet P)=\Hor^\bullet P$$  for every qgt $\F$. In particular  $$\Delta_\Hor(\F(\varphi))=\F(\varphi^{(0)})\otimes \varphi^{(1)}$$  with $\Delta_\Hor(\varphi)=\varphi^{(0)}\otimes \varphi^{(1)}$ for all $\varphi$ $\in$ $\Hor^\bullet P$.
\end{Remark} 

Let $\zeta=(P,B,\Delta_P)$ be a quantum principal $\G$--bundle with a differential calculus. The set of all characters of $H$ 
\begin{equation}
\label{4.f7}
H_{cl}:=\{ \chi:H:\longrightarrow \C\mid \chi \mbox{ is a character} \}
\end{equation}
has a group structure with multiplication  $\chi_1\ast\chi_2:= (\chi_1\otimes\chi_2)\circ \Delta,$ unity  $\epsilon,$  and inverses defined by $\chi^{-1}:=\chi \circ S$ (\cite{woro1, micho1}). In agreement with the Gelfand--Naimark theorem, this group can be interpreted as the group of all {\it classical} points of $\G$ and it is isomorphic to a compact subgroup of $U(n)$ for some $n$ $\in$ $\N$, according to \cite{micho1} (remember that $\G$ is a compact matrix quantum group). Every character $\chi$ can be extended to 
\begin{equation}
\label{4.f8}
\chi:\Gamma^\wedge\longrightarrow \C
\end{equation} 
by $$\chi|_{H}:=\chi\qquad \mbox{ and }\qquad \chi|_{\Gamma^{\wedge k}}:=0\qquad \mbox{ for } \qquad k\geq 1.$$ Consider $$\F_{\chi}:=(\id_{\Omega^\bullet(P)}\otimes \chi)\circ \Delta_{\Omega^\bullet(P)}: \Omega^\bullet(P) \longrightarrow \Omega^\bullet(P).$$ This map is a graded differential $\ast$--algebra isomorphism with inverse  $\F^{-1}_{\chi}:=\F_{\chi^{-1}}$. Indeed, 
\begin{eqnarray*}
\F_{\chi^{-1}}\circ \F_{\chi}&=&  (\id_{\Omega^\bullet(P)}\otimes \chi^{-1})\circ \Delta_{\Omega^\bullet(P)} \circ (\id_{\Omega^\bullet(P)}\otimes \chi)\circ \Delta_{\Omega^\bullet(P)} 
 \\
 &=& 
(\id_{\Omega^\bullet(P)}\otimes \chi^{-1}\otimes \chi)\circ (\Delta_{\Omega^\bullet(P)} \otimes \id_{\Gamma^\wedge})\circ \Delta_{\Omega^\bullet(P)}
 \\
 &=& 
(\id_{\Omega^\bullet(P)}\otimes \chi^{-1}\otimes \chi)\circ (\id_{\Omega^\bullet(P)}\otimes \Delta)\circ \Delta_{\Omega^\bullet(P)}
 \\
 &=& 
(\id_{\Omega^\bullet(P)}\otimes ((\chi^{-1}\otimes \chi)\circ \Delta))\circ \Delta_{\Omega^\bullet(P)}
\\
 &=&(\id_{\Omega^\bullet(P)}\otimes \epsilon)\circ \Delta_{\Omega^\bullet(P)} 
 \\
 &=&
\id_{\Omega^\bullet(P)},
\end{eqnarray*}
and a similar calculation proves that $$\F_{\chi}\circ \F_{\chi^{-1}}=\id_{\Omega^\bullet(P)}.$$  Due to the fact that $$\mu \,\in \, \Omega^\bullet(B) \;\Longleftrightarrow \; \Delta_{\Omega^\bullet(P)}(\mu)=\mu\otimes \mathbbm{1} $$ it is clear that $$\F_{\chi}|_{\Omega^\bullet(B)}=\id_{\Omega^\bullet(B)}.$$  Finally, a direct calculation as before, proves that $\F_{\chi}$ is a qgt if and only if 
\begin{equation}
\label{4.f9}
(\id_{\Gamma^\wedge}\otimes \chi)\circ \Delta=(\chi\otimes \id_{\Gamma^\wedge})\circ \Delta.
\end{equation}

\noindent In this way, if one considers the submonoid  ${H}'_{cl}$  of $H_{cl}$ such that equation (\ref{4.f9}) holds, then it is possible to define the monoid morphism
\begin{equation}
\label{4.f10}
Y: {H}'_{cl} \longrightarrow \qGG,\qquad \chi \longmapsto \F_{\chi}.
\end{equation}

If $\Gamma^\wedge$ is cocommutative, then ${H}'_{cl}=H_{cl}$ and $Y$ is a group morphism. This is the {\it quantum} counterpart of the {\it classical} fact that, for a given principal $G$--bundle $\pi:P\longrightarrow B$ with $G$ abelian, the diffeomorphism 
\begin{equation*}
r_A:P \longrightarrow P,\qquad x  \longmapsto xA
\end{equation*}
is a gauge transformation for all $A$ $\in$ $G$.\\

As we have mentioned at the beginning of this section, we are following the work developed in \cite{br} but in the context of Durdevich's formulation. In the literature, there are other papers that deal with the $0$--degree quantum gauge group. A particularly important analysis can be found in \cite{landi}. 

Since the quantum gauge group is quite large, it is really difficult to work with (even in the $0$--degree case). For that reason, it is natural to work only with special subgroups for certain cases, as we did in equation (\ref{4.f10}) when $\Gamma^\wedge$ is cocommutative, and as the authors of \cite{landi} also did. In \cite{landi}, the authors determined when it is {\it natural} to consider algebra morphisms as quantum gauge transformations. Specifically, when the $\ast$--Hopf algebra $H^\infty$ is coquasitriangular and the total quantum space $P$ is quasi--commutative (and hence the quantum base space $B$ lies in the center of $P$), there is a natural definition of the quantum gauge group of the qpb by $H$--equivariant algebra maps $\underline{H}\longrightarrow P$ \cite{landi}. 

Unfortunately, in general, it is not possible to extend the theory developed in \cite{landi} to the level of differential algebras. For instance, in Example \ref{e.1} we showed explicitly that there is no way to extend the coquasitriangular structure to the universal graded differential calculus for the quantum group $\G$ associated with $\Z_2$. Continuing with this case, it follows that the theory of \cite{landi} cannot be applied for degrees greater than or equal to $1$. However, it is worth noting that the authors of \cite{landi} did not claim to propose an approach to quantum gauge transformations that applies to any degree.

On the other hand, Definition \ref{qgg} allows us to work with these cases for every degree. This should not be a surprise, since we have defined $\qGG$ in the most general way.

The following example illustrates the theory developed in this subsection. In particular, we compute the quantum translation map in all degrees and exhibit several quantum gauge transformations  together with their corresponding convolution invertible maps.

\begin{Example}
\label{e.1.1}
Let $(B,\cdot,\mathbbm{1}_B,\ast)$ be a quantum space. Considering $\G$ as in Example \ref{e.1}, we define the quantum principal $\G$--bundle $$\zeta=(P:=B\otimes H, B, \Delta_P:=\id_B\otimes \Delta).$$ 

The next step is to endow $\zeta$ with a differential calculus. Let us consider $(\Omega^\bullet(B),d,\ast)$ any graded differential $\ast$--algebra generated by $\Omega^0(B)=B$. Furthermore, let us take the universal differential envelope $\ast$--calculus $(\Gamma^\wedge,d,\ast)$ of Example \ref{e.1}. By defining (using the corresponding tensor products) 
\begin{equation}
    \label{caltriv0}
\Omega^\bullet(P):=\Omega^\bullet(B)\otimes \Gamma^\wedge, \quad \Delta_{\Omega^\bullet(P)}:=\id_{\Omega^\bullet(B)}\otimes \Delta
\end{equation}
we obtain a differential calculus for $\zeta$.

Let us calculate the $0$--degree quantum translation map. Notice that $$\T=\{\delta^\C_\triv,\; \delta^\C_\alt\}$$ is a complete set of mutually non--equivalent irreducible $\G$--corepresentations, where $$\delta^\C_\triv:\C\longrightarrow \C\otimes H,\qquad z\longmapsto z\otimes \beta^0_1$$ and $$\delta^\C_\alt:\C\longrightarrow \C\otimes H,\qquad z\longmapsto z\otimes \beta^0_2,$$ with  (the unital element of $H$)  $$\beta^0_1:=\mathbbm{1}_H=\phi_0+\phi_1,$$  and  $$\beta^0_2:=\phi_0-\phi_1.$$ In this way, $$\beta^0:=\{\beta^0_1,\beta^0_2 \}$$ is the linear basis of $H$ given by Theorem \ref{rep}. Since $\{1\}$ is an orthonormal basis of $\C$ with respect to the inner product that makes $\delta^\C_\triv$ and $\delta^\C_\alt$ unitary, the set of left $B$--generators of Remark (\ref{rema}) are given by $$T^\l_\triv: \C\longrightarrow B\otimes H, \qquad z\longmapsto z\mathbbm{1}_B\otimes \beta^0_1$$ and  $$T^\l_\alt: \C\longrightarrow B\otimes H, \qquad z\longmapsto z\mathbbm{1}_B\otimes \beta^0_2.$$ By equation (\ref{qtrs0}), we have  
\begin{equation}
    \label{eq.qtrsb1}
    \qtrs(\beta^0_1)=T^\l_\triv(1)^\ast\otimes_B T^\l_\triv(1)=(\mathbbm{1}_B\otimes \beta^0_1)\otimes_B (\mathbbm{1}_B\otimes \beta^0_1)=\mathbbm{1}_P\otimes_B \mathbbm{1}_P=:[\beta^0_1]_1\otimes_{B}[\beta^0_1]_2
\end{equation}
and 
\begin{equation}
    \label{eq.qtrsb2}
    \qtrs(\beta^0_2)=T^\l_\alt(1)^\ast\otimes_B T^\l_\alt(1)=(\mathbbm{1}_B\otimes \beta^0_2)\otimes_B (\mathbbm{1}_B\otimes \beta^0_2)=:[\beta^0_2]_1\otimes_{B}[\beta^0_2]_2.
\end{equation}

Now, let us calculate the quantum translation map for higher degrees. Since $(\Gamma^\wedge,d,\ast)$ is the universal graded differential calculus (see Example \ref{e.1}) and remembering that  $$\mathfrak{qg}^\#={\Ker(\epsilon)\over \mathcal{R}}=\Ker(\epsilon)=\mathrm{span}_\C\{\pi(\phi_1)\},$$  there is a canonical linear basis of $\Gamma^\wedge$ given by 
$$\beta=\beta^0\oplus \beta^1\oplus\beta^2\oplus\cdots \oplus \beta^n \oplus \cdots,$$ where  $$\beta^j=\{\beta^j_1:=\beta^0_1 \overbrace{\pi(\phi_1)\cdots  \pi(\phi_1)}^{j-\mathrm{times}}=\overbrace{\pi(\phi_1)\cdots  \pi(\phi_1)}^{j-\mathrm{times}},\;\;\beta^j_2:=\beta^0_2 \overbrace{\pi(\phi_1)\cdots\pi(\phi_1)}^{j-\mathrm{times}}=\beta^0_2\,\beta^j_1\}$$ is a linear basis of $\Gamma^{\wedge\,j}$ for $j$ $\in$ $\N$.

Consider the real qpc
\begin{equation}
\label{qpctrz}
\omega^\triv:\mathfrak{qg}^\# \longrightarrow \Omega^1(P),\qquad \theta\longmapsto \mathbbm{1}_B\otimes \theta.
\end{equation}
By equation (\ref{adjoinphi}), we have 
$$\ad(\beta^1_1)=\ad(\pi(\phi_1))=\pi(\phi_1)\otimes \mathbbm{1}_H=\beta^1_1\otimes \beta^0_1=:\beta^{1\,(0)}_1\otimes \beta^{1\,(1)}_1$$ and by equation (\ref{6.f1.6.1}), we obtain
\begin{eqnarray}
\label{eq.qtrsb2.1}
\quad \qtrs(\beta^1_1) =\mathbbm{1}_P\otimes_{\Omega^\bullet(B)}(\mathbbm{1}_B\otimes \beta^1_1)-(\mathbbm{1}_B\otimes \beta^1_1)\otimes_{\Omega^\bullet(B)} \mathbbm{1}_P 
=
[\beta^1_1]_1\otimes_{\Omega^\bullet(B)} [\beta^1_1]_2. 
\end{eqnarray}
Furthermore, according to equation (\ref{6.f1.7}), we get
\begin{eqnarray*}
\qtrs(\beta^1_2)=\qtrs(\beta^0_2\beta^1_1)&=&[\beta^1_1]_1\, \qtrs(\beta^0_2)\,[\beta^1_1]_2  
 \\
 &=&
 (\mathbbm{1}_B\otimes \beta^0_2)\otimes_{\Omega^\bullet(B)} (\mathbbm{1}_B\otimes \beta^1_2)+(\mathbbm{1}_B\otimes \beta^1_2)\otimes_{\Omega^\bullet(B)} ( \mathbbm{1}_B\otimes \beta^0_2)
 \\
 &=& [\beta^1_2]_1\otimes_{\Omega^\bullet(B)}[\beta^1_2]_2,
\end{eqnarray*}
Finally, by equation (\ref{6.f1.10}) we get for $j$ $\geq$ $2$
\begin{equation}
    \label{eq.qtrsb3}
    \qtrs(\beta^j_1)=\qtrs(\beta^{j-1}_1\beta^1_1)=(-1)^{(j-1)\,\partial [\beta^1_1]_1} \,[\beta^1_1]\, \qtrs(\beta^{j-1}_1) \, [\beta^1_1]_2=:[\beta^j_1]_1 \otimes_{\Omega^\bullet(B)} [\beta^j_1]_2,
\end{equation}
and 
\begin{equation}
    \label{eq.qtrsb4}
    \qtrs(\beta^j_2)=\qtrs(\beta^0_2\beta^{j}_1)=[\beta^j_1]_1 \,\qtrs(\beta^0_2)\,[\beta^{j}_1]_2=:[\beta^j_2]_1 \otimes_{\Omega^\bullet(B)} [\beta^j_2]_2.
\end{equation}
These recursive formulas completely characterize $\qtrs$ in all degrees.

The purpose of this example now shifts to showing some quantum gauge transformations. However, first, we need to calculate  the coproduct $\Delta$ for all degrees in terms of the linear basis $\beta$. A straightforward calculation proves that 
\begin{equation}
    \label{coprodbeta1}
    \Delta(\beta^0_1)=\beta^0_1\otimes\beta^0_1, \qquad \Delta(\beta^0_2)=\beta^0_2\otimes\beta^0_2
\end{equation}
and we know that (see equation (\ref{deltapi}))
\begin{equation}
    \label{porqueyo5}
    \Delta(\beta^1_1)=\beta^1_1\otimes \beta^0_1+\beta^0_1\otimes \beta^1_1.
\end{equation}
Since $\Delta$ is multiplicative, we obtain
\begin{equation}
    \label{porqueyo6}
\Delta(\beta^1_2)=\Delta(\beta^0_2)\Delta(\beta^1_1)=\beta^0_2\,\beta^1_1\otimes \beta^0_2\,\beta^0_1+\beta^0_2\,\beta^0_1\otimes \beta^0_2\,\beta^1_1 =\beta^1_2\otimes \beta^0_2+\beta^0_2\otimes \beta^1_2.
\end{equation}
 In general, the multiplicativity property of $\Delta$ implies that
\begin{equation}
    \label{finaldeltapi1}
    \Delta(\beta^j_1)=\overbrace{\Delta(\beta^1_1)\cdots \Delta(\beta^1_1)}^{j-\mathrm{times}}=\Delta(\beta^1_1)\Delta(\beta^{j-1}_1), \qquad \Delta(\beta^j_2)=\Delta(\beta^0_2)\Delta(\beta^j_1)
\end{equation}
for $j\geq 2$. Since the elements $\beta^1_1\otimes \beta^0_1$, $\beta^0_1\otimes\beta^1_1$ are elements of a tensor product of the graded differential $\ast$--algebras, it is easy to check that they  anti--commute each other. Hence, by equation (\ref{finaldeltapi1}), it follows that
\begin{equation}
    \label{finaldeltapi2}
    \Delta(\beta^j_1)=\sum^j_{k=0}\begin{pmatrix}
j \\
k \\
\end{pmatrix}_{-1}\, \beta^{j-k}_1\otimes \beta^k_1
\end{equation}
and
\begin{equation}
    \label{finaldeltapi3}
\Delta(\beta^j_2)=\Delta(\beta^0_2)\Delta(\beta^j_1)=\sum^j_{k=0}\begin{pmatrix}
j \\
k \\
\end{pmatrix}_{-1}\, \beta^0_2\,\beta^{j-k}_1\otimes \beta^0_2\,\beta^k_1=\sum^j_{k=0}\begin{pmatrix}
j \\
k \\
\end{pmatrix}_{-1}\, \beta^{j-k}_2\otimes \beta^k_2
\end{equation}
for $j\geq 2$, where  $\begin{pmatrix}
j \\
k \\
\end{pmatrix}_{-1}$  is the $q$--binomial coefficient for $q=-1$ (\cite{libro,ks}). Equations (\ref{coprodbeta1}), (\ref{porqueyo5}), (\ref{porqueyo6}), (\ref{finaldeltapi2}), (\ref{finaldeltapi3}) completely characterize $\Delta$ in all degrees.

Let $\F$ be a graded linear map of the form
\begin{equation}
    \label{gaugetran1}
    \F:=\id_{\Omega^\bullet(B)}\otimes C: \Omega^\bullet(P)\longrightarrow \Omega^\bullet(P),
\end{equation}
where $$C:\Gamma^\wedge\longrightarrow \Gamma^\wedge$$ is a graded linear isomorphism such that $$C(\beta^0_1)=\beta^0_1\qquad \mbox{  and } \qquad C\, \in \,\Mor(\Delta,\Delta).$$ Since $\Delta_{\Omega(P)}:=\id_{\Omega^\bullet(B)}\otimes \Delta$, it directly follows that $\F$ is a qgt.  Using mathematical induction and our previous characterization of the coproduct $\Delta$ in all degrees, one can easily proven that such isomorphisms $C$ satisfy $$C(\beta^j_1)=\beta^j_1,\qquad C(\beta^j_2)=z\,\beta^j_2 $$ with $z$ $\in$ $\C-\{ 0\}$, for all $j$ $\in$ $\N_0=\N \cup \{0\}$. Of course, we have
\begin{equation}
    \label{gaugetran11}
    \F^{-1}:=\id_{\Omega^\bullet(B)}\otimes C^{-1}: \Omega^\bullet(P)\longrightarrow \Omega^\bullet(P) 
\end{equation}
where $$C^{-1}(\beta^j_1)=\beta^j_1, \;\;C^{-1}(\beta^j_2)=z^{-1}\,\beta^j_2.$$

It is worth mentioning that not every qgt has the form shown in equation (\ref{gaugetran1}). In fact, let $$C:\Gamma^\wedge \longrightarrow \Gamma^\wedge$$ be a graded linear isomorphism  as above, and let $$A:\Gamma^\wedge\longrightarrow \Omega^\bullet(B)$$ be a graded linear map such that $$A(\beta^1_1)\not=0\qquad \mbox{ and }\qquad  A(\beta^j_l)=0,$$ with $\beta^j_l$ any other element of the basis $\beta$. Then, the graded left $\Omega^\bullet(B)$--module morphism 
\begin{equation}
    \label{z2guage}
    \F:\Omega^\bullet(P)\longrightarrow \Omega^\bullet(P)
\end{equation}
such that $$\F(\mu\otimes \vartheta) = \mu\otimes C(\vartheta) +\mu\, A(\vartheta)\otimes \beta^0_1$$
 where $\mu$ $\in$ $\Omega^\bullet(B)$, $\vartheta$ $\in$ $\Gamma^\wedge$, is a qgt with inverse given by 
 \begin{equation}
     \label{z2guage.1}
     \F^{-1}(\mu\otimes \vartheta) = \mu\otimes C^{-1}(\vartheta) -\mu A(\vartheta)\otimes \beta^0_1.
 \end{equation}

To finalize our analysis of quantum gauge transformations, we will characterize all  $0$--degree qgt. We claim that every $0$--degree qgt 
\begin{equation}
        \label{0degreeqgt1}
        \F: P\longrightarrow P
    \end{equation}
is given by
    \begin{equation*}
        \F(b\otimes \beta^0_1)=b\otimes \beta^0_1,\qquad \F(b\otimes \beta^0_2)=b\,\widetilde{b}\otimes \beta^0_2,
    \end{equation*}
    for every $b$ $\in$ $B$, where $\widetilde{b}$ is an invertible element of $B$. Indeed, let $\F$ be a $0$--degree qgt. Since $\F(\mathbbm{1}_P)=\mathbbm{1}_P$ and $\F$ is a left $B$--module morphism, it follows that $$\F(b\otimes \beta^0_1)=b\,\F(\mathbbm{1}\otimes \beta^0_1)=b\,\F(\mathbbm{1}_P)=b\,\mathbbm{1}_P=b\otimes \beta^0_1.$$ On the other hand, if $$\F(b\otimes \beta^0_2)=b\,\widetilde{a}\otimes \beta^0_1+b\,\widetilde{b}\otimes \beta^0_2$$ with $\widetilde{a}$, $\widetilde{b}$ $\in$ $B$, by equation (\ref{coprodbeta1}) we have $$\Delta_P(\F(b\otimes \beta^0_2))= b\,\widetilde{a}\otimes \Delta(\beta^0_1)+b\,\widetilde{b}\otimes \Delta(\beta^0_1)=b\,\widetilde{a}\otimes \beta^0_1\otimes \beta^0_1+b\,\widetilde{b}\otimes \beta^0_2\otimes \beta^0_2$$ and $$(\F\otimes \id_H)\Delta_P(b\otimes \beta^0_2)= \F(b\otimes \beta^0_2)\otimes \beta^0_2=b\,\widetilde{a}\otimes \beta^0_1\otimes \beta^0_2+b\,\widetilde{b}\otimes \beta^0_2\otimes \beta^0_2.$$
    Since $$\Delta_P(\F(b\otimes \beta^0_2))=(\F\otimes \id_H)\Delta_P(b\otimes \beta^0_2)$$ and the elements $\beta^0_1\otimes \beta^0_1$, $\beta^0_2\otimes \beta^0_2,$ $\beta^0_1\otimes \beta^0_2$ are linear independent, we get $\widetilde{a}=0$. Following the same strategy, we can prove that 
    \begin{equation}
        \label{0degreeqgtinverse}
        \F^{-1}(b\otimes \beta^0_1)=b\otimes \beta^0_1,\qquad \F^{-1}(b\otimes \beta^0_2)=b\,\widetilde{c}\otimes \beta^0_2, 
    \end{equation}
    for some $\widetilde{c}$ $\in$ $B$. Finally, the identity  $\F^{-1}\circ \F=\F\circ \F^{-1}=\id_P$  directly implies that $\widetilde{b}$ is invertible and $\widetilde{c}=\widetilde{b}^{-1}$. 
This proves our claim.

As we have previously mentioned, the quantum gauge group is, in general, quite large. For the purposes of this work, it is not necessary to provide a complete characterization of all qgt's. It is sufficient to focus on those given in equations (\ref{gaugetran1}), (\ref{z2guage}) and (\ref{0degreeqgt1}), since these qgt's will be enough to illustrate both the theory developed so far and the results of the following sections.

For example, it follows from equations (\ref{coprodbeta1})--(\ref{porqueyo6}), (\ref{finaldeltapi2}), (\ref{finaldeltapi3})  that $\Gamma^\wedge$ is cocommutative. Therefore (see equations (\ref{4.f7}), (\ref{4.f10})) $${H}'_{cl}=H_{cl}=\Z_2.$$ It is straightforward to verify that the elements of $Y(H_{cl})$ are qgt's as in equation (\ref{gaugetran1}) for graded linear isomorphisms $C$ with $z=1$ or $z=-1$.

To finalize this example, let us show the convolution invertible maps of some of the qgt's previously presented. In accordance with Proposition \ref{4.1}, all invertible convolution maps are given by
\begin{equation*}
        \f_\F:=m_{\Omega^\bullet}\circ (\id_{\Omega^\bullet(P)}\otimes_{\Omega^\bullet(B)}\F)\circ \qtrs: \Gamma^\wedge\longrightarrow \Omega^\bullet(P)
    \end{equation*}
for $\F$ $\in$ $\qGG$. If $\F$ has the form showed in equation (\ref{z2guage}), by equation (\ref{eq.qtrsb1}) we get
\begin{equation}
\label{convo0}
    \f_\F(\beta^0_1)=m_{\Omega^\bullet} (\id_{\Omega^\bullet(P)}\otimes_{\Omega^\bullet(B)}\F)\qtrs(\beta^0_1)=m_{\Omega^\bullet} (\mathbbm{1}_P\otimes_{\Omega^\bullet(B)}\F(\mathbbm{1}_P))=\mathbbm{1}_P
\end{equation}
and by equation (\ref{eq.qtrsb2}) we obtain
\begin{eqnarray}
\label{convo1}
    \f_\F(\beta^0_2)&=&m_{\Omega^\bullet} (\id_{\Omega^\bullet(P)}\otimes_{\Omega^\bullet(B)}\F)\qtrs(\beta^0_2)
    \\
    &=&
    m_{\Omega^\bullet} ((\mathbbm{1}_B\otimes \beta^0_2)\otimes_{\Omega^\bullet(B)}\F(\mathbbm{1}_B\otimes \beta^0_2))=z\,\mathbbm{1}_P.\nonumber
\end{eqnarray}
Notice that (see equation (\ref{z2guage.1}))
\begin{equation}
\label{convo0.1}
    \f^{-1}_\F(\beta^0_1)=m_{\Omega^\bullet} (\id_{\Omega^\bullet(P)}\otimes_{\Omega^\bullet(B)}\F^{-1})\qtrs(\beta^0_1)=m_{\Omega^\bullet} (\mathbbm{1}_P\otimes_{\Omega^\bullet(B)}\F^{-1}(\mathbbm{1}_P))=\mathbbm{1}_P
\end{equation}
and 
\begin{eqnarray}
\label{convo1.1}
    \f^{-1}_\F(\beta^0_2)&=&m_{\Omega^\bullet} (\id_{\Omega^\bullet(P)}\otimes_{\Omega^\bullet(B)}\F^{-1})\qtrs(\beta^0_2)
    \\
    &=&
    m_{\Omega^\bullet} ((\mathbbm{1}_B\otimes \beta^0_2)\otimes_{\Omega^\bullet(B)}\F^{-1}(\mathbbm{1}_B\otimes \beta^0_2))=z^{-1}\,\mathbbm{1}_P.\nonumber
\end{eqnarray}
In light of equation (\ref{eq.qtrsb2.1}) and since all elements of $\mathfrak{qg}^\#$ are of the form $\theta=w\,\pi(\phi_1)=w\,\beta^1_1$ with $w$ $\in$ $\C$, we have that
\begin{eqnarray}
    \label{convo2}
    \f_\F(\theta)=w\,m_{\Omega^\bullet} (\id_{\Omega^\bullet(P)}\otimes_{\Omega^\bullet(B)}\F)\qtrs(\beta^j_1)&=&w\,(\F(\mathbbm{1}_B\otimes \beta^1_1) -\mathbbm{1}_B\otimes \beta^1_1)\nonumber
    \\
    &=&
    w\,(\mathbbm{1}_B\otimes C(\beta^1_1)+A(\beta^1_1)\otimes \beta^0_1-\mathbbm{1}_B\otimes \beta^1_1)
    \\
    &=&
   w\,(\mathbbm{1}_B\otimes \beta^1_1+A(\beta^1_1)\otimes \beta^0_1-\mathbbm{1}_B\otimes \beta^1_1)\nonumber
   \\
    &=&
   w\,A(\beta^1_1)\otimes \beta^0_1=A(\theta)\otimes \beta^0_1.\nonumber 
\end{eqnarray} 
Finally, by equations (\ref{eq.qtrsb1}), (\ref{eq.qtrsb2}), (\ref{0degreeqgt1}), every 0--degree invertible convolution map is given by
\begin{equation}
    \label{0convolution1}
    \f_\F(\beta^0_1)=m_{\Omega^\bullet} (\id_{\Omega^\bullet(P)}\otimes_{\Omega^\bullet(B)}\F)\qtrs(\beta^0_1)=m_{\Omega^\bullet} (\mathbbm{1}_P\otimes_{\Omega^\bullet(B)}\F(\mathbbm{1}_P))=\mathbbm{1}_P
\end{equation}
\begin{eqnarray}
\label{0convolution2}
    \f_\F(\beta^0_2)&=&m_{\Omega^\bullet} (\id_{\Omega^\bullet(P)}\otimes_{\Omega^\bullet(B)}\F)\qtrs(\beta^0_2)
    \\
    &=&
    m_{\Omega^\bullet} ((\mathbbm{1}_B\otimes \beta^0_2)\otimes_{\Omega^\bullet(B)}\F(\mathbbm{1}_B\otimes \beta^0_2))=\widetilde{b}\nonumber
\end{eqnarray}
and by equations (\ref{eq.qtrsb1}), (\ref{eq.qtrsb2}), (\ref{0degreeqgtinverse}), we have 
\begin{equation}
    \label{0convolution1.1}
    \f^{-1}_\F(\beta^0_1)=m_{\Omega^\bullet} (\id_{\Omega^\bullet(P)}\otimes_{\Omega^\bullet(B)}\F^{-1})\qtrs(\beta^0_1)=m_{\Omega^\bullet} (\mathbbm{1}_P\otimes_{\Omega^\bullet(B)}\F^{-1}(\mathbbm{1}_P))=\mathbbm{1}_P
\end{equation}
\begin{eqnarray}
\label{0convolution2.1}
    \f^{-1}_\F(\beta^0_2)&=&m_{\Omega^\bullet} (\id_{\Omega^\bullet(P)}\otimes_{\Omega^\bullet(B)}\F^{-1})\qtrs(\beta^0_2)
    \\
    &=&
    m_{\Omega^\bullet} ((\mathbbm{1}_B\otimes \beta^0_2)\otimes_{\Omega^\bullet(B)}\F^{-1}(\mathbbm{1}_B\otimes \beta^0_2))=\widetilde{b}^{-1}.\nonumber
\end{eqnarray}
\end{Example}

It follows from Examples \ref{e.1}, \ref{e.1.1} that

\begin{Corollary}
\label{coro}
    $\Gamma^{\wedge\,\infty}$ is a graded differential cocommutative $\ast$--Hopf algebra that is not coquasitriangular.
\end{Corollary}

\subsection{Action on Quantum Connections}

In differential geometry, one of the central ideas of gauge theory is to study classes of objects transformable from one to another by gauge transformations.  Probably, one of the most important examples of this arises when the gauge group acts on the set of principal connections, since this reverberates in an action on associated linear connections. The purpose of this subsection is to develop the {\it non--commutative geometrical} counterpart of these actions.

The proof of the following theorem is straightforward. However, since this theorem introduces an action of $\qGG$ on the space of quantum principal connections that differs from the one in \cite{br,ha}, we will present its proof. 

\begin{Theorem}
\label{4.3}
Let $\zeta=(P, B,\Delta_P)$ be a qpb with a differential calculus, and $\F$ $\in$ $\qGG$. Then, for every qpc $\omega$
\begin{enumerate}
\item The linear map $\F^\circledast \omega$ defined by
\begin{equation}
\label{4.f11}
\F^\circledast \omega:=\F\circ \omega
\end{equation}
is again a qpc, and this defines a group action of $\qGG$ on the set of all qpc's $\mathfrak{qpc}(\zeta)$.
\item If $\F$ preserves the $\ast$ operation, $\F^\circledast \omega$ is real if and only if $\omega$ is real.
\item If $\F$ is a graded algebra morphism, $\F^\circledast \omega$ is regular when $\omega$ is regular.
\item If $\F$ is a graded algebra morphism, $\F^\circledast \omega$ is multiplicative when $\omega$ is multiplicative.
\end{enumerate}
\end{Theorem}

\begin{proof}
     
\begin{enumerate}
    \item It is clear that $$\F^\circledast \omega:\mathfrak{qg}^\#\longrightarrow \Omega^1(P)$$ is a linear map. Furthermore, by the equation (\ref{4.f4}), for every $\theta$ $\in$ $\mathfrak{qg}^\#$ we have
    \begin{eqnarray*}
\Delta_{\Omega^\bullet(P)}((\F^\circledast \omega)(\theta))  \;=\; (\F\otimes \id_{\Gamma^\wedge})\Delta_{\Omega^\bullet(P)}(\omega(\theta)) &=& (\F\otimes \id_{\Gamma^\wedge})((\omega\otimes \id_H)\ad(\theta)+\mathbbm{1}\otimes \theta)
 \\
 &=&
((\F^\circledast \omega)\otimes \id_H)\ad(\theta)+\mathbbm{1}\otimes\theta
\end{eqnarray*}
and hence, $\F^\circledast \omega$ is a qpc. This defines a right group action of $\qGG$ on $\mathfrak{qpc}(\zeta)$ (the product on $\qGG$ is given by $(\F_1\cdot \F_2)(w)=\F_2(\F_1(w))$).
\item Assume that $\F$ preserves the $\ast$ operation. Then  $$(\F^\circledast \omega)(\theta)^\ast=\F(\omega(\theta))^\ast=\F(\omega(\theta)^\ast).$$ Now, it directly follows that $(\F^\circledast \omega)(\theta)^\ast=(\F^\circledast \omega)(\theta^\ast)$ if and only $\omega(\theta)^\ast=\omega(\theta^\ast)$.
\item  Assume that $\F$ is a graded algebra morphism. Let $ \omega$ be a regular qpc (see equation (\ref{2.f25})) and let $\varphi$ $\in$ $\Hor^k P$. Thus, by Remark \ref{remahor} we have
$$\F^{-1}(\varphi)\,\in\,\Hor^k\,P\qquad \mbox{ and }\qquad \Delta_\Hor(\F^{-1}(\varphi))=\F^{-1}(\varphi^{(0)})\otimes \varphi^{(1)}.$$ Hence
\begin{eqnarray*}
((\F^\circledast \omega)(\theta))\varphi=\F(\omega(\theta))\varphi  \;=\; \F(\omega(\theta)\F^{-1}(\varphi))  &=& (-1)^k\,\F(\F^{-1}(\varphi^{(0)})\,\omega(\theta\diamondsuit \varphi^{(1)}))
 \\
 &=&
(-1)^k\,\varphi^{(0)}\, \F(\omega(\theta\diamondsuit \varphi^{(1)}))
 \\
 &=&
(-1)^k\,\varphi^{(0)}\, (\F^\circledast\omega)(\theta\diamondsuit \varphi^{(1)}).
\end{eqnarray*}
Therefore, $\F^\circledast \omega$ is regular.
\item  Assume that $\F$ is a graded algebra morphism and let $\omega$ be a multiplicative qpc (see equation (\ref{2.f26})). So, for all $g$ $\in$ $\mathcal{R}$ we obtain $$0  = \omega(\pi(g^{(1)}))\omega(\pi(g^{(2)})),$$ which implies that
\begin{eqnarray*}
0  = \F(\omega(\pi(g^{(1)}))\omega(\pi(g^{(2)}))) &=&\F(\omega(\pi(g^{(1)})))\F(\omega(\pi(g^{(2)})))
 \\
 &=&
(\F^\circledast\omega)(\pi(g^{(1)}))(\F^\circledast\omega)(\pi(g^{(2)})).
\end{eqnarray*}
Hence, $\F^\circledast\omega$ is multiplicative.
\end{enumerate}  
\end{proof}

In the {\it classical} case, if $\pi:P\longrightarrow B$ is a principal $G$--bundle, the gauge group acts on the space of principal connections by 
\begin{equation}
    \label{gaugeconn}
    F^\#\omega:=\omega \circ dF,
\end{equation}
where $\omega$ is a principal connection, $F: P\longrightarrow P$  is a gauge transformation viewed as a vertical principal bundle automorphism, and $dF: TP\longrightarrow TP$ is its differential \cite{nodg}.  The {\it dualization} of the equation (\ref{gaugeconn}) corresponds to equation (\ref{4.f11}), taking into account that in the {\it non--commutative geometrical} case, the map $\F|_{\Omega^1(P)}:\Omega^1(P)\longrightarrow \Omega^1(P)$ plays the role of $dF$. Moreover, in differential geometry one can prove that
\begin{equation}
    \label{gaugeconn1}
    F^\#\omega=\ad_{f_F}\circ \omega+f_F^\#\theta,
\end{equation}
where $\ad$ is the right adjoint action of $G$ on its Lie algebra $\mathfrak{g}$, $f_F:P\longrightarrow G$ is the $G$--equivariant map associated with $F$, and in this case, $\theta$ is the Maurer--Cartan form \cite{gtvp}. In addition, the curvature of $\omega$ satisfies 
\begin{equation}
    \label{gaugeconn2}
    F^\# R^\omega=\ad_{f_F}\circ R^\omega=R^{F^\#\omega}.
\end{equation}
The following proposition presents the {\it non--commutative geometrical} counterparts of equations (\ref{gaugeconn1}), (\ref{gaugeconn2}).

\begin{Proposition}
\label{4.4}
Given a qpc $\omega$, we get
\begin{equation}
\label{f73}
\F^{\circledast}\omega=m_{\Omega}\circ(\omega\otimes \f_{\F})\circ \ad +\f_{\F},
\end{equation}
where $\f_\F:\Gamma^\wedge\longrightarrow \Omega^\bullet(P)$ is the convolution invertible map associated with $\F$ (see equation (\ref{4.f5})). Furthermore, the curvature of $\omega$ (see equation (\ref{2.f28})) satisfies 
\begin{equation}
\label{f74.1}
\F\circ R^{\omega}=m_{\Omega}\circ(R^{\omega}\otimes \f_\F)\circ \ad.
\end{equation}
In addition, if $\F$ is a graded differential algebra morphism, we have
\begin{equation}
\label{f74.2}
\F\circ R^{\omega}=R^{\F^{\circledast}\omega}.
\end{equation}
\end{Proposition}

\begin{proof}
By Proposition \ref{4.1} we know that $\F=\F_{\f_\F}$. Thus, for every $\theta$ $\in$ $\mathfrak{qg}^\#$ we obtain
\begin{eqnarray*}
(\F^{\circledast}\omega)(\theta)=(\F\circ \omega)(\theta)=(\F_{\f_\F} \circ \omega)(\theta)&=& (m_{\Omega} (\id_{\Omega^\bullet(P)}\otimes \f_\F) \Delta_{\Omega^\bullet(P)} ({\omega}(\theta))
\\
&=& 
(m_{\Omega} (\id_{\Omega^\bullet(P)}\otimes \f_\F)((\omega\otimes \id_H)\ad(\theta)+\mathbbm{1}\otimes\theta)
\\
&=&
 m_{\Omega}(\omega\otimes \f_\F)\ad(\theta)+\f_\F(\theta).
\end{eqnarray*}
This proves the first statement of this proposition.

On the other hand,  since $\Im(R^\omega)\subseteq \Hor^2 P$ and $R^\omega$ $\in$ $\Mor(\ad,\Delta_\Hor)$ (\cite{micho2,stheve}), it follows that
\begin{eqnarray*}
\F\circ R^\omega=\F_{\f_\F} \circ R^\omega&=& m_{\Omega} \circ (\id_{\Omega^\bullet(P)}\otimes \f_\F) \circ \Delta_{\Omega^\bullet(P)} \circ R^{\omega}
\\
&=& 
m_{\Omega} \circ (\id_{\Omega^\bullet(P)}\otimes \f_\F) \circ \Delta_{\Hor} \circ R^{\omega}
\\
&=&
m_{\Omega}\circ (R^{\omega}\otimes \f_\F)\circ \ad,
\end{eqnarray*}
so the second statement of this proposition has been proven.

Finally, let us assume that $\F$ is a graded differential algebra morphism.  If $\Theta(\theta)=\displaystyle \sum_{i,j}\theta_i\otimes \theta'_j$, then
\begin{eqnarray*}
\F(R^{\omega}(\theta)) \;=\; \F(d\omega(\theta))-\F(\langle \omega,\omega\rangle(\theta))) &=& 
d\,\F(\omega(\theta))-\displaystyle \sum_{i,j}\F(\omega(\theta_i))\F(\omega( \theta'_j)) 
\\
 &=&
d\,\F^{\circledast}\omega(\theta)-\langle \F^{\circledast}\omega, \F^{\circledast}\omega\rangle(\theta)
 \\
 &=&
R^{\F^{\circledast}\omega}(\theta)
\end{eqnarray*}
which completes the proof.
\end{proof}

\begin{Example}
    \label{e.1.2.1}
    Continuing with Example \ref{e.1.1}, suppose $\F$ is of the form given in equation (\ref{z2guage}). Consider the qpc $\omega^\triv$ (see equation (\ref{qpctrz})) and let $\theta$ $\in$ $\mathfrak{qg}^\#$. Then, there exists $u$ $\in$ $\C$ such that $\theta=u\,\pi(\phi_1)=u\,\beta^1_1$. This implies that
    \begin{eqnarray}
        \label{qpctrz1}
        \F^\circledast\omega^\triv(\theta)=\F(\mathbbm{1}_B \otimes \theta)&=&\mathbbm{1}_B\otimes C(\theta)+A(\theta)\otimes \beta^0_1 \nonumber
        \\
        &=&
        \mathbbm{1}_B\otimes u\,C(\beta^1_1)+A(\theta)\otimes \beta^0_1
        \\
        &=&
        \mathbbm{1}_B\otimes u\, \beta^1_1+ A(\theta)\otimes \beta^0_1 \nonumber
        \\
        &=&
        \mathbbm{1}_B\otimes \theta+A(\theta)\otimes \beta^0_1. \nonumber
    \end{eqnarray}
    
  On the other hand, by equations (\ref{adjoinphi}), (\ref{convo0}), (\ref{convo2}) we obtain 
\begin{eqnarray*}
 m_{\Omega}(\omega^\triv\otimes \f_\F)\ad(\theta)+\f_\F(\theta) 
 &=&
 m_{\Omega}(\omega^\triv\otimes \f_\F)(\theta\otimes \beta^0_1)+A(\theta)\otimes \beta^0_1
 \\
 &=&
 m_{\Omega}(\omega^\triv(\theta)\otimes \f_\F(\beta^0_1))+A(\theta)\otimes \beta^0_1
 \\
 &=&
 m_{\Omega}((\mathbbm{1}_B \otimes \theta)\otimes \mathbbm{1}_P)+A(\theta)\otimes \beta^0_1
 \\
 &=&
 \mathbbm{1}_B\otimes \theta+A(\theta)\otimes \beta^0_1,
\end{eqnarray*}
which explicitly shows that equation (\ref{f73}) holds.

Now consider the linear map $$\Theta: \mathfrak{qg}^\#\longrightarrow \mathfrak{qg}^\#\otimes \mathfrak{qg}^\#$$ given by $$\Theta(\beta^1_1)=2\,\beta^1_1\otimes \beta^1_1.$$ It is straightforward to prove that $\Theta$ is actually an embedded differential (see Definition \ref{embeddeddifferential}). In this way, for every $\theta=u\,\beta^1_1$ $\in$ $\mathfrak{qg}^\#$ we obtain that
\begin{eqnarray*}
  R^{\omega^\triv}(\theta)=u\,R^{\omega^\triv}(\beta^1_1)&=&u\,(d\omega^\triv(\beta^1_1)-\langle \omega^\triv,\omega^\triv\rangle(\beta^1_1))
  \\
  &=& 
  u\,(\mathbbm{1}_B\otimes d\beta^1_1-2\, \omega^\triv(\beta^1_1)\omega^\triv(\beta^1_1))
  \\
  &=&
  u\,(2\mathbbm{1}_B\otimes\beta^1_1\,\beta^1_1-2\, \mathbbm{1}_B\otimes \beta^1_1\,\beta^1_1)=0.
\end{eqnarray*}
Therefore, $R^{\omega^\triv}=0$ and equation (\ref{f74.1}) holds trivially. 
\end{Example}

In literature, for example \cite{libro,ha},  there is a commonly accepted action of the ($0$--degree) quantum gauge group on qpc's given by
\begin{equation}
    \label{falseact}
    \omega \longmapsto \f \,\widetilde{\ast} \, \omega \,\widetilde{\ast} \, \f^{-1}+ \f \,\widetilde{\ast} \, (d \f^{-1})
\end{equation}
for convolution invertible maps $\f:H\longrightarrow P$.  Here, 
 $\omega$ is considered a map from $H$ to $\Omega^1(P)$ and $\widetilde{\ast}$ is convolution product (see equation (\ref{convolutionproduct})). 

In differential geometry, it is well--known that for a matrix Lie group \cite{nodg} $$F^\#\omega=(f_F)\, \omega\, (f_F)^{-1}+f\,d(f_F)^{-1},$$ so that in the \emph{classical} case, equations (\ref{4.f11}) and (\ref{falseact}) agree (for matrix Lie groups). However, as we have mentioned earlier in this subsection, in the \emph{non--commutative geometrical} case, equations (\ref{4.f11}) and (\ref{falseact}) do not generally agree.

\begin{Proposition}
    \label{noteq}
    In Example \ref{e.1.1}, equations (\ref{4.f11}), (\ref{falseact}) define different actions of $\qGG$ on the space of qpc's $\mathfrak{qpc}(\zeta)$. 
\end{Proposition}
\begin{proof}
    Let us define an action of $\qGG$ on $\mathfrak{qpc}(\zeta)$ by  $$\F^\times \omega:= \f_\F\,\widetilde{\ast} \, \omega \,\widetilde{\ast} \, \f^{-1}_\F+ \f_\F\,\widetilde{\ast} \, (d \f^{-1}_\F),$$ where we have adopted the following abuses of notation: $$\f_\F:=\f_\F|_H:H\longrightarrow P$$ 
    and
    \begin{equation*}
\omega:H \longrightarrow \Omega^1(P),\qquad g\longmapsto \omega (\pi(g)).
\end{equation*}
    By  Proposition \ref{4.1}, every $0$--degree convolution invertible map is of the form $\f_\F:=\f_\F|_H$  for a unique map  $\F:=\F|_P:P\longrightarrow P$. 

   Let $\omega$ be a qpc and let $\F$ be a qgt of the form given in equation (\ref{z2guage}). By equations (\ref{convo0.1})   (\ref{convo1.1}), it immediately follows that $d\f^{-1}_\F=0$. Thus $$\f_\F\,\widetilde{\ast} \, (d \f^{-1}_\F)=0.$$ 

On the other hand, by equation (\ref{coprodbeta1})  we have
$$\beta^{0\,(1)}_1\otimes \beta^{0\,(2)}_1\otimes \beta^{0\,(3)}_1 =\beta^{0}_1\otimes \beta^{0}_1\otimes \beta^{0}_1,$$
$$\beta^{0\,(1)}_2\otimes \beta^{0\,(2)}_2\otimes \beta^{0\,(3)}_2 =\beta^{0}_2\otimes \beta^{0}_2\otimes \beta^{0}_2.$$  Then, by equations (\ref{convo0}), (\ref{convo0.1}) we get
\begin{eqnarray*}
    (\f_\F\,\widetilde{\ast} \, \omega \,\widetilde{\ast} \,\f^{-1}_\F)(\beta^0_1)=\f_\F(\beta^{0\,(1)}_1)\, \omega(\beta^{0\,(2)}_1)\,  \f^{-1}_\F(\beta^{0\,(3)}_1)&=& \f_\F(\beta^{0}_1)\, \omega(\pi(\beta^{0}_1))\,  \f^{-1}_\F(\beta^{0}_1)
    \\
    &=&
    \mathbbm{1}_P\, \omega(\pi(\beta^{0}_1))\,  \mathbbm{1}_P
    \\
    &=&
    \omega(\pi(\beta^{0}_1))=\omega(\pi(\mathbbm{1}))=0
\end{eqnarray*}
and by equations (\ref{convo1}), (\ref{convo1.1}) we obtain
\begin{eqnarray*}
    (\f_\F\,\widetilde{\ast} \, \omega \,\widetilde{\ast} \, \f^{-1}_\F)(\beta^0_2)=\f_\F(\beta^{0\,(1)}_2)\, \omega(\beta^{0\,(2)}_2)\,  \f^{-1}_\F(\beta^{0\,(3)}_2)&=& \f_\F(\beta^{0}_2)\, \omega(\pi(\beta^{0}_2))\,  \f^{-1}_\F(\beta^{0}_2)
    \\
    &=&
   z\,\mathbbm{1}_P\,\omega(\pi(\beta^{0}_2))\,z^{-1}\, \mathbbm{1}_P
    \\
&=&\omega(\pi(\beta^{0}_2)). 
\end{eqnarray*}
Therefore 
\begin{equation}
    \label{wrongaction1}
    \F^\times\omega=\omega.
\end{equation}
Nevertheless,  equation (\ref{qpctrz1}) clearly shows that 
\begin{equation}
    \label{trivact}
    \F^\circledast\omega^\triv\not=\omega^\triv.
\end{equation}
We conclude that both actions are different. 
\end{proof}

\begin{Remark}
    \label{remaaction}
    In the context of the previous proposition, if all invertible elements of $B$ are of the form  $$\{z\,\mathbbm{1}_B\mid z\, \in \,\C-\{0\} \},$$ then, by using the characterization of every $0$--degree invertible convolution map given at the end of Example \ref{e.1.1} (see equations (\ref{0convolution1})--(\ref{0convolution2.1})), we can repeat exactly the same calculations of the proof of Proposition \ref{noteq} and conclude that the action of $\qGG$ on $\mathfrak{qpc}(\zeta)$ under equation (\ref{falseact}) is trivial. However, the action of $\qGG$ on $\mathfrak{qpc}(\zeta)$ under equation (\ref{4.f11}) cannot be trivial (see equation (\ref{trivact})).
\end{Remark} 

It is worth mentioning that the covariant derivatives satisfy (see equations (\ref{2.f30}), (\ref{2.f30}))
\begin{equation}
\label{f75}
D^{\F^{\circledast}\omega}(\varphi)=d\varphi-(-1)^k\varphi^{(0)} \,\F^{\circledast}\omega(\pi(\varphi^{(1)})),
\end{equation}
\begin{equation}
\label{f75.1}
\widehat{D}^{\F^{\circledast}\omega}(\varphi)=d\varphi+ \widehat{\F^{\circledast}\omega}(\pi(\varphi^{(1)}))\varphi^{(0)},
\end{equation}
where $\varphi$ $\in$ $\Hor^k P$.\\

For the rest of this section we shall assume that $\F$ is a graded differential $\ast$--algebra morphism. This happens, for example, for qgt's induced by elements of $\widehat{H_{cl}}$.

\subsection{On Induced Quantum Linear Connections}

In differential geometry, the gauge group acts on associated vector bundles via vector bundle isomorphisms. Specifically, given a principal $\G$--bundle $\pi:P\longrightarrow B$ and an associated vector bundle $\pi_{\alpha}: E^V \longrightarrow B$ for the linear representation $\alpha:H\longrightarrow GL(V)$, a gauge transformation $F$ induces a vector bundle isomorphism defined by \cite{sald1}
\begin{equation*}
\begin{aligned}
\mathbf{A}_F: E^V &\longrightarrow E^V \\
[x,v] &\longmapsto [F(x),v].
\end{aligned}
\end{equation*}
If $\omega$ is a principal connection, then $\mathbf{A}_F$ is a parallel vector bundle isomorphism between ($E^V,\nabla^\omega_\alpha$) and $(E^V,\nabla^{F^\ast\omega}_{\,\alpha})$, and we have $$\nabla^{F^\ast \omega}_{\,\alpha}=(\id_{\Omega^1(B)}\otimes_{C^\infty(B)} \widetilde{\mathbf{A}}_F)\circ \nabla^{\omega}_\alpha \circ \widetilde{\mathbf{A}}^{-1}_F, $$ where  $\widetilde{\mathbf{A}}_F$ denotes the isomorphism of sections of $E^V$ induced by $\mathbf{A}_F$ \cite{sald1}. The curvature satisfies a similar formula \cite{gtvp}. These facts motivate the following theorem.

\begin{Theorem}
\label{4.5}
Let $(\zeta,\omega)$ be a qpb with a qpc, and let $\F$ be a qgt that is also a graded differential $\ast$--algebra morphism. Then, $\F$ defines a left $B$--module automorphism $\mathbf{A}_\F$ of $E^V_\l$ such that  $$(\id_{\Omega^\bullet(B)}\otimes_{B} \mathbf{A}_\F)\,\circ\, \nabla^\omega_V= \nabla^{\F^{\circledast}\omega}_V \, \circ \, \mathbf{A}_\F$$ for a fixed finite--dimensional $\G$--corepresentation $\delta^V$. Furthermore, we have (see equation (\ref{3.f10})) $$(\mathbf{A}_\F \otimes_B \id_{\Omega^\bullet(B)})\circ \sigma_{V}=\sigma_{V}\circ (\id_{\Omega^\bullet(B)}\otimes_B \mathbf{A}_\F).$$  
\end{Theorem}

\begin{proof}
By equation (\ref{dirsumqvb}), it is enough to prove the theorem for $\delta^V$ $\in$ $\T$. Let us start noticing that by equation (\ref{f75}) and Remark \ref{remahor}, we have $$D^{\F^{\circledast}\omega}(\F(\varphi))=d\F(\varphi)-(-1)^k\,\F(\varphi^{(0)}) \,\F(\omega(\pi(\varphi^{(1)})))=\F(d\varphi-(-1)^k\,\varphi^{(0)}\,\omega(\pi(\varphi^{(1)}))) = \F(D^{\omega}(\varphi)),$$ for every $\varphi$ $\in$ $\Hor^k\,P$. Hence $$D^{\F^{\circledast}\omega}\circ \F=\F\circ D^{\omega}.$$ In addition, considering the associated left qpb $E^V_\l$ of $\delta^V$, the linear map 
\begin{equation}
\label{f76}
\mathbf{A}_\F : E^V_\l  \longrightarrow E^V_\l,\qquad
T  \longmapsto \F\circ T
\end{equation}
is well--defined because 
\begin{eqnarray*}
    \Delta_P\circ (\F\circ T)=(\F\otimes \id_H)\circ \Delta_P\circ T &=& (\F\otimes \id_H)\circ (T\otimes \id_H)\circ \delta^V
    \\
    &=&((\F\circ T)\otimes \id_H)\circ \delta^V.
\end{eqnarray*}
This shows that indeed, $\F\,\circ T$ $\in$ $E^V_\l$.  Furthermore, the map $\mathbf{A}_\F$ is a $B$--bimodule isomorphism, and its inverse is  given by 
\begin{equation}
\label{f76.1}
\mathbf{A}^{-1}_\F : E^V_\l  \longrightarrow E^V_\l,\qquad
T  \longmapsto \F^{-1}\circ T.
\end{equation}
Notice that $$\mathbf{A}^{-1}_\F=\mathbf{A}_{\F^{-1}}.$$

In this way, for all $T$ $\in$ $E^V_\l$ we obtain
\begin{eqnarray*}
(\nabla^{\F^{\circledast}\omega}_{\,V} \circ \mathbf{A}_\F)(T) \;=\; \nabla^{\F^{\circledast}\omega}_{\,V}(\F\circ T) &=& \sum^{d_{V}}_{k=1}\mu^{D^{\F^{\circledast}\omega}\circ \F\circ T}\otimes_B T^\l_k
 \\
 &=&
\sum^{d_{V}}_{k=1}\mu^{\F \circ D^{\omega}\circ T}\otimes_B T^\l_k;
\end{eqnarray*}
so $$(\Upsilon^{-1}_{V}\circ \nabla^{\F^{\circledast}\omega}_{\,V} \circ \mathbf{A}_\F)(T)=\F\circ D^{\omega}\circ T.$$ On the other hand,
\begin{eqnarray*}
((\id_{\Omega^\bullet(B)}\otimes_B \mathbf{A}_\F)\circ \nabla^{\omega}_V)(T)  &=& \sum^{d_{V}}_{k=1} (\id_{\Omega^\bullet(B)}\otimes_B \mathbf{A}_\F) (\mu^{D^{\omega}\circ T}\otimes_B T^\l_k)
 \\
 &=& 
\sum^{d_{V}}_{k=1} \mu^{D^{\omega}\circ T}\otimes_B \mathbf{A}_\F(T^\l_k)  \;=\;\sum^{d_{V}}_{k=1}\mu^{D^{\omega}\circ T}\otimes_B \F\circ T^\l_k;
\end{eqnarray*}
thus $$(\Upsilon^{-1}_{V} \circ (\id_{\Omega^\bullet(B)}\otimes_B \mathbf{A}_\F)\circ \nabla^{\omega}_{V})(T)=\sum^{d_{V}}_{k=1}\mu^{D^{\omega}\circ T} \F\circ T^\l_k= \F\circ \sum^{d_{V}}_{k=1}\mu^{D^{\omega}\circ T} T^\l_k= \F\circ D^{\omega}\circ T.$$ By using the fact  that $\Upsilon^{-1}_{V}$ is bijective, we conclude that $\mathbf{A}_\F$ satisfies the first part of the statement.

Let us take $\psi$ $\in$ $\Omega^\bullet(B)\otimes_B E^V_\l$. Then, if $\Upsilon^{-1}_{V}(\psi)=\displaystyle \sum_k T^\r_k\, \widetilde{\mu}_k$, we get $$(\mathbf{A}_\F\otimes_B \id_{\Omega^\bullet(B)})\sigma_{V}(\psi)=\displaystyle\sum_k\mathbf{A}_\F(T^\r_k)\otimes_B \widetilde{\mu}_k =\sum_k \F\circ T^\r_k\otimes_B \widetilde{\mu}_k;$$ thus $$\widehat{\Upsilon}^{-1}_{V}(\mathbf{A}_\F\otimes_B \id_{\Omega^\bullet(B)})\sigma_{V}(\psi)=(\F \circ \Upsilon^{-1}_{V})(\psi).$$ On the other hand, $$\sigma_{V}(\id_{\Omega^\bullet(B)}\otimes_B \mathbf{A}_\F)(\psi)=\displaystyle\sum_k T^\r_k\otimes_B \mu'_k$$ if $(\F\circ \Upsilon^{-1}_{V})(\psi)=\displaystyle \sum_k T^\r_k\,\mu'_k$, because of $(\Upsilon^{-1}_{V}\circ (\id_{\Omega^\bullet(B)}\otimes_B \mathbf{A}_\F))(\psi)=(\F\circ \Upsilon^{-1}_{V})(\psi) $. Hence, $$\widehat{\Upsilon}^{-1}_{V}\sigma_{V}(\id_{\Omega^\bullet(B)}\otimes_B \mathbf{A}_\F)(\psi)= (\F\circ \Upsilon^{-1}_{V})(\psi)$$ and the theorem follows because $\widehat{\Upsilon}^{-1}_{V}$ is bijective.
\end{proof}

\begin{Corollary}
\label{4.6}
The following formula holds: $$\nabla^{\F^{\circledast}\omega}_{V}=(\id_{\Omega^\bullet(B)}\otimes_B \mathbf{A}_\F)\circ \nabla^{\omega}_{V}\circ \mathbf{A}^{-1}_\F=(\id_{\Omega^\bullet(B)}\otimes_B \mathbf{A}_\F)\circ \nabla^{\omega}_{V}\circ \mathbf{A}_{\F^{-1}}.$$
\end{Corollary}

On the other hand, we have

\begin{Proposition}
\label{4.7}
If $\F$ is qgt that is also a graded differential $\ast$--algebra morphism, we obtain $\mathbf{A}_\F$ $\in$ $U(E^V_\l)$ (the space of unitary operators of $E^V_\l$).
\end{Proposition}

\begin{proof}
As before, it is enough to prove the proposition for $\delta^V$ $\in$ $\T$. Then, by taking $T_1$, $T_2$ $\in$ $E^V_\l$ we get
\begin{eqnarray*}
\langle \mathbf{A}_\F(T_1),T_2\rangle_\l\,= \, \sum^{n_{V}}_{k=1} \F(T_1(e_k))T_2(e_k)^\ast &=& \sum^{n_{V}}_{k=1} \F(T_1(e_k)\F^{-1}(T_2(e_k))^\ast)
 \\
 &=& 
\sum^{n_{V}}_{k=1} T_1(e_k)\F^{-1}(T_2(e_k))^\ast  
 \\
 &=&
\langle T_1,\F^{-1}\circ T_2\rangle_\l
 \\
 &=& 
\langle T_1,\mathbf{A}^{-1}_\F(T_2)\rangle_\l, 
\end{eqnarray*}
where in the third equality we have used that $\F(b)=b$ for all $b$ $\in$ $B$ ($\F$ is a graded left $\Omega^\bullet(B)$--module morphism with $\F(\mathbbm{1})=\mathbbm{1}$) and the fact that $$\displaystyle\sum^{n_{V}}_{k=1} T_1(e_k)\F^{-1}(T_2(e_k))^\ast\;\in\; B.$$ We conclude that $\mathbf{A}_\F$ is adjointable with respect to $\langle-,-\rangle_\l$ and $\mathbf{A}^\dagger_\F=\mathbf{A}^{-1}_\F$.
\end{proof}

\noindent Combining Proposition \ref{4.7} and Corollary \ref{4.6} we obtain the following equation 
\begin{equation}
    \label{linargaugetrans}
    \nabla^{\F^{\circledast}\omega}_{V}=(\id_{\Omega^\bullet(B)}\otimes_B \mathbf{A}_\F)\circ \nabla^{\omega}_{V}\circ \mathbf{A}^{\dagger}_\F.
\end{equation}

\begin{Proposition}
    \label{exteriorgauge}
    If $\F$ is qgt that is also a graded differential $\ast$--algebra morphism, the exterior covariant derivative of $\nabla^{\F^{\circledast}\omega}_{V}$ is given by
    \begin{equation}
\label{f77.1}
d^{\nabla^{\F^{\circledast}\omega}_{V}}=(\id_{\Omega\bullet(B)}\otimes_B \mathbf{A}_\F)\circ d^{\nabla^\omega_V}\circ (\id_{\Omega\bullet(B)}\otimes_B \mathbf{A}^\dagger_\F).
\end{equation}
\end{Proposition}
\begin{proof}
    By equation (\ref{dirsumqvb}), it is enough to prove the theorem for $\delta^V$ $\in$ $\T$. In accordance with equations (\ref{3.f10.1}), (\ref{linargaugetrans}), for every $\mu$ $\in$ $\Omega^K(B)$, $T$ $\in$ $E^V_\l$, we obtain
    \begin{eqnarray*}
        d^{\nabla^{\F^{\circledast}\omega}_{V}}(\mu\otimes_B T)&=&d\mu\otimes_B T+(-1)^k \mu\,\nabla^{\F^{\circledast}\omega}_{V}T
        \\
        &=&
        d\mu\otimes_B T+(-1)^k \mu\,(\id_{\Omega^\bullet(B)}\otimes_B \mathbf{A}_\F) \nabla^{\omega}_{V}(\mathbf{A}^{\dagger}_\F(T))
        \\
        &=&
        d\mu\otimes_B T+(-1)^k \sum^{d_V}_{j=1} \mu\,\mu^{D^\omega \circ \mathbf{A}^{\dagger}_\F(T)}_j  \otimes_B \mathbf{A}_\F(T^\l_j)
        \\
        &=&
        (\id_{\Omega^\bullet(B)}\otimes_B \mathbf{A}_\F)(d\mu\otimes_B \mathbf{A}^{\dagger}_\F(T)+(-1)^k \sum^{d_V}_{j=1} \mu\,\mu^{D^\omega \circ \mathbf{A}^{\dagger}_\F(T)}_j  \otimes_B T^\l_j)
        \\
        &=&
        (\id_{\Omega^\bullet(B)}\otimes_B \mathbf{A}_\F)\,d^{\nabla^\omega_V}(\mu\otimes \mathbf{A}^{\dagger}_\F(T))
        \\
        &=&
        ((\id_{\Omega^\bullet(B)}\otimes_B \mathbf{A}_\F)\circ d^{\nabla^\omega_V}\circ (\id_{\Omega^\bullet(B)}\otimes_B \mathbf{A}^\dagger_\F))(\mu\otimes T)
    \end{eqnarray*}
    and hence, equation (\ref{f77.1}) holds.
\end{proof}

Combining equations (\ref{3.f10.2}), (\ref{linargaugetrans}) and the last proposition, it immediately follows that
\begin{Proposition}
    \label{curvaturegauge}
    The curvature of a qpc $\omega$ satisfies the following equation
    \begin{equation}
\label{f77}
R^{\nabla^{\F^{\circledast}\omega}_{V}}=(\id_{\Omega^\bullet(B)}\otimes_B \mathbf{A}_\F)\circ R^{\nabla^{\omega}_{V}}\circ \mathbf{A}^\dagger_\F.
\end{equation}
\end{Proposition}

Of course, there are similar results for $(E^V_\r, \widehat{\nabla}^{\omega}_{V})$, $(E^V_\r, \widehat{\nabla}^{\F^{\circledast} \omega}_{\,V})$, and

\begin{equation}
\label{f78}
\widehat{\mathbf{A}}_\F : E^V_\r  \longrightarrow E^V_\r,\qquad
T  \longmapsto \widehat{\F}\circ T,
\end{equation}
where $\widehat{\F}:=\ast \circ \F \circ \ast$.

\begin{Remark}
\label{4.9}
Notice that to define $\mathbf{A}_\F$ and $\widehat{\mathbf{A}}_\F$, it is not necessary for $\F$ to be a graded differential $\ast$--algebra morphism; our definition works for any qgt $\F$, which in turn induces a natural group action of $\qGG$ on $E^V_\l$ and $E^V_\r$. However, as we have verified in this subsection, the properties of $\mathbf{A}_\F$ (and $\widehat{\mathbf{A}}_\F$) with respect to the induced qlc's require that $\F$ be a graded differential $\ast$--algebra morphism.
\end{Remark}

\section{Examples}

Just as we have mentioned earlier, the primary purpose of this paper is to present the canonical Hermitian structure, study its properties and its relationship with qlc's. Our second goal is to develop an  {\it ad hoc} definition of the quantum gauge group for a given qpb with a differential calculus.  These topics were addressed in the previous sections, and we have provided an elementary example to illustrate our theory concerning the quantum gauge group.  It is worth mentioning that the only assumption made throughout this paper was in Remark \ref{rema}, where we assumed that the quantum base space $B$ can be completed to a $C^\ast$--algebra. This assumption makes it very easy to find qpb's for which our theory applies. In this section, we present two different classes of examples. In \cite{appendix}, the reader can find another class of examples of our theory applied to special kinds 
of {\it classical/quantum hybrid} principal bundles involving Dunkl operators as covariant derivatives of qpc's (\cite{ds}).

\subsection{Trivial Quantum Principal Bundles}

Mirroring the {\it classical} case, trivial quantum principal bundles are perhaps the first examples that come to mind.  We consider the theory of trivial qpb's as developed in \cite{micho2,stheve}. In Durdevich's 
formulation, a quantum principal $\G$--bundle of the form $$\zeta^\triv=(P:=B\otimes H,B,\Delta_P:=\id_B\otimes \Delta)$$
is called {\it trivial} \cite{micho2,stheve}. In the previous section, we have already worked on these kinds of qpb's, specifically in Examples \ref{e.1.1}, \ref{e.1.2.1} and Proposition \ref{noteq}. In this way, in this subsection we will present some generalizations. First of all, we have

\begin{Proposition}
    \label{gentrv1}
    Let $\T$ be a complete set of mutually non--equivalent irreducible $\G$-- corepresentations with $\delta^\C_\triv$ $\in$ $\T$. If $\delta^V$ $\in$ $\T$ coacts on a $\C$--vector space of dimension $n_{V}$, then there exists a left/right $B$--basis $$\{T^{\l}_k \}^{n_{V}}_{k=1} \;\subseteq\; \Mor(\delta^V,\Delta_P)$$ such that equation (\ref{generators}) holds. In particular, associated left/right qvb's are always free modules (by the Serre--Swan theorem, they can be considered as trivial qvb's) for every $\delta^V$ $\in$ $\T$. 
\end{Proposition}
\begin{proof}
    Let $\{g^{V}_{ij}\}^{n_{V}}_{i,j=1}$ be the linear basis of $H$ given in Theorem \ref{rep}. Consider the matrix $H^{V}=(g^{V}_{ij})$ $\in$ $M_{n_{V}}(H)$. Then the linear maps $$T^{\l}_k : V\longrightarrow B\otimes H$$ defined by $$T^{\l}_k(e_i)=\mathbbm{1}_B\otimes g^{V}_{ki}=:x^V_{ki},$$ form a left/right $B$--basis of $\Mor(\delta^V,\Delta_P)$ and they satisfy equation (\ref{generators}), because in accordance with equation (\ref{2.f9}), we have $H^{V\,\dagger}H^{V}= \Id_{n_{V}}.$
\end{proof}

The last proposition shows the specific form of the maps  $\{T^{\l}_k \}^{n_{V}}_{k=1}$ for every trivial qpb. Notice that Proposition \ref{gentrv1} also implies that $\varrho^V(\mathbbm{1})=\Id_{n_V}$ (see equation (\ref{3.f2})) and
hence, we explicitly obtain that the canonical Hermitian structures  on associated left/right qvb's are non--degenerate (see the proof of Proposition \ref{caher}). 

In Durdevich's formulation, if a differential calculus on $\zeta^\triv$ is given by (using the corresponding tensor products) $$\Omega^\bullet(P)=\Omega^\bullet(B)\otimes \Gamma^\wedge\quad  \mbox{ and }\quad \Delta_{\Omega^\bullet(P)}:=\id_{\Omega^\bullet(B)}\otimes \Delta, $$ where $\Omega^\bullet(B)$ is some graded differential $\ast$--algebra generated by $\Omega^0(B)=B$ and $\Gamma^\wedge$ is the universal differential envelope $\ast$--calculus of some bicovariant $\ast$--FODC of $\G$, then the differential calculus is commonly referred to as \emph{trivial}.

Let us take any trivial qpb $\zeta^\triv=(P,B,\Delta)$ with a trivial differential calculus. Then, the linear map
\begin{equation}
\label{5.f1.1}
\omega^\triv:\mathfrak{qg}^\# \longrightarrow \Omega^1(P),\qquad
\theta \longmapsto \mathbbm{1}_B\otimes \theta.
\end{equation}
is a real, regular and multiplicative qpc (see equations (\ref{2.f24.6}), (\ref{2.f25}), (\ref{2.f26})) and it is referred to as the {\it trivial qpc}.  According to \cite{micho2,stheve}, there is a bijection between $$\mathrm{Hom}(\mathfrak{qg}^\#,\Omega^1(B))=\{ A:\mathfrak{qg}^\#\longrightarrow \Omega^1(B)\mid  A \mbox{ is linear} \} $$ and the set of all qpc's of $\zeta^\triv$. This bijection is based on the fact that every qpc $\omega$
can be uniquely written in the form
\begin{equation}
\label{5.f1}
\omega=(A\otimes \id_H)\circ \ad+\omega^\triv,
\end{equation}
where $A:\mathfrak{qg}^\#\longrightarrow \Omega^1(B)$ is a linear map. The map $A$ can be interpreted as the {\it non--commutative gauge potential} of $\omega$. This bijection extends naturally to the curvature by  (\cite{micho2}) $$R^{\omega}=(F\otimes \id_H)\circ \ad,$$ where $F:\mathfrak{qg}^\#\longrightarrow \Omega^2(B)$ is the linear map defined as  $F:=dA-\langle A,A\rangle.$  The map $F$ can be interpreted as the {\it non--commutative field strength} of $\omega$.

\begin{Proposition}
\label{4.10}
Let $\delta^V$ $\in$ $\T$. Continuing with Proposition \ref{gentrv1} we have 
\begin{enumerate}
\item For the trivial qpc, the covariant derivatives $D^{\omega^\triv}$ (see equation (\ref{2.f30})) and $\widehat{D}^{\omega^\triv}$ (see equation (\ref{2.f30.1}))   satisfy   $D^{\omega^\triv}\circ T^\l_k=0$, $\widehat{D}^{\omega^\triv}\circ T^\l_k=0$  for all $k$.  
\item Taking into account equation (\ref{3.f5}), the induced qlc's for $\omega^\triv$  can be expressed by 
\begin{equation}
\label{f79}
\nabla^{\omega^\triv}_{V}(T)=\sum^{n_{V}}_{k=1}db^{_T}_k\otimes_B  T^\l_k\,,\quad \widehat{\nabla}^{\omega^\triv}_{V}(T)=\sum^{n_{V}}_{k=1}T^\l_k\otimes_B db^{_T}_k 
\end{equation}
\item The exterior covariant derivatives are given by $$d^{\nabla^{\omega^\triv}_{V}}(\mu\otimes_B T)=\sum^{n_{V}}_{k=1}d(\mu\,b^{_T}_k)\otimes_B T^\l_k,\quad \mbox{ and }\quad d^{\widehat{\nabla}^{\omega^\triv}_{V}}(T\otimes_B\mu)=\sum^{n_{V}}_{k=1}T^\l_k\otimes_B d(\mu\,p^{_T}_k)$$
for all $\mu$ $\in$ $\Omega^\bullet(B)$. In particular $R^{\nabla^{\omega^\triv}_{V}}=0$, $R^{\widehat{\nabla}^{\omega^\triv}_{V}}=0$.
\end{enumerate}
\end{Proposition}

\begin{proof}
\begin{enumerate}
\item Because of $dg=g^{(1)}\pi(g^{(2)})$ for all $g$ $\in$ $H$, we get
\begin{eqnarray*}
D^{\omega^\triv}(T^\l_k(e_i))\,=\, D^{\omega^\triv}(\mathbbm{1}_B\otimes g^{V}_{ki}) &= &  \mathbbm{1}_B\otimes dg^{V}_{ki}-(\mathbbm{1}_B\otimes g^{V\,(1)}_{ki})\,\omega^\triv(\pi(g^{V\,(2)}_{ki}))
  \\
  &= &
 \mathbbm{1}_B\otimes dg^{V}_{ki}-\mathbbm{1}_B\otimes g^{V\,(1)}_{ki}\pi(g^{V\,(2)}_{ki})\;=\;0.
\end{eqnarray*}
This shows that $D^{\omega^\triv}\circ T^\l_k=0$ for $k=1,...,n_{V}$. A similar calculation proves that $\widehat{D}^{\omega^\triv}\circ T^\l_k=0$.
\item By equation (\ref{3.f5}) we obtain $$T=\displaystyle\sum^{n_{V}}_{k=1}b^{_T}_k\,T^\l_k=\sum^{n_{V}}_{k=1}T^\l_k\,b^{_T}_k$$ for all $T$ $\in$ $\Mor(\delta^V,\Delta_P)$, where $\delta^V$ $\in$ $\T$. Part (1) of this proposition together with equations  (\ref{2.f32}), (\ref{2.f32.1}) show that $$D^{\omega^\triv}\circ T=\displaystyle\sum^{n_{V}}_{k=1}db^{_T}_k\, T^\l_k,\qquad \widehat{D}^{\omega^\triv}\circ T=\displaystyle\sum^{n_{V}}_{k=1} T^\l_k\,db^{_T}_k.$$ Hence equation (\ref{f79}) follows directly from equations (\ref{3.f8}), (\ref{3.f8.1}).
\item This follows directly from part (2) of this proposition and equations (\ref{3.f10.1}), (\ref{3.f10.4}).
\end{enumerate}
\end{proof}

Since $$T=\displaystyle\sum^{n_{V}}_{k=1}b^{_T}_k\,T^\l_k=\sum^{n_{V}}_{k=1}T^\l_k\,b^{_T}_k$$ for all $T$ $\in$ $\Mor(\delta^V,\Delta_P)$, the canonical Hermitian structures are given by $$\langle T_1,T_2\rangle_\l=\sum^{n_v}_{k=1}b^{_{T_1}}_k\,(b^{_{T_2}}_k)^\ast\quad\mbox{ and }\quad \langle T_1,T_2\rangle_\r =\sum^{n_v}_{k=1}(b^{_{T_1}}_k)^\ast\, b^{_{T_2}}_k.$$ Therefore, explicitly we have $$\langle \nabla^{\omega^\triv}_{V}(T_1),T_2\rangle_\l+ \langle T_1,\nabla^{\omega^\triv}_{V}(T_2)\rangle_\l=\sum^{n_v}_{k=1}db^{_{T_1}}_k\, (b^{_{T_2}}_k)^\ast+b^{_{T_1}}_k\,d (b^{_{T_2}}_k)^\ast=d\langle T_1,T_2\rangle_\l$$ and $$\langle \widehat{\nabla}^{\omega^\triv}_{V}(T_1),T_2\rangle_\r+ \langle T_1,\widehat{\nabla}^{\omega^\triv}_{V}(T_2)\rangle_\r=\sum^{n_v}_{k=1}d(b^{_{T_1}}_k)^\ast\, b^{_{T_2}}_k+(b^{_{T_1}}_k)^\ast\,db^{_{T_2}}_k=d\langle T_1,T_2\rangle_\r.$$

Since every finite--dimensional $\G$--corepresentation is the direct sum of a finite number of elements of $\T$, Propositions \ref{gentrv1}, \ref{4.10} and the last two equalities are naturally generalized to every $\delta^V$ $\in$ $\FD(\Rep_{\G})$ and they show explicitly, the behavior of every trivial qpb with a trivial differential calculus under the geometric structures introduced in Section $3$ for $\omega^\triv$. 

For degree zero elements, we have (see equation (\ref{qtrs0})) $$\qtrs(g^V_{ij})=\sum^{n_V}_{k=1}(\mathbbm{1}\otimes g^{V\,\ast}_{ki})\otimes_B (\mathbbm{1}\otimes g^{V}_{kj}).$$ Furthermore, considering $\omega^\triv$ in equation (\ref{6.f1.6.1}), we obtain $$\qtrs(\theta)=\mathbbm{1}_P\otimes_{\Omega^\bullet(B)}(\mathbbm{1}_B\otimes \theta)-\omega(\theta^{(0)})[\theta^{(1)}]_1\otimes_{\Omega^\bullet(B)} [\theta^{(1)}]_1$$ with $\ad(\theta)=\theta^{(0)}\otimes \theta^{(1)}.$

 It is worth remembering that in Examples \ref{e.1.1}, \ref{e.1.2.1} and Proposition \ref{noteq}, we have provided an explicit example of the form of the quantum gauge group and its action on the space of qpc's in a trivial qpb with a trivial differential calculus. In the context of these examples, since $$\ad(\theta)=\theta\otimes \mathbbm{1}$$ for all $\theta$ $\in$ $\mathfrak{qg}^\#$, equation (\ref{5.f1}) implies that the action of $\qGG$ on the space of qpc's is transitive (see equation (\ref{qpctrz1})). In other words, and according to Remark \ref{remaaction}, the action of $\qGG$ on qpc's via  $\F^\times\omega$ can be trivial; while the action of $\qGG$ on qpc's via  $\F^\circledast \omega$ is always transitive. 
 
 We decided to use trivial qpb's with $\G$ the quantum group associated with $\Z_2$ in our concrete examples because their simplicity allows us to clearly show the novel aspects of the theory and highlight its differences from other approaches (see, for example, \cite{landi,ha}). Moreover, their simplicity enables the reader to easily follow all the calculations.

 \subsection{Homogeneous Quantum Principal Bundles}

Homogeneous quantum principal bundles are one of the most well--studied examples of qpb's and the reader can check the basics in, for example, \cite{micho2,stheve,libro}. Let  $\G$ be a quantum group with dense $\ast$--Hopf algebra $$(P,\cdot,\mathbbm{1},\Delta,\epsilon,S,\ast)$$ and $\mathcal{H}$ be a quantum subgroup with dense $\ast$--Hopf algebra $$(H,\cdot,\mathbbm{1},\Delta',\epsilon',S',\ast).$$ The quantum subgroup structure implies the existence of a $\ast$--Hopf algebra epimorphism $$j:P\longrightarrow H$$ and consider the linear map $$\Delta_P:=(\id_P\otimes j)\circ \Delta: P\longrightarrow P\otimes H.$$ Defining $$B:=\{b \in P \mid \Delta_P(b)=b\otimes \mathbbm{1} \},$$ the triple $$\zeta=(P,B,\Delta_P)$$ is a qpb called {\it homogeneous quantum principal bundle} \cite{micho2,stheve,libro}.

The following proposition shows that our theory can be applied to this class of quantum bundles. 

\begin{Proposition}
    \label{homo}
    Let $\zeta=(P,B,\Delta_P)$  be a homogeneous qpb. Then, for every $\delta^V$ $\in$ $\T'$, there exists $$\{ T^V_k\}^{d_V}_{k=1}\;\subseteq\;\Mor(\delta^V,\Delta_P)$$ such that equation (\ref{generators}) holds. In addition, we have 
    \begin{equation}
        \label{generatorshomo}
        \sum^{d_V}_{k=1} x^V_{ki}\,x^{V\,\ast}_{kj}=\delta_{ij}\mathbbm{1}.
    \end{equation}
    Here, $\T'$ is a complete set of mutually non--equivalent irreducible $\mathcal{H}$--corepresentation with $\delta^\C_\triv$ $\in$ $\T'$.
\end{Proposition}

\begin{proof}
    Let $\delta^{V}$ $\in$ $\T'$  and let $\T$ be a complete set of mutually non--equivalent irreducible $\G$--corepresentation. For any $\delta^W_P$ $\in$ $\T$, consider the elements $\{ p^W_{il}\}^{n_W}_{i,l=1}$ $\subseteq$ $P$ of equation (\ref{2.f8}). Also notice that $(\id_W\otimes j)\circ \delta^W_P$ is always a $\mathcal{H}$--corepresentation, but it can be a reducible corepresentation. 
    
    Since $$j:P\longrightarrow H$$ is the $\ast$--Hopf algebra epimorphism of the corresponding quantum subgroup structure, there exists $\delta^W_P$ $\in$ $\T$ such that  
    $$(\id_W\otimes j)\circ\delta^W_P\cong \bigoplus^m_{r=1}\delta^{V_r},$$ where $\delta^{V_1},...,\delta^{V_s}=\delta^{V},...., \delta^{V_m}$ $\in$ $\T'$ for some $s$ $\in$ $\{1,...,m\}$.  Without loss of generality, assume that $s=1$ and the indices $i$ and $l$ are ordered in such a way that $$j(p^{W}_{il})=h^{V}_{il}$$ for $i,l=1,...,{\mathrm{dim}}_\C(V)=n_{V}$ and for the other indices $$j(p^{W}_{il})=0,$$ where $\{ h^{V}_{il}\}^{n_{V}}_{i,l=1}$ $\subseteq$ $H$ are the elements of equation (\ref{2.f8}) for $\delta^{V}$. In this way, we define the linear maps $$T^\l_k:V=V^s \longrightarrow P$$ such that $$T^\l_k(e_l)=p^{W}_{kl}=:x^{V}_{kl}$$ for $k$ $\in$ $\{1,...,d_V=n_W\}$. By equation (\ref{2.f9}), it is easy to verify that $T^\l_i$ $\in$ $\Mor(\delta^{V},\Delta_P)$ and that the equation (\ref{generators}) is satisfied. Also, by equation (\ref{2.f9}) it is easy to verify that equation (\ref{generatorshomo}) is fulfilled. 
\end{proof}

It is worth mentioning that last proposition shows the specific form of the generators $\{T^\l_s \}^{d_{V}}_{s=1}$ for homogeneous qpb's. Furthermore, since the elements $\{x^{V}_{si}\}$ are columns of the matrix $P^W=(p^W_{ij})$, it follows that $\varrho^V(\mathbbm{1})=\Id_{n_{W}}$ and
therefore, we explicitly obtain that the canonical Hermitian structures on the associated left/right qvb's are non--degenerate (see the proof of Proposition \ref{caher}). Moreover, this proves that associated left/right qvb's are trivial (they are free left/right modules).

For homogeneous qpb's, $P$ and $H$ are derived from quantum groups. Therefore, quantum differential forms of $P$ and quantum differential forms of $H$ will be given by the universal differential envelope $\ast$--calculus corresponding to some $\ast$--FODC of $P$ and some bicovariant $\ast$--FODC of $H$. Following this idea, reference \cite{micho2} presents a method to construct a differential calculus on homogeneous qpb's that we can use to apply our theory.

Let us focus in a concrete example. Let $\G$ be the quantum  $SU(2)$ group and $\mathcal{H}$ the quantum group associated with the Lie group $U(1)$.  It is worth mentioning that we are going to use the original Woronowicz's notation presented in \cite{woro1} for the quantum  $SU(2)$ group. In other words, the dense $\ast$--Hopf algebra of the quantum  $SU(2)$ group consists of the $\ast$--algebra $$(P:=\mathrm{SU}_q(2),\cdot,\mathbbm{1},\ast)$$ generated by two symbols $\{\alpha,\gamma\}$ satisfying
\begin{equation}
\label{2.f1.6}
\begin{aligned}
\alpha^{\ast}\alpha+\gamma^{\ast}\gamma=\mathbbm{1},\qquad \alpha\alpha^{\ast}+q^{2}\gamma\gamma^{\ast}=\mathbbm{1},\qquad \gamma\gamma^{\ast}=\gamma^{\ast}\gamma\\
q\gamma\alpha=\alpha\gamma,\quad q\alpha^\ast\gamma^\ast=\gamma^\ast\alpha^\ast, \quad q\gamma^{\ast}\alpha=\alpha\gamma^{\ast},\quad q\alpha^\ast\gamma=\gamma\alpha^\ast
\end{aligned}
\end{equation}
 and the following relations for the coproduct, the counit and the antipode:   
\begin{equation}
\label{2.f1.7}
\begin{split}
\Delta(\alpha)&=\alpha\otimes\alpha-q\gamma^\ast\otimes\gamma, \quad \Delta(\gamma)=\gamma\otimes\alpha+\alpha^\ast\otimes\gamma,\quad \epsilon(\alpha)=1, \quad \epsilon(\gamma)=0,\\
&S(\alpha)=\alpha^\ast, \qquad S(\alpha^\ast)=\alpha,\qquad S(\gamma)=-q\gamma,\qquad S(\gamma^\ast)=-q^{-1}\gamma^\ast.
\end{split}
\end{equation}
In the same way, let us take the quantum group $\mathcal{U}(1)$ naturally associated to the Lie group $\mathrm{U}(1)$. Thus, its dense $\ast$--Hopf algebra is given by the Laurent polynomial algebra, i.e., $$(H:=\C[z,z^\ast]=\C[z,z^{-1}],\cdot,\mathbbm{1},\ast)$$ and following the relations for the coproduct, the counit and the antipode:
\begin{equation}
\label{3.f4.1.1}
\Delta'(z)=z\otimes z,\quad \epsilon'(z)=1, \quad S'(z)=z^\ast, \quad S'(z^\ast)=z.
\end{equation}
Notice that this algebra is commutative.

The homogeneous qpb $$\zeta=(P,B,\Delta_P)$$   for the $\ast$--Hopf algebra epimorphism  $$j: P\longrightarrow H$$   given by $$j(\alpha)=z, \quad j(\gamma)=0,$$ is called the $q$--Dirac monopole bundle or the quantum Hopf fibration \cite{micho2,sald3}. In differential geometry, the Hopf fibration is one of the most basic and well–established examples of principal bundles. We will use the {\it non--commutative geometrical} version of this bundle to illustrate our theory.

It is well--known that any complete set of mutually non--equivalent irreducible $\mathcal{H}$-- corepresentations $\T'$ is  in bijection with $\Z$. In addition, for every $n$ $\in$ $\Z$,  the corresponding corepresentation is $1$--dimensional.  Then, for $n$ $\in$ $\N_0=\N\cup \{0\}$ and according to the previous proposition, the elements $\{x^V_{i+1\,1}\}^{n}_{i=0}$ are exactly the first column of the corepresentation matrix of $P$ for spin $l=n/2$ \cite{sald3,su(2)}. On the other hand, for $n< 0$, the elements $\{x^{V}_{i+1\,1}\}^{|n|}_{i=0}$ given in Proposition (\ref{homo}), without assuming $s=1$, are exactly the last column of the corepresentation matrix of $P$ for spin $l=|n|/2$ \cite{sald3,su(2)}.

To obtain a differential calculus on $\zeta$, let us consider the Woronowicz $3D$ $\ast$--FODC $(\Psi,d)$ of $P$ \cite{woro2}. This $\ast$--FODC is left covariant (see equation (\ref{2.f9.2})) and it is defined by the right $P$--ideal $\mathcal{R} \subseteq \Ker(\epsilon)$ generated by $$A=\{ \gamma^2, \, \gamma^{\ast\,2},\,\gamma\gamma^\ast,\, \alpha\gamma-\gamma,\,\alpha\gamma^\ast-\gamma^\ast,\,q^2\alpha+\alpha^\ast-(1+q^2)\mathbbm{1}\}.$$ The $\C$--vector space $${_\inv}\Psi={\Ker(\epsilon)\over \mathcal{R}}$$ is $3$--dimensional, $$\beta:=\{\eta_3:=\pi(\alpha-\alpha^\ast),\;\eta_+=\pi(\gamma),\;\eta_-=\pi(\gamma^\ast) \}$$ is a linear basis of ${_\inv}\Psi$ and it is also a left $P$--basis of $\Psi$ (\cite{stheve}), where $$\pi:P\longrightarrow {_\inv}\Psi$$ is the quantum germs map of $(\Psi,d)$ (see equation (\ref{2.f14})) \cite{woro2,stheve}. In this way, the universal differential envelope $\ast$--calculus $$(\Psi^\wedge,d,\ast)$$ of $(\Psi,d)$  will be the space of quantum differential forms of $P$.

On the other hand, consider the right $H$--ideal $$\mathcal{R}'\subseteq \Ker(\epsilon')$$ generated by $$j(A)=\{q^2z+z^\ast-(1+q^2)\mathbbm{1} \}.$$ Thus, $\mathcal{R}'$ satisfies $S'(\mathcal{R}')^\ast\subseteq \mathcal{R}'$ and $\Ad'(\mathcal{R}')\subseteq \mathcal{R}'\otimes H$. According to Proposition \ref{0U}, $\mathcal{R}'$ defines a bicovariant $\ast$--FODC $$(\Gamma,d)$$ of $\mathcal{H}$. It is worth mentioning that $$\mathfrak{qh}^\#:={\Ker(\epsilon')\over \mathcal{R}'}=\mathrm{span}_\C \{\varsigma:= \pi'(z-z^\ast)\},$$ where $$\pi':H\longrightarrow \mathfrak{qh}^\#$$ is the quantum germs map of $(\Gamma,d)$ and $\beta'=\{\varsigma\}$ is also a left $H$--basis of $\Gamma$ \cite{woro2,stheve}. In this way, the universal differential envelope $\ast$--calculus $$(\Gamma^\wedge,d,\ast)$$ of $(\Gamma,d)$ will be the space of quantum differential forms of $H$ and it is worth mentioning that $\Gamma^\wedge$ does not have elements for $n \geq 2$. Additionally, it is not graded--commutative because (\cite{micho2}) $$\pi'(z)\,z=q^2 z\,\pi'(z).$$  

Finally, define the map $$\Delta_{\Psi^\wedge}|_{P}:= \Delta_P;$$ and on $\Psi$ set $$\Delta_{\Psi^\wedge}|_{\Psi}(p\,\eta_+)=\Delta_P(p)\,(\eta_+\otimes z^2),\qquad  \Delta_{\Psi^\wedge}|_{\Psi}(p\,\eta_-)=\Delta_P(p)\,(\eta_+\otimes z^{\ast 2}),$$ $$\Delta_{\Psi^\wedge}|_{\Psi}(p\,\eta_3)=\Delta_P(p)\,(\eta_3\otimes \mathbbm{1}+\mathbbm{1}\otimes \varsigma)$$ for all $p$ $\in$ $P$. This map can be extended to a graded differential $\ast$--algebra morphism  $$\Delta_{\Psi^{\wedge}}: \Psi^{\wedge}\longrightarrow \Psi^{\wedge} \otimes \Gamma^\wedge$$ thereby endowing $\zeta$ with a differential calculus \cite{micho2,sald3}.

There is a canonical qpc $$\omega^c:\mathfrak{qh}^\#  \longrightarrow \Psi $$ given by $$\omega^c(\varsigma)=\eta_3.$$ This qpc is the {\it non--commutative geometrical} counterpart of the principal connection on the Hopf fibration associated with the Levi--Civitta connection \cite{sald3}. 

Let us compute the induced qlc's. Considering that $\T' \longleftrightarrow \Z$, take $n$ $\in$ $\N$. For  $T$ $\in$ $\Mor(n,\Delta_P)$ such that $T(1)=\alpha^n$, by definition we have 
\begin{equation}
    \label{leftalpha}
    \nabla^{\omega^\c}_n(T)=-{(1-q^{2n})q^{3-2n}\over 1-q^2}\sum^n_{k=0}\alpha^{n-1}\gamma^\ast \eta_+ x^{V\,\ast}_{k+1\,1}\otimes_{B} T^\l_{k+1\,1},
\end{equation}
\begin{equation}
    \label{rightalpha}
    \widehat{\nabla}^{\omega^\c}_n(T)=-{(1-q^{2n})q^{3-2n}\over 1-q^2}\sum^n_{k=0} q^{2k}\,T^\l_{k+1\,1}\otimes_{B} x^{V\,\ast}_{k+1\,1} \alpha^{n-1}\gamma^\ast \eta_+.
\end{equation}

The left and right canonical Hermitian structure are given by $$\langle T_1,T_2\rangle_\l=T_1(1)T_2(1)^\ast\,, \quad \langle T_1,T_2\rangle_\r=T_1(1)^\ast T_2(1)$$ and, for example $$\langle T,T\rangle_\l=\alpha^n \alpha^{\ast\,n},\qquad \langle T,T\rangle_\r=\alpha^{\ast\,n} \alpha^n.$$ Moreover, by using equation (\ref{generators}) and the commutativity relations of $P$  and its calculus, we get (\cite{micho2,woro1})
\begin{eqnarray*}
\langle \nabla^{\omega^\c}_{n} (T),T \rangle_\l &= & -{(1-q^{2n})q^{3-2n}\over 1-q^2}\sum^n_{k=0} \alpha^{n-1}\gamma^\ast \eta_+ x^{V\,\ast}_{k+1\,1} \langle T^n_{k+1\,1},T \rangle_\l
  \\
  &= &
-{(1-q^{2n})q^{3-2n}\over 1-q^2}\sum^n_{k=0} \alpha^{n-1}\gamma^\ast \eta_+ x^{V\,\ast}_{k+1\,1} x^{V}_{k+1\,1} \alpha^{\ast\,n}
  \\
  &= &
-{(1-q^{2n})q^{3}\over 1-q^2} \alpha^{n-1}\alpha^{\ast\,n}\gamma^\ast\eta_+
\end{eqnarray*}
and  
\begin{eqnarray*}
\langle T,\nabla^{\omega^\c}_{n} (T) \rangle_\l  &= &  -{(1-q^{2n})q^{3-2n}\over 1-q^2}\sum^n_{k=0} \langle T,T^n_{k+1\,1}\rangle_\l\, x^{V}_{k+1\,1}(\alpha^{n-1}\gamma^\ast \eta_+)^\ast
  \\
  &= &
-{(1-q^{2n})q^{3-2n}\over 1-q^2}\sum^n_{k=0} \alpha^n x^{V\,\ast}_{k+1\,1} \,x^{V}_{k+1\,1}(\alpha^{n-1}\gamma^\ast \eta_+)^\ast
  \\
  &= &
-{(1-q^{2n})q \over 1-q^2} \alpha^n\alpha^{\ast\,n-1}\gamma \eta_-;
\end{eqnarray*}
which shows that $$\langle \nabla^{\omega^\c}_{n} (T),T \rangle_\l+ \langle T,\nabla^{\omega^\c}_{n} (T) \rangle_\l=d(\alpha^n \alpha^{\ast\,n})=d\langle T,T\rangle_\l.$$

In the same way,  we obtain 
\begin{eqnarray*}
\langle \widehat{\nabla}^{\omega^\c}_{n} (T),T \rangle_\r &= & 
-{(1-q^{2n})q^{3-2n}\over 1-q^2}\sum^n_{k=0} q^{2k} (\alpha^{n-1}\gamma^\ast \eta_+)^{\ast} x^{V}_{k+1\,1} \langle T^n_{k+1\,1},T \rangle_\r
  \\
  &= &
-{(1-q^{2n})q^{3-2n}\over 1-q^2}\sum^n_{k=0} q^{2k} (\alpha^{n-1}\gamma^\ast \eta_+)^{\ast} x^{V}_{k+1\,1} x^{V\,\ast}_{k+1\,1}\alpha^n
  \\
  &= &
-{(1-q^{2n})q^{1-2n}\over 1-q^2} \alpha^{\ast\,n-1}\alpha^n\gamma\eta_-,
\end{eqnarray*}
and 
\begin{eqnarray*}
\langle T,\widehat{\nabla}^{\omega^\c}_{n} (T) \rangle_\r  &= &  -{(1-q^{2n})q^{3-2n}\over 1-q^2}\sum^n_{k=0} q^{2k} \langle T,T^n_{k+1\,1}\rangle_\r\, x^{V\,\ast}_{k+1\,1}\alpha^{n-1}\gamma^\ast \eta_+
  \\
  &= &
-{(1-q^{2n})q^{3-2n}\over 1-q^2}\sum^n_{k=0} q^{2k} \alpha^{\ast\,n}  x^{V}_{k+1\,1}x^{V\,\ast}_{k+1\,1}\alpha^{n-1}\gamma^\ast \eta_+
  \\
  &= &
-{(1-q^{2n})q^{3-2n} \over 1-q^2} \alpha^{\ast\,n} \alpha^{n-1}\gamma^\ast \eta_.
\end{eqnarray*}
Hence $$\langle \widehat{\nabla}^{\omega^\c}_{n} (T),T \rangle_\r+ \langle T,\widehat{\nabla}^{\omega^\c}_{n} (T) \rangle_\r=d(\alpha^{\ast\,n}\alpha^n)=d\langle T,T\rangle_\r.$$

For degree zero elements, we have (see equation (\ref{qtrs0})) $$\qtrs(z^n)=\sum^n_{k=0} 
\begin{bmatrix}
n \\
k
\end{bmatrix}_{q^{-2}} \gamma^{\ast\,k}\alpha^{\ast\,n-k} \otimes_{B} \alpha^{n-k} \gamma^k.$$ Moreover, considering $\omega^c$ in equation (\ref{6.f1.6.1}), it follows that $$\qtrs(\varsigma)=\mathbbm{1}\otimes_{\Omega^\bullet(B)} \eta_3-\eta_3\otimes_{\Omega^\bullet(B)}\mathbbm{1}.$$

Due to the way we have defined the quantum gauge group, it is very large and, in general, quite challenging to calculate its explicit form. However, it is possible to prove the following proposition.

\begin{Proposition}
    \label{qhfqgg}
    In the quantum Hopf fibration with the differential calculus defined above, the action of $\qGG$ on the space $\mathfrak{qpc}(\zeta)$ is transitive. 
\end{Proposition}

\begin{proof}
    Let $\omega$ $\in$ $\mathfrak{qpc}(\zeta)$. According to \cite{micho2, stheve}, $\mathfrak{qpc}(\zeta)$ is an affin space modeled by  $\overrightarrow{\mathfrak{qpc}(\zeta)}$ (see equation (\ref{2.f24.2})). Consequently, there exists $\lambda$ $\in$ $\overrightarrow{\mathfrak{qpc}(\zeta)}$ such that $$\omega=\omega^c+\lambda.$$ By equations (\ref{2.f15.1}), (\ref{3.f4.1.1}),  we get
    \begin{equation}
        \label{porqueyo}
        \ad(\varsigma)=(\pi'\otimes \id_{H})\Ad(z-z^\ast)=\pi'(z-z^\ast)\otimes \mathbbm{1}=\varsigma\otimes \mathbbm{1},
    \end{equation}
     and it follows that $$\lambda(\varsigma)=:\mu\;\in\;\Omega^1(B)=\{\mu \in \Psi \mid \Delta_{\Psi^\wedge}(\mu)=\mu\otimes\mathbbm{1} \}.$$ Hence, we have $$\omega(\varsigma)=\eta_3+\mu.$$

    On the other hand, consider the graded linear map $$\f: \Gamma^\wedge\longrightarrow \Psi^\wedge$$  defined by $$\f|_H=\epsilon'\, \mathbbm{1}\quad \mbox{ and }\quad \f|_{\Gamma}(h\varsigma)=\epsilon'(h)\mu.$$ Since $\Ad(\vartheta)=\vartheta\otimes \mathbbm{1}$ for all $\vartheta$ $\in$ $\Gamma^\wedge$, a straightforward calculation shows that $\f$ satisfies equation (\ref{4.f3}). Moreover, $\f$ is a convolution invertible map with convolution inverse $$\f^{-1}: \Gamma^\wedge\longrightarrow \Psi^\wedge$$ given by $$\f^{-1}|_H=\epsilon'\, \mathbbm{1} \quad \mbox{ and } \quad \f^{-1}|_{\Gamma}(h\varsigma)=-\epsilon'(h)\mu.$$ Thus, by Proposition \ref{4.1} we obtain a  qgt $$\F_\f:\Psi^\wedge\longrightarrow \Psi^\wedge $$ such that $$\F_\f|_P=\id_P \quad \mbox{ and } \quad \F_\f(\eta_3)=\eta_3+\mu .$$ According to Proposition \ref{4.4} and equation (\ref{porqueyo}), we have
    \begin{eqnarray*}
        \F^{\circledast}_\f\omega^\c(\varsigma)=m_{\Omega}(\omega^\c\otimes \f )\ad(\varsigma) +\f(\varsigma)
        &=&
        m_{\Omega}(\omega^\c\otimes \f )(\varsigma\otimes \mathbbm{1}) +\f(\varsigma)
         \\
        &=&
\omega^\c(\varsigma)\f(\mathbbm{1})+f(\varsigma)
        \\
        &=&
        \mathbbm{1}\otimes \varsigma+\mu\otimes \mathbbm{1}
        \\
        &=&
        \omega(\varsigma).
    \end{eqnarray*}
    We conclude that the action is transitive.
\end{proof}

As we have mentioned in the previous section, it is natural to work with {\it ad hoc} subgroups of $\qGG$ in each situation, as in \cite{landi}. As another example, in \cite{sald3} we work with the subgroup $\qGG_{\mathrm{YM}}$ of $\qGG$ that leaves invariant the {\it non--commutative geometrical} Yang--Mills functional in the quantum Hopf fibration with the differential calculus introduced above.  In this way, by Proposition \ref{qhfqgg} one can conclude that
$$\qGG_{\mathrm{YM}}=\{\F \in \qGG \mid \F^\circledast \omega=\omega^\c+\lambda\;\mbox{ with }\; d\lambda(\varsigma)=0\}.$$ Since $\omega^\c$ is a critical point of the non--commutative geometrical Yang--Mills functional, the last characterization of $\qGG_{\mathrm{YM}}$ implies that, up elements of this subgroup, $\omega^\c$ is the unique Yang--Mills qpc, exactly as in the {\it classical} case \cite{sald3}.

\section{Concluding Comments}

This paper extends the work presented in \cite{sald1} by considering general qpc's rather than only the real and regular ones, and by defining additional geometrical structures, which we will review.

First of all, we would like to highlight the importance of the universal differential envelope $\ast$--calculus $$(\Gamma^\wedge,d,\ast)$$ as quantum differential forms on $\G$. This space not only allows us to extend the $\ast$--Hopf algebra structure of $H^\infty$ to $\Gamma^{\wedge\,\infty}$, but it is also maximal with this property \cite{micho1}; moreover, it generalizes the $\ast$--algebra of $\C$--valued differential forms of a (compact matrix) Lie group (\cite{appendix}). In addition, $(\Gamma^\wedge,d,\ast)$ allows us to define the quantum translation map at the level of differential calculus, which in turn leads to Proposition \ref{4.1}. It is worth remembering that not all conditions in $H^\infty$ can be extended to $\Gamma^{\wedge\,\infty}$ (for example, see Example \ref{e.1}). This opens the door to an exciting research project in which one could explore ways to generalize those properties that do not naturally extend, so that they encompass $\Gamma^{\wedge\,\infty}$ in each situation.

Since $\Mor(\delta^V,\Delta_P)$ is a $B$--bimodule in a natural way, we decided to deal with the left and right structures, and Durdevich's theory allows us to develop the theory for the left/right associated qvb's using the covariant derivatives $$D^\omega \quad \mbox{ and }\quad \widehat{D}^\omega.$$ In \cite{sald2, sald3, sald4} one can appreciate more explicitly the importance of taking into account both associated qvb's and their induced qlc's. For smooth compact manifolds, both associated qvb's are the same, and since every qpc that comes from the {\it dualization} of a {\it classical} principal connection is regular and real, both induced qlc's are the same. 

It is worth mentioning that there are other papers dealing with Hermitian structures on quantum spaces, for example \cite{Heck}, in which the author presented a notion of spin geometry on quantum groups. In \cite{Heck}, quantum differential forms of $\G$ are given by the braided
exterior calculus (\cite{woro2}) instead of the universal differential envelope $\ast$--calculus that we used. Nevertheless, there is a surjectived morphism between these two spaces (\cite{micho2}), and with that one could try to integrate both ideas in order to develop a theory for spinor quantum bundles.

The main two reasons to use the word {\it canonical} in Definitions \ref{lhs}, \ref{rhs} are the facts that at the end, these $B$--valued inner products do not depend on the choice of the generators $\{T^\l_k\}$ and of course, the other reason is the result presented in Theorem \ref{fgs}.  

As we have mentioned at the end of Section $3$, this theorem is the core of this paper since it recreates an important {\it classical} result in the most general framework of non--commutative geometry, which will allow to define formally adjoint operators of $d^{\nabla^\omega_V}$, $d^{\widehat{\nabla}^\omega_V} $ and with that, we will able to define Laplacians for associated qvb's and study field theory on them, like the reader can check in \cite{sald2}. In addition, we want to emphasize Theorems \ref{hil}, \ref{hil1}. These theorems show a link between associated qvb's and the well--known theory of Hilbert $C^\ast$--modules by the canonical Hermitian structure.

Now, let focus on the quantum gauge group. As we have mentioned before, Definition \ref{qgg} is the one presented in \cite{br} but at the level of differential calculus, and of course, this definition does not recreate the {\it classical} case: it is to large. 

In Durdevich's theory there have been some attempts to get a definition of the quantum gauge group, for example in \cite{micho5, micho8}. To accomplish the purpose of this paper, the definition of $\qGG$ presented in \cite{micho5} is not useful because it does not create an action on the space of qpc's. On the other hand, the formulation showed in \cite{micho8} is only for the special case $B:=C^\infty_\C(M)$, where $M$ is a compact smooth manifold, and for a special graded differential $\ast$--algebra on $\G$: the minimal admissible calculus. This is why we decided to use Definition \ref{qgg}, despite the fact that it does not recreate the {\it classical} case.

One possible option to recover the {\it classical} case is to define $\qGG$ as the group of all graded differential $\ast$--algebra isomorphisms $$\F:\Omega^\bullet(P)\longrightarrow \Omega^\bullet(P)$$ that satisfy equation (\ref{4.f4}). However, depending on the qpb, the quantum gauge group would not have enough elements. This is a problem, for example, when we talk about Yang--Mills theory in non--commutative geometry, since the orbit of  Yang--Mills qpc's could be trivial \cite{sald2, sald3, sald4}. From a physical point of view, this implies that there could be too many non--gauge--equivalent boson fields. This issue may also arise if we define $\qGG$ with one more condition than the ones presented in Definition \ref{qgg}. For example, by requiring that the elements of $\qGG$ commute with the differential of $\Omega^\bullet(P)$. To prevent this from happening, we have decided to define $\qGG$ in the most general way, use equation (\ref{4.f11}) for the action on qpc's and work with {\it ad hoc} subgroups of $\qGG$ in each situation. 

In Subsection $4.2$, we have mentioned that in the literature, for example \cite{ha}, the commonly accepted action of the ($0$--degree) quantum gauge group on qpc's is given by
\begin{equation*}
    \F^\times\omega= \f_\F \ast \omega \ast \f_\F^{-1}+\f_\F\ast (d \f_\F^{-1}),
\end{equation*}
where $\omega$ is considered a map from $H$ to $\Omega^1(P)$. Similarly, the curvature, which in these references is defined by equation (\ref{falsecur}), satisfies (\cite{ha})
\begin{equation}
    \label{wrongactioncur}
    r^{\F^\times\omega}= \f_\F \ast r^\omega \ast \f_\F^{-1}.
\end{equation}
 Nevertheless, equation (\ref{wrongactioncur}) is not well--defined in Durdevich's framework because, as we have checked in Section $2.2$, in Durdevich's formulation the curvature is defined from the quantum dual Lie algebra $\mathfrak{qg}^\#$ to $\Omega^2(P)$ as in the {\it dualization} of the {\it classical} case, and only for multiplicative qpc's the curvature can be defined as $r^\omega$ (\cite{micho2,stheve}). In contrast, the action given in equation (\ref{4.f11}) is always well--defined in Durdevich's formulation. It is worth remembering that equation (\ref{4.f11}) is simply the {\it dualization} of the {\it classical} action of the gauge group on the space of principal connections via the pull--back (see equation (\ref{gaugeconn})). 

On the other hand, although this work has been developed in the framework of non--commutative geometry, the quantum gauge group is a {\it classical} group. Therefore, an exciting research project would be to explore a way to define $\qGG$ as a quantum group, although there would be a coaction on the space of qpc's instead of an action.

The reader is invited to notice the remarkable geometric--dual similarity of this theory with differential geometry, particularly in equations (\ref{3.f8}), (\ref{3.f8.1}), (\ref{3.f10.3}), (\ref{3.f10.6}); Definitions \ref{lhs}, \ref{rhs}, \ref{qgg}; Proposition \ref{4.4} and Theorems \ref{hil}, \ref{hil1}, \ref{fgs}, \ref{4.3}. This similarity even permits the development of a {\it non--commutative geometrical} version of electromagnetic field theory on the Moyal--Weyl algebra as shown in \cite{sald5}. In that case, the {\it non–-commutative Maxwell equations} are no longer identically zero in the vacuum, meaning that qpc's can represent photon fields that generate electric and magnetic charges and currents, even in the vacuum.\\

{\bf Conflict of Interests:\\} 

This paper was created under non--financial conflict of interests.

\end{document}